\documentclass[reqno,12pt]{amsart} 
\usepackage{amsmath,amssymb,amsfonts,eucal}
\usepackage[all]{xy}
\usepackage{graphicx,psfrag}

\newtheorem{thm}{Theorem}[section]
\newtheorem{cor}[thm]{Corollary}
\newtheorem{prop}[thm]{Proposition}
\newtheorem{lem}[thm]{Lemma}

\theoremstyle{remark}
\newtheorem{prob}[thm]{Problem}

\newtheorem{exam}{Example}[section]

\newcommand \myexam[1]{\smallskip\begin{exam}[\emph{#1}]}
\newcommand \myprob[1]{\smallskip\begin{prob}[\emph{#1}]}

\numberwithin{equation}{section}
\numberwithin{figure}{section}
\numberwithin{table}{section}

\renewcommand \phi{\varphi}
\renewcommand \epsilon{\varepsilon}
\newcommand \eset{\emptyset}
\renewcommand \emptyset{\varnothing}
\newcommand \setm{\setminus}

\newcommand \inv{^{-1}}
\newcommand \id{\mathrm{id}}

\newcommand \Aut{\operatorname{Aut}}
\newcommand \Aff{\operatorname{PsAut}}
\newcommand \clos{\operatorname{clos}}
\newcommand \bcl{\operatorname{bcl}}

\newcommand \ind{{:}}

\newcommand \bgr[1]{\langle#1\rangle}
\newcommand \full{^{{}^{{}_{{}_\bullet}}}}

\newcommand \cB{\mathcal{B}}
\newcommand \cC{\mathcal{C}}
\newcommand \cF{\mathcal F}

\newcommand \bbZ{\mathbb{Z}}

\newcommand \fG{\mathfrak G}
\newcommand \fH{\mathfrak H}
\newcommand \fQ{\mathfrak Q}

\newcommand \fV{\mathfrak V}

\newcommand \te{\tilde e}
\newcommand \tf{\tilde f}
\newcommand \ts{\tilde s}
\newcommand \tA{\tilde A}
\newcommand \tB{\tilde B}
\newcommand \tC{\tilde C}
\newcommand \tD{\tilde D}
\newcommand \tH{\tilde H}
\newcommand \tK{\tilde K}
\newcommand \tP{\tilde P}
\newcommand \tQ{\tilde Q}
\newcommand \tR{\tilde R}
\newcommand \tS{\tilde S}
\newcommand \tT{\tilde T}
\newcommand \teta{\tilde\eta}
\newcommand \tDelta{\tilde\Delta}

\newcommand \bbeta{\bar\beta}

\newcommand \uniquegains{Theorem V.2.1}

\begin{document}


\title{Associativity in Multary Quasigroups:  \\
	The Way of Biased Expansions}
\author{Thomas Zaslavsky}
\address{Binghamton University,
	Binghamton, New York 13902-6000}

\date{Jan.--Mar., 2002; Aug., 2003; Rev.\ Apr.--May, 2004.  Version of \today.}

\begin{abstract}
 A \emph{biased expansion} of a graph is a kind of branched covering graph 
with additional structure related to combinatorial homotopy of circles.  
Some but not all biased expansions are constructed from groups 
(\emph{group expansions}); these include all biased expansions of complete 
graphs (assuming order at least four), which correspond to Dowling's 
lattices of a group and encode an iterated group operation.  A biased 
expansion of a circle with chords encodes a multary (polyadic, $n$-ary) 
quasigroup, the chords corresponding to factorizations, i.e., associative 
structure.  We show that any biased expansion of a 3-connected graph (of 
order at least four) is a group expansion, and that all 2-connected biased 
expansions are constructed by expanded edge amalgamation from group 
expansions and irreducible multary quasigroups.  If a 2-connected biased 
expansion covers every base edge at most three times, or if every 
four-node minor is a group expansion, then the whole biased expansion is a 
group expansion. In particular, if a multary quasigroup has a 
factorization graph that is 3-connected, if it has order $3$, or if every 
residual ternary quasigroup is an iterated group isotope, it is isotopic 
to an iterated group.  We mention applications to generalizing Dowling 
geometries and to transversal designs of high strength.
 \end{abstract}



\thanks{Research partially assisted by grant DMS-0070729 from the National Science Foundation.}

\maketitle
\pagestyle{myheadings}
\markboth{\sc Thomas Zaslavsky}{\sc Associativity in Multary Quasigroups}

\pagebreak[4]\tableofcontents


\section{Biased graphs and the associative law}\label{bgal}

A \emph{multary quasigroup} is a set $\fQ$ with a multary operation 
$(\cdots): \fQ^n \to \fQ$, where $n \geq 2$, such that the equation $x_0 = 
(x_1\cdots x_n)$ is uniquely solvable for any one variable $x_i$ given the 
values of the remaining variables.  Multary quasigroups were introduced by 
Belousov and Sandik \cite{BelSand}.
 A \emph{Dowling geometry of a group}, $Q_n(\fG)$, is a certain matroid of rank $n\geq 1$ 
associated with a group $\fG$; it was invented by Dowling \cite{CGL} and shown by Kahn and 
Kung \cite{VCG} to have a central role in matroid theory. These two structures are both 
equivalent to particular kinds of the same general object, something I call a \emph{biased 
expansion graph}.  Associativity in multary quasigroups, and quasigroup generalizations of 
Dowling geometries, both depend on and can be analyzed through the structure of biased 
expansions.


\subsection{The associative law}\label{assoc}

The customary view of the associative law is that it describes a relationship between two 
different ways of carrying out a binary operation on three arguments:
 $$
(xy)z = x(yz).
 $$
 We look at it differently: we regard associativity as a property of factorizability or 
reducibility of a multary product.  For instance, letting $(\cdots)$ denote a ternary or binary product, we think of ordinary associativity as the combination of
\begin{equation}	\label{lefta}
(xyz) = ((xy)z)
\end{equation}
and
\begin{equation}	\label{righta}
(xyz) = (x(yz)).
\end{equation}
 The two factorizations \eqref{lefta} and \eqref{righta} constitute associativity in the 
usual sense (if all the binary multiplications are the same), but to us this is a 
secondary phenomenon.  We are more interested in association as factorization.  Our 
approach is based on representing the product by a kind of graph called a biased 
expansion; it treats multary operations but only of a certain kind: the equation 
$(xy\cdots z) = t$ must be uniquely solvable for any variable given all the others.  Sets 
with such operations, which were introduced by Belousov and Sandik \cite{BelSand}, are 
known as \emph{multary} or \emph{polyadic} \emph{quasigroups}.\footnote{The most common term is ``$n$-ary'', where $n$ is left unspecified, which is unsatisfactory.  Since 
``polyadic'' has not become popular, I propose ``multary'' as a generic adjective.  One 
occasionally sees ``multiary'', but etymologically this seems less correct.  We should 
mention that ternary quasigroups had previously been defined by Rado \cite{Rado}.}  We are 
investigating associativity in multary quasigroups by means of biased expansion graphs.

Formally, a \emph{multary quasigroup} (of \emph{type} or \emph{arity} $k$) is a set $\fQ$ with a multary operation $f : \fQ^k \to \fQ$ such that the equation
$$
f(x_1,\ldots,x_k) = x_0
$$
is uniquely solvable for any one $x_i$ if the other $k$ $x$'s are given.  (We assume $2 \leq k < \infty$, except that infinite arity is allowed in Section \ref{group}.)  
Our generalized associativity is (consecutive) factorization of $f$:
\begin{equation}\label{E:factor}
f(x_1,\ldots,x_k) = g(x_1,\ldots,x_{i},h(x_{i+1},\ldots,x_j),x_{j+1},\ldots,x_k),
\end{equation}
where $1 < i+1 < j \leq k$ and $g$ and $h$ are multary quasigroups of suitable arity.  Such a factorization is sometimes called an \emph{$(i+1,j)$-reduction} of $f$.  
As an associativity property we call it \emph{reductive associativity}.  
One of our main results (Theorem \ref{T:3connq}) is that, if the graph of factorizations of $f$ is 3-connected, $f$ is essentially an \emph{iterated group}: a multary quasigroup whose operation has the form $f(x_1,\ldots,x_k) = x_1 x_2 \cdots x_k$ computed in a group.  (That was known in the case that the graph is complete; see Section \ref{group}.)
This is one of the two extreme kinds of multary quasigroup.  
At the other extreme is an \emph{irreducible} multary quasigroup, whose operation has no factorizations at all.  (Any binary quasigroup is irreducible.  It is an iterated group just when it is a group.)

By ``essentially'' we mean up to \emph{isotopy}.  Multary operations $f:\fQ^k\to \fQ$ and $g:\fQ_1^k\to \fQ_1$ are called \emph{isotopic} if there exist bijections $\alpha_0,\alpha_1,\ldots,\alpha_k : \fQ \to \fQ_1$ such that
$$
g(x_1^{\alpha_1},\ldots,x_k^{\alpha_k}) = f(x_1,\ldots,x_k)^{\alpha_0}.
$$
From our graph-theoretic standpoint we cannot distinguish isotopic 
operations.  Nor can we distinguish between operations that are circularly parastrophic, where \emph{circular parastrophy} means circular permutation of the $k+1$ variables, i.e., replacing the operation $a_0=f(a_1,a_2,\ldots,a_k)$ by any operation $g$ defined by $a_i=g(a_{i+1},a_{i+2},\ldots,a_k,a_0,\ldots,a_{i-1})$ or $a_i=g(a_{i-1},a_{i-2},\ldots,a_0,$ $a_k,\ldots,a_{i+1})$ for some $i$, where the subscripts are taken modulo $k+1$.  (This must not be confused with arbitrary permutation of the operands, which is \emph{general parastrophy}.)
Since there is no name for the combination of isotopy and circular parastrophy, we shall call it \emph{circular allotopy}: multary operations are \emph{circularly allotopic} if one can be obtained from the other by a combination of isotopy and circular parastrophy.  We reiterate: our method does not distinguish multary quasigroups that are circularly allotopic.

Belousov treated Equation \eqref{E:factor} as representing a binary operation on functions, written $f = g +^{i+1} h$.  The resulting algebra led to many theorems.  A simple example is the property Belousov called \emph{$(i,j)$-associativity}:
\begin{equation}\label{E:ijassoc}
g +^i h = g' +^j h' .
\end{equation}
 Belousov solved this equation in \cite{BIAQ} by axiomatizing the algebra 
of multary quasigroups with composition operations $+^i$.  (See Corollary 
\ref{C:qfactor}.)  However, he left outstanding an important question:  
whether, whenever the factorization graph of a multary quasigroup is 
3-connected, the operation is isotopic to an iterated group 
operation.\footnote{Dudek \cite{Dpc} has heard that Belousov conjectured 
this to be true, but I have not been able to confirm that statement.}
 We prove this, as well as reproducing the solution of 
\eqref{E:ijassoc}, as corollaries of our structure theorem for biased 
expansions.

Another problem left incomplete by Belousov is that of small multary 
quasigroups.  If the order $|\fQ|$ is very small, $\fQ$ may have no room 
to fail to factor, i.e., to be other than an iterated group isotope.  For 
example, it is well known that a binary quasigroup of order at most 4 is a 
group isotope.
 According to information provided by Dudek \cite{Dpc}, Belousov and 
collaborator(s) proved that $\fQ$ of any arity is isotopic to an iterated 
group when $|\fQ| = 2$ and stated that the same holds true when $|\fQ|=3$ 
but that their proof of the latter was too long to be 
published.\footnote{I have not been able to locate any publication 
of these results.}
 These facts are a simple consequence of our second main criterion for a 
multary quasigroup to be an iterated group isotope.  A multary quasigroup 
obtained by fixing the values of some set of independent variables in 
$\fQ$ is called a \emph{residual} multary quasigroup of $\fQ$.  We show in 
Section \ref{4minors}, as interpreted in Theorem \ref{T:qinduced3}, 
that $\fQ$ is isotopic to an iterated group if its arity is at least three 
and every residual ternary quasigroup is an iterated group isotope.

The traditional multary generalization of associativity, due to D\"ornte 
\cite{Dornte} and extensively studied (see, e.g., 
\cite{Post,Gluskin,Hosszu,Dudek}), which is a stronger form of 
$(i,j)$-associativity, is a special case of our general picture.  A 
$k$-ary operation $f : \fQ^k \to \fQ$ is called \emph{associative} if it 
satisfies all the $k$ identities
 \begin{equation}	\label{E:assoc}
\hat f(x_1,\ldots,x_{2k-1}) = f(x_1,\ldots,x_{i},f(x_{i+1},\ldots,x_{i+k}),x_{i+k+1},\ldots,x_{2k-1})
\end{equation}
for $i = 1,\ldots,k$, where $\hat f$ is defined by any one of the identities.  (That is, \eqref{E:assoc} consists of $k-1$ identities and one definition.)  We might call this \emph{substitutive associativity} by way of contrast with reductive associativity.  A multary quasigroup with substitutive associativity for all $i$ is called a \emph{$k$-group} (or \emph{multary} or \emph{polyadic group}), a $2$-group being an ordinary group.  Evidently, $\hat f$ is an example of a multary quasigroup operation that is reducible in a multiplicity of ways.  
By our general theorem just mentioned, $f$ is isotopic to an iterated group.  
Hossz\'u \cite{Hosszu} and Gluskin \cite{Gluskin} have a far more complete result: an explicit construction of all multary groups in terms of groups, which was later much generalized.  
It should be possible to supplement our method so as to obtain their theorem, but we do not explore that question here.


\subsection{Introduction to expansions}\label{expan}

Biased expansions of circles are cryptomorphic to circular allotopy classes of multary quasigroups.  But first we must define biased expansions.

A \emph{biased graph} $\Omega = (\|\Omega\|, \cB)$ consists of a graph $\|\Omega\|$, which may be finite or infinite, and a \emph{linear class} $\cB$ of circles (circuits, cycles) of $\|\Omega\|$, meaning that in each theta subgraph the number of circles that belong to $\cB$ is different from two.  
The circles in $\cB$ are called the \emph{balanced circles} of $\Omega$.  

A \emph{biased expansion} of a graph $\Delta$ is a biased graph $\Omega$ together with a \emph{projection} mapping $p: \|\Omega\| \to \Delta$ that is surjective, is the identity on nodes, maps no balanced digon to a single edge, and has the \emph{circle lifting property}: for each circle $C = e_1 e_2 \cdots e_l$ in $\Delta$ and each $\tilde e_1 \in p\inv(e_1)$, \ldots, $\tilde e_{l-1} \in p\inv(e_{l-1})$, there is a unique $\tilde e_l \in p\inv(e_l)$ for which $\tilde e_1 \tilde e_2 \cdots \tilde e_l$ is balanced.  
In addition, no edge fiber $p\inv(e)$ may contain a balanced digon; but this is implied by the other properties if $e$ is not an isthmus.
One can think of $\Omega$ as a kind of branched covering of $\Delta$.  
We write $\Omega\downarrow_p\Delta$ to mean that $\Omega$ is a biased expansion of $\Delta$ with projection $p$; though usually we omit $p$ from the notation.
We call $\Omega$ a \emph{regular} or \emph{$\gamma$-fold} biased expansion if each $p\inv(e)$ has the same cardinality $\gamma$; then $\gamma$ is the \emph{multiplicity} of the expansion.  $\gamma\cdot\Delta$ denotes a $\gamma$-fold biased expansion of $\Delta$.  
Clearly, a biased expansion of an inseparable graph must be regular.  
A biased expansion is \emph{trivial} if it is regular with multiplicity $1$. 
In defining a biased expansion of a circle it is not necessary to require $\Omega$ to be a biased graph; that follows from the rest of the definition because a theta graph exists only by containing a digon.

A simple kind of biased expansion is a \emph{group expansion} \cite[Example I.6.7]{BG}.  
The expansion of a graph $\Delta$ by a group $\fG$, in brief the \emph{$\fG$-expansion} of $\Delta$, is the gain graph, denoted by $\fG\Delta$, whose node set is $N(\Delta)$ and whose edge set is $\fG \times E(\Delta)$, the endpoints of an edge $ge$ (this is shorthand for $(g,e)$) being the same as those of $e$.  
The projection $p : \fG\Delta \to \Delta$ maps $ge$ to $e$.  
We associate with $ge$ the group element $g$, called the \emph{gain} of $ge$; in order to define gains in a technically correct manner we orient $\Delta$ arbitrarily and orient $ge$ similarly to $e$, so the gain of $ge$ in the chosen direction is $g$ and in the opposite direction $g\inv$.  
A circle in $\fG\Delta$ is balanced if the product of the gains of its edges, taken in a consistent direction around the circle, equals $1$, the group identity.  This defines a biased graph, which we write $\bgr{\fG\Delta}$.  
If $\Delta$ is simple with $n$ nodes, then $\fG\Delta$ is contained in $\fG K_n$; thus group expansions of complete graphs are basic.  

A very different kind of biased expansion is the expansion of a circle $C_{k+1}$ of length $k+1$ by a $k$-ary quasigroup $\fQ$.  
In the \emph{quasigroup expansion} $\bgr{\fQ C_{k+1}}$, the nodes are $v_0, v_1, \ldots, v_k$.  There is an edge $ae_{i-1,i}$ for every $a\in \fQ$ and $i=1,2,\ldots,k$ as well as an edge $ae_{0k}$.  The balanced circles are the circles $\{a_1e_{12}, a_2e_{12}, \ldots, a_ke_{k-1,k}, a_0e_{0k} \}$ such that $(a_1a_2\cdots a_k) = a_0$ in $\fQ$.  
A quasigroup expansion need not be contained in an expansion of a complete graph; see Theorem \ref{L:q}.

Conversely, from a biased expansion $\gamma\cdot C_{k+1}$ we obtain a circular allotopy class of $k$-ary quasigroups.  Let $C_{k+1}$ have nodes $v_0,v_1,\ldots,v_k$ and edges $e_{01},e_{12},\ldots,e_{k-1,k},e_{0k}$.  
Set $E_{ij} = E(\gamma\cdot C_{k+1})\ind\{v_i,v_j\}$ and fix bijections $\beta_i : \fQ \to E_{i-1,i}$ for $i=1,\ldots,k$ and $\beta_0 : \fQ \to E_{0k}$, where $\fQ$ is a set that will be the set of elements of the $k$-ary quasigroup.  
The multary operation is $(x_1\cdots x_k) = \beta_0\inv(\te)$ where $\te$ is the unique edge in $E_{0k}$ that forms a balanced circle with $\beta_1(x_1),\ldots,\beta_k(x_k)$.  
The arbitrary choice of the bijections is what makes $\fQ$ defined only up to isotopy.  The arbitrariness of the distinguished edge $e_{0k}$ and the direction of reading the circle is what leaves $\fQ$ well defined only up to circular parastrophy.  Thus we have the first two parts of Proposition \ref{P:qgroup}.  
The third part is proved at Theorem \ref{L:q}.  We say $\bgr{\fQ C_{k+1}}$ \emph{extends} to $e_{ij}$ if there is a biased expansion $\Omega\downarrow(\Delta\cup\{e_{ij}\})$ such that $p\inv(C_{k+1}) = \bgr{\fQ C_{k+1}}$.  (Section \ref{extension} has a fuller treatment.)

\begin{prop}\label{P:qgroup}
A $k$-ary quasigroup expansion $\bgr{\fQ C_{k+1}}$ is a biased expansion of $C_{k+1}$.  Conversely, every biased expansion of $C_{k+1}$ has the form $\bgr{\fQ C_{k+1}}$ for a $k$-ary quasigroup $\fQ$. 

Furthermore, two $k$-ary quasigroup expansions $\bgr{\fQ_1 C_{k+1}}$ and $\bgr{\fQ_2 C_{k+1}}$ are isomorphic if and only if $\fQ_1$ and $\fQ_2$ are circularly allotopic.

Moreover, $\bgr{\fQ C_{k+1}}$ extends to a chord $e_{ij}$ of $C_{k+1}$ if and only if the operation $f$ of $\fQ$ factors as in \eqref{E:factor}.
\end{prop}

Taking $k=2$, we see that biased expansions of a triangle are (as Dowling knew implicitly; see \cite[pp.\ 78--79]{CGL}) graph-theoretic realizations of circular allotopy classes of binary quasigroups and, therefore, of Latin squares and 3-nets.  The quadrangle criterion of Latin squares \cite[Theorem 1.2.1(2)]{Latin} tells us when a binary quasigroup is isotopic to a group; its translation into the language of expansions is the following, applicable to any biased expansion of a triangle:
\begin{quote}
\emph{Quadrangle Criterion.}  
For any twelve distinct edges $e^{\alpha,\beta}_{12}$, $e^{\alpha,\beta}_{23}$, $e^{\alpha,\beta}_{13}$, where $e^{\alpha, \beta}_{ij} \in p\inv(e_{ij})$ for $\alpha, \beta=1,2$ and $ij = 12, 23, 13$, if seven of the triangles of the form $e^{\alpha,\delta}_{12} e^{\beta,\delta}_{23} e^{\alpha,\beta}_{13}$ are balanced, then so is the eighth.
\end{quote}

\begin{prop} \label{PD:bgexpan}
A biased expansion of $K_3$ has the form $\bgr{\fG K_3}$ for some group $\fG$ if and only if it satisfies the quadrangle criterion.
\hfill\qedsymbol
\end{prop}

A $\gamma$-fold biased expansion of $C_{k+1}$, where $k\geq 3$, giving us a $k$-ary quasigroup according to Proposition \ref{P:qgroup}, can be interpreted as the $k$-dimensional generalization of a Latin square that is called a \emph{permutation hypercube} \cite[p.\ 181]{Latin} or \emph{Latin hypercube}, defined up to isotopy and parastrophy.  
As far as I know, no analog of the quadrangle criterion has been formulated for such objects and therefore for biased expansions of larger graphs such as $C_{k+1}$ for $k \geq 3$.

\bigskip

In a sense indicated by Proposition \ref{P:qgroup}, biased expansion graphs are a graphical generalization of groups and multary quasigroups.
That they truly are a generalization is shown by the fact that a biased expansion need not have a Hamiltonian circle.  If the base graph is Hamiltonian, then we have a multary operation---which in general depends on the choice of Hamiltonian circle---from Proposition \ref{P:qgroup}, but if it is not, then we have something that, from an algebraic standpoint, is more complicated; it might be thought of as a combinatorial complex of multary quasigroups.

If biased expansion graphs generalize groups, it is natural to ask how far a given biased expansion is from being a group.  
 Biased expansions of complete graphs are group expansions, essentially 
because $K_4$, as the base of a biased expansion, encodes the associative 
law (see Section \ref{group}).
 Thus, a more precise version of the question is:  How far is the base 
graph from being complete?  But this is still not quite right, because it 
might be possible to extend the expansion to new edges between nonadjacent 
nodes.  If the expansion extends to a complete graph, then it is contained 
in a group expansion so it itself is a group expansion and any 
corresponding multary quasigroup is an iterated group isotope.  In 
general, there is always a maximal extension of the given biased expansion 
graph, that has the most pairs of adjacent nodes (see Section 
\ref{extension}); it is of this extended expansion that we should ask 
the refined question, and indeed it makes sense to think of the number of 
nonadjacent node pairs in the base graph as a measure of how much a biased 
expansion fails to represent a group.

 It is perhaps noteworthy that nongroup expansions of large incomplete 
graphs, and in particular irreducible multary quasigroups, exist at all.  
However, all examples are 2-separable, for, as we prove in Sections 
\ref{extension} and \ref{group}, every biased expansion of a 
3-connected simple graph having more than three nodes is a group 
expansion.  From this and other work we can deduce the complete structure 
of a biased expansion (Section \ref{bx}) and answer the question raised 
in \cite[Example III.3.8]{BG} of exactly which graphs have a nongroup 
biased expansion (Corollary \ref{C:basenongroup}).


\subsection{Dowling geometries}\label{dowl}

Biased expansion graphs were inspired by Dowling's matroids of a group---though no matroids were used in the preparation of this article.  
One way to construct the rank-$n$ Dowling matroid (or ``geometry'') $Q_n(\fG)$ of a group $\fG$ is to take the group expansion $\fG K_n$, adjoin a half edge to each node, and take the bias matroid \cite[Section II.2]{BG}.  I sketch this construction (from \cite[Example III.5.7 or Part V]{BG}) to suggest how biased expansions of simple graphs, especially those that are maximal on the given node set, are a natural generalization of Dowling geometries.  The unexplained terms can be found in \cite[Part III or V]{BG}.

Given any biased expansion $\Omega\downarrow\Delta$, one can add an unbalanced loop at each node and take the bias matroid; call this matroid $G\full(\Omega)$.  The operator $G\full$, applied to maximal biased expansions, generalizes the construction of Dowling geometries.  The Dowling geometries are the only examples derived from groups because the only maximal group expansions are those of the complete graphs $K_n$.  
Given that nongroup biased expansions exist, one asks what other matroids can be obtained from maximal biased expansions; they are natural candidates for generalized Dowling geometries.  That question motivated this work. 
We will see (via Theorem \ref{L:q}) that an expansion $\fQ C_{n+1}$ of $C_{n+1}$ by an irreducible $n$-ary quasigroup $\fQ$ is maximal, since it can have no chords; thus these are part of the answer.  The question is completely answered if we can classify all maximal biased expansions.  That is our Theorem \ref{T:structuremax}---the solution of a problem that had puzzled me since 1976.


\subsection{Transversal designs}\label{td}

A final way to look upon a multary quasigroup, or a biased expansion of a circle, is as a kind of transversal design.  A \emph{transversal $t$-design} consists of a set of points partitioned into $l$ \emph{point classes} $L_i$ (usually called ``groups'', but they have nothing to do with algebra) of $k$ points each, and a class of \emph{blocks}, which are subsets of points satisfying
\bigskip
\begin{enumerate}
\item[(TD1)] no two points in a class are contained in a common block, and
\bigskip
\item[(TD2)] any $t$ points, no two in a class together, are contained in exactly $\lambda$ common blocks,
\end{enumerate}
\bigskip
where $\lambda$, the \emph{index}, is a fixed positive integer.  As $t$ is called the \emph{strength}, we refer to an $(l-1)$-design as having \emph{high strength}.  
A $k$-fold biased expansion $\Omega$ of $C_l$ is equivalent to a transversal design $T$ of high strength with $\lambda=1$.  The points of $T$ are the edges of $\Omega$, the class $L_i$ consists of all edges with endpoints $v_{i-1}$ and $v_i$, and the blocks are the balanced circles.

A group expansion thus generates a design based on the group.  
The construction of the design is easy to describe directly.  
The classes are copies of the group and a block is any set $\{ x_{1}, x_{2}, \ldots, x_{l} \}$, composed of one element of each class $L_{j}$, such that $x_1 x_2 \cdots x_l = 1$.
Here we need to assume the classes are ordered; we let the first class be $L_1$, the second $L_2$, etc.  
The analog of factorization is consecutive composition: supposing $T'$ and $T''$ are two transversal designs of high strength with index 1, we form their \emph{$i$-composition} $T$ by identifying $L'_i$ with $L''_1$ and then defining the classes of $T$ to be those of $T'$ and $T''$, with the exception of $L'_i$ ($=L''_1$), and the blocks to be the sets of the form $B'\oplus B''$ where $B'$ and $B''$ are overlapping blocks of $T'$ and $T''$.  The inverse operation to $i$-composition is \emph{$i$-decomposition}.  The analog of the factorization graph is defined by the existence of $i$-decompositions of $T$.  We have the theorems, for instance, that if this graph is 3-connected and the number of classes is at least four, then the design is derived from a group, and that if every four-class transversal design induced by $T$ and including the class $L_1$ is derived from a group, then $T$ is so derived.  (Precise statements can be obtained by translating results of Section \ref{qgroup}.)

One is inspired by the design interpretation to wonder about generalizing to larger values of $\lambda$.  In terms of biased expansions:

\begin{prob}\label{Pr:lambda}
 Suppose for each lift $\tP$ of $C\setm e$ there are exactly $\lambda$ 
edges $\te$ that make a balanced circle, where $\lambda>1$.  Do our 
theorems generalize to this situation?
 \end{prob}

The operational view of this generalization is that we have a $\lambda$-valued $n$-ary operation, each of whose $n$ inverse operations is also $\lambda$-valued.
One has to modify the definition of biased expansion:  $\cB$ need no longer be a linear class; instead, only each separate lift of a base theta subgraph would be subject to the linearity condition that its number of balanced circles be different from two.  
The value of $\lambda$ cannot be a constant, because, in operational language, the composition of $\lambda_1$-valued and $\lambda_2$-valued operations is $\lambda_1\lambda_2$-valued.


\subsection{Overview}\label{ov}

A summary of our main results:
\medskip

\emph{Biased expansion graphs.}
  \begin{itemize}
  \item A 3-connected expansion is a group expansion (Theorem \ref{T:3connbg}).
  \item Edge amalgamation and edge sum for 2-separable expansions (Theorem \ref{T:twist}).
  \item Decomposition into group and multary quasigroup expansions (Theorems \ref{T:structuremax}, \ref{T:structureexp}).
  \item Characterization of base graphs having nongroup biased expansions (Corollary \ref{C:basenongroup}).
  \item Uniqueness and structure of maximal biased expansions (Theorems \ref{TX:max}, \ref{T:structuremax}, Corollary \ref{C:basemax}).
  \item An expansion of multiplicity at most 3 is a group expansion (Theorem \ref{T:thin}).
  \item A 2-connected biased expansion with at least four nodes is a group expansion if every minor of order four is gainable (Theorem \ref{T:gminor4}).
  \end{itemize}
\medskip

\emph{Multary quasigroups.}
  \begin{itemize}
  \item The factorization graph corresponds to maximal extension (Theorem \ref{L:q}).
  \item A nonbinary multary quasigroup whose factorization graph is 3-connected is an iterated group isotope (Theorem \ref{T:3connq}).
  \item Characterization of possible factorization graphs (Theorem \ref{L:qtheta}).
  \item A multary quasigroup of which every residual ternary quasigroup is an iterated group isotope is itself isotopic to an iterated group (Theorem \ref{T:qinduced3}).
  \end{itemize}


\section{Preliminary remarks}\label{prelim}

Here we collect a few old and new definitions and some elementary observations about expansions.


\subsection{Basic concepts}\label{concepts}

Formally, a graph $\Gamma$ is a pair $(N,E)$ consisting of a node set and an edge set.  
Edges are \emph{links} (two distinct endpoints) or \emph{loops} (two coincident endpoints).  
A graph all of whose edges are links is a \emph{link graph}; a link graph without parallel edges is \emph{simple}.  
A \emph{circle}, also known as a polygon, circuit, or cycle, is the graph or edge set of a simple closed path.  A \emph{theta graph} is the union of three internally disjoint paths with the same two endpoints; the paths are the \emph{constituent paths} of the theta graph.  
The sum (symmetric difference) of sets is written $S \oplus T$.  This applies in particular to circles, regarded as edge sets, and especially to circles whose union is a theta graph.
A graph is \emph{2-connected} or \emph{inseparable} if any two edges lie in a common circle.  A \emph{block} of a graph is a maximal inseparable subgraph.  A 3-connected graph is assumed to be inseparable.
An \emph{induced subgraph} is a subgraph $\Gamma\ind X = (X,E\ind X)$, where $X \subseteq N$, whose edges are all those of $\Gamma$ with both endpoints in $X$.

Suppose $\Gamma_1$ and $\Gamma_2$ are two graphs that have in common a link $e$.  An \emph{(edge) amalgamation} of $\Gamma_1$ and $\Gamma_2$ along $e$, written $\Gamma_1 \cup_e \Gamma_2$, is a graph obtained by identifying the two copies of $e$.  (It is not usually unique, since the copies can be identified in two ways.)  An \emph{edge sum} (or \emph{2-sum}), written $\Gamma_1 \oplus_e \Gamma_2$, is $\big(\Gamma_1 \cup_e \Gamma_2\big) \setm e$.

Suppose $\Gamma$ is a graph and $\Xi$ a subgraph.  A \emph{bridge of} $\Xi$ in $\Gamma$ is a maximal subgraph $B$ of $\Gamma$ with the properties that $E(B)\cap E(\Xi) = \eset$, $B\not\subseteq \Xi$, and any node common to $B$ and an edge not in $B$ lies in $\Xi$.  (See Tutte \cite[Section I.8]{Tbook}.)  

Suppose $\Delta$ is a graph and $\gamma$ is a nonzero cardinal number, possibly infinite: then, by $\gamma\Delta$ we mean $\Delta$ with every edge replaced by $\gamma$ copies of itself.  Thus the underlying graph of a regular biased expansion $\gamma\cdot\Delta$ is $\gamma\Delta$; note the importance of the dot in the notation.

Biased graphs were defined in the introduction.  Some additional notions:  
A subgraph or edge set in a biased graph $\Omega$ is called \emph{balanced} if every circle in it is balanced.  
A balanced biased or gain graph should be thought of as like an ordinary graph and the bias, i.e., the choice of balanced circles, as a kind of skewing; so the less balanced, the more biased.

A \emph{gain graph} $\Phi = (\|\Phi\|,\phi)$ with \emph{gain group} $\fG$ is a graph $\|\Phi\|$ together with an orientable \emph{gain function} $\phi : E(\Phi) \to \fG$; that is, $\phi$ is defined on oriented edges and, letting $e\inv$ denote $e$ with the opposite orientation, $\phi(e\inv) = \phi(e)\inv$.  A group expansion, obviously, is a gain graph.  A circle in $\Phi$ is called \emph{balanced} if the product of the gains of its edges is $1$, the group identity; thus $\Phi$ produces a biased graph $\bgr{\Phi}$.  
\emph{Switching} a gain graph $\Phi$ by a \emph{switching function} $\eta: N \to \fG$ means replacing $\phi$ by a new gain map, $\phi^\eta$, defined by $\phi^\eta(e) = \eta(v)\inv \phi(e) \eta(w)$ if $e$ is oriented from endpoint $v$ to endpoint $w$.  Switching gains does not change bias: $\bgr{\Phi}=\bgr{\Phi^\eta}$.  
Not all biased graphs are obtainable from gains.  An expansion of $C_{k+1}$ by a multary quasigroup that is not isotopic to an iterated group is one example; we shall see others in Section \ref{amalg}.  
Biased graphs and gain graphs are from \cite[Part I]{BG}.  

A \emph{minor} of a graph, biased graph, or gain graph is a subgraph or 
contraction of a subgraph. 
 Since contraction of biased and gain graphs is complicated and 
plays a minor role in this article, we omit the definitions, referring the 
reader to \cite[Sections I.2 and I.5]{BG}.

We shall have use for a theorem about chains of paths.  If $A, B \subseteq \Delta$, an \emph{$AB$-path} is a path with one endpoint in $A$, the other in $B$, and otherwise disjoint from $A \cup B$.

\begin{lem}[Path Lemma]\label{C:epath}
In a 2-connected graph $\Delta$ let $A$ and $B$ be disjoint paths.  If $P_0$ and $P$ are two $AB$-paths, then there exist $AB$-paths $P_1, \ldots, P_k = P$ such that each $P_{i-1} \cup P_i \cup A \cup B$ contains exactly one circle.
\end{lem}

\begin{proof}
This lemma can be deduced from Tutte's Path Theorem {\cite[Theorem 4.34]{Tlect}}, but we give a direct proof suggested by Marcin Mazur \cite{Mazur}.  
The result is trivial if $P_0$ and $P$ are internally disjoint.  
Otherwise, let the nodes of $P_0$, in order from $A$ to $B$, be $v_0, v_1, \ldots, v_l$; define $P_0(v_j)$ to be that portion of $P_0$ from $v_{j+1}$ to $v_{l-1}$; and make a similar definition for $P$. 
Let $x_0$ be the first node of $P$ when traced from $A$ to $B$.  Let $x_1$ be the first node of $P(x_0)$ that lies in $P_0(v_0)$, $x_2$ the first node of $P(x_1)$ that lies in $P_0(x_1)$, and in general $x_i$ the first node of $P(x_{i-1})$ that lies in $P_0(x_{i-1})$. 
Define $k-1$ as the last value of $i$ for which an $x_i$ exists. 
For $0<i<k$, $P_i$ is obtained by tracing $P$ from $x_0$ to $x_i$ and then $P_0$ from $x_i$ to $v_l$; and $P_k$ is $P$.  
Then $P_{i-1} \cup P_i$ contains the circle formed by the segments of $P_0$ and $P$ from $x_{i-1}$ to $x_i$, and no other circle; except that the unique circle in $P_0\cup P_1$ consists of $P_0$ and $P_1$ up to $x_1$ along with $A$ from $v_0$ to $x_0$, and the unique circle in $P_{k-1} \cup P_k$ consists of $P_0$ and $P$ from $x_{k-1}$ to their endpoints in $B$ together with $B$ between those endpoints.
\end{proof}

A \emph{homomorphism} (synonym: \emph{mapping}) of graphs is an incidence-preserving mapping of node and edge sets.
A \emph{homomorphism of biased graphs} is a homomorphism of the underlying graphs that preserves balance, but not necessarily imbalance, of edge sets.

In a biased graph there is a kind of closure called the \emph{balance-closure} (not ``balanced closure''), defined for any edge set by 
$$
\bcl S = S \cup \{ e \notin S : \text{ there is a balanced circle $C \ni e$ such that } C\setm e \subseteq S \} .
$$
This is not an abstract closure operator, nor is it true that $\bcl S$ must be balanced; the essential property of balance-closure is

\begin{lem}[{\cite[Proposition I.3.1]{BG}}]\label{L:bcl}
For $S \subseteq E(\Omega)$, $\bcl S$ is balanced if and only if $S$ is balanced.
\end{lem}


\subsection{Basics of expansions}\label{basicexp}

Elementary facts about expansions let us confine our attention to simple, inseparable base graphs.  First, it is clear that a biased expansion of a graph $\Delta$ is the union of arbitrary biased expansions of the blocks of $\Delta$.  Second, biased expansion of a loop is uninteresting.  Third, suppose $e$ and $f$ are parallel links in $\Delta$.  In a biased expansion $\Omega$ of $\Delta$, there is a unique bijection between $p\inv(e)$ and $p\inv(f)$ such that, if $\te$ and $\tf$ correspond, then $\te\tf$ is balanced and, for any set $\tP \subseteq E(\Omega)$ that contains neither $\te$ nor $\tf$, $\tP \cup \{\te\}$ is a balanced circle if and only if $\tP \cup \{\tf\}$ is.  Thus, $\Omega$ is completely determined by $\Omega \setm p\inv(f)$.  
Moreover, any gains $\phi$ for $\Omega$ are completely determined by the gains on $\Omega \setm p\inv(f)$ by the equation $\phi(\tf)=\phi(\te)$.  (See Example \ref{XX:parallel}.)

Moreover, any gains $\phi$ for $\Omega$ are completely determined by the gains on $\Omega \setm p\inv(f)$ by the equation $\phi(\tf)=\phi(\te)$.
Still further, if a biased expansion graph $\Omega$ has gains in a group $\fG$, then it is a group expansion by a subgroup of $\fG$ \cite[\uniquegains(a)]{BG}.  Thus the nongroup biased expansions are the same as the nongainable biased expansions.

A \emph{homomorphism} (or \emph{mapping}) $\Omega \to \Omega'$ of biased expansions is a biased-graph homomorphism along with a homomorphism of base graphs such that the two mappings commute with projection.  We shall have occasion to use only homomorphisms that are injective.

The \emph{restriction of $\Omega$ to $\Delta'\subseteq\Delta$}, written $\Omega\big|_{\Delta'}$, is the subgraph $p\inv(\Delta')$ with the bias and projection mapping inherited from $\Omega$.  (This should not be confused with restricting $\Omega$ to an arbitrary subgraph of itself; $\Omega\big|_{\Delta'}$ is one such restriction, but not all restrictions are of that form.)

A basic property of expansions is the existence of balanced copies of $\Delta$ or any subgraph, extending any balanced subgraph of the expansion.  A \emph{lift} of an edge set $S \subseteq E(\Delta)$ is a subset $\tS \subseteq E(\Omega)$ for which $p\big|_{\tS}$ is a bijection onto $S$.  We shall always mean by $\tS$ a lift of $S$.

\begin{lem}
\label{L:balsub}
Let $\Omega$ be a biased expansion of a graph $\Delta$.  
Given any subsets $A \subseteq B \subseteq E(\Delta)$ and a balanced lift $\tA$, there is a balanced lift $\tB$ that contains $\tA$.
\end{lem}

\begin{proof}
Extend $A$ to a maximal subgraph $S$ of $B$ that has no additional circles besides those in $A$.  Take any lift $\tS\supseteq\tA$; it is balanced because $\tA$ is balanced.  Then $\bcl \tS$ projects to $\clos S = \clos B$, where $\clos$ is the ordinary graphic matroid closure
$$
\clos S = S \cup \{ e \notin S : \text{ there is a circle $C \ni e$ such that } C\setm e \subseteq S \} .
$$
Thus, $\bcl\tS$ is balanced by Lemma \ref{L:bcl}, and it contains a lift of $B$.  Take $\tB = p\inv(B) \cap \bcl\tS$.
\end{proof}

One can apply the lemma, for example, when $A$ is a forest, since any lift of a forest is balanced.  It is also the basis for an alternative definition of biased expansions; see \cite[Part V]{BG}.

The nongroup biased expansions are the same as the nongainable biased expansions, because if a biased expansion graph $\Omega$ has gains in a group $\fH$, then it is a group expansion by a subgroup $\fG$ of $\fH$.  Moreover, $\fG$ is unique up to isomorphism \cite[\uniquegains]{BG}.  
Furthermore, if $\Omega\downarrow\Delta$ is a group expansion by $\fG$, one can choose the gain mapping $\phi: E(\Omega) \to \fG$ so that $\phi\inv(1)$ is any desired balanced lift of $E(\Delta)$ (a consequence of \cite[Lemma I.5.3]{BG}), and then $\phi$ is determined up to automorphisms of $\fG$ \cite[\uniquegains]{BG}.
(I interpret the choosability of $\phi\inv(1)$ to mean that selecting a balanced lift of $\Delta$ is the expansion-graph analog of isotoping a quasigroup to a loop.)


\subsection{Expansion minors}\label{minor}

Certain minors of a biased expansion are themselves expansions.  An instance is a restriction of a biased expansion to a subset of the fibers, that is, $\Omega\big|_{\Delta'}$ where $\Delta'$ is any subgraph of $\Delta$; analogously, the restriction of $\fG\Delta$ is an expansion $\fG\Delta'$.  
Similarly, a contraction of a group or biased expansion by a balanced edge set is again a group or biased expansion save for possibly having extra balanced or unbalanced loops; for instance, if $\tS$ is a balanced edge set in a biased expansion $\Omega$ of $\Delta$, then $\Omega/\tS$ without loops is a biased expansion of $\Delta/p(\tS)$ without loops.  We want a notion that combines both of these kinds of minors, that of an `expansion minor'.

Let $\Omega$ be a biased expansion of a graph $\Delta$.  An \emph{expansion minor} of $\Omega$ is any minor $\Omega'$ of $\Omega$ (without loose or half edges) whose edge set is a union of fibers $p\inv(e)$ of $\Omega$; that is, $E(\Omega') = p\inv(S)$ for some $S\subseteq E$.  An expansion minor of a group expansion is similar.

\begin{prop}\label{P:expminor}
Let $\Delta$ be a graph.
\begin{enumerate}
\item[(a)]  An expansion minor $\Omega'$ of a biased expansion $\Omega$ of $\Delta$ is a biased expansion of a minor $\Delta'$ of $\Delta$.  If $\Omega$ is regular of multiplicity $\gamma$, then so is $\Omega'$.  
\item[(b)]  An expansion minor of a group expansion $\fG\Delta$ is a group expansion $\fG\Delta'$ of a minor $\Delta'$ of $\Delta$, and conversely.
\item[(c)]  An expansion minor of an expansion minor of $\Omega$ or $\fG\Delta$ is an expansion minor of $\Omega$ or $\fG\Delta$, respectively.
\end{enumerate}
\end{prop}

Part (b) is especially significant.  It says that we can tell something about the gainability of a biased expansion from its \emph{triangular expansion minors}, that is, expansion minors that are expansions of $K_3$.  We apply this idea in Section \ref{amalg}.

The proof is in the more precise description of expansion minors contained in two lemmas.  First we define a construction method for expansion minors.  

\begin{quote}
\emph{Construction XM.}  Given a biased expansion $\Omega$ of $\Delta$, take $S\subseteq E(\Delta)$, a weak partition $E(\Delta) =S\cup T\cup D$ (that is, $S$, $T$, and $D$ are pairwise disjoint sets whose union is $E(\Delta)$; some of them may be void), and a balanced lift $\tT$ of $T$ into $\Omega$.  From $[\Omega\setm p\inv(D)]/\tT$ delete all loops that belong to $p\inv(T)$, and delete an arbitrary subset of the isolated nodes (if there are any).  Call this $\Omega'$.
\end{quote}

\begin{lem}\label{L:bexpminor}
\begin{enumerate}
\item[(a)]  The biased graph $\Omega'$ of Construction XM is an expansion minor of $\Omega$, and every expansion minor of $\Omega$ is formed in this way.
\item[(b)]  $\Omega'$ is a biased expansion of a graph $\Delta'$ which is a minor of $\Delta$ formed from $\Delta\setm D/T$ by deleting some subset of its isolated nodes (if any).
\item[(c)]  $E(\Omega') = p\inv(S)$; the projection mapping $p' = p\big|_{p\inv(S)}$; and $(p')\inv(e) = p\inv(e)$ for each $e\in S = E(\Delta')$.  
\item[(d)]  If $\Omega$ is regular, then $\Omega'$ is regular with the same multiplicity.
\end{enumerate}
\end{lem}

\begin{proof}  We assume that the reader is acquainted with the definitions and notation of contraction and minors in \cite[Sections I.2 and I.5]{BG}.

(a)\quad  It is clear that $\Omega'$ is an expansion minor; the task is to prove the converse.  A minor of $\Omega$ is formed by contracting an edge set $A$, then deleting a subset of $A^c$.  (We ignore isolated nodes as a triviality.)  Let $A_0 = A\ind N_0(A)$ and $\tT = A\setm A_0$.  Some of the edges after contraction may be half or loose edges if $A_0 \neq \eset$.  The half edges come in entire fibers $p\inv(e)$, where $e$ joins $N_0(A)$ to its complement.  The loose edges come in fibers $p\inv(e)$ where $e\in E(\Delta)\ind N_0(A)$ but $e\not\in p(A_0)$, or in partial fibers $p\inv(e)\setm A_0$ where $e \in p(A_0)$.  In either case we may simply delete the entire fiber; at worst this leaves extra isolated nodes.  Thus we delete $p\inv(D_1)$ where $D_1 = E(\Delta)\setm E(\Delta\ind N_0(A)^c)$.  

This leaves us contracting only the balanced part $\tT$ of $A$; a process that results in no half or loose edges.  To get $\Omega'$ we must delete all the remaining edges in $p\inv(T)$, where $T = p(\tT)$; all of these are loops.  The remaining graph $\Omega''$ now meets the definition of an expansion minor of $\Omega$; it differs from $\Omega'$ only in that the latter may require deleting more edges, which must be whole fibers $p\inv(e)$ for $e\in D_2\subseteq E(\Delta)$.  Thus $D = D_1\cup D_2$ and $S = E(\Delta)\setm (T\cup D)$ in Construction XM.

(b)\quad  We have to prove that, for any circle $C$ in $\Delta' = (\Delta\setm D)/T$, edge $e\in C$, and lift $\tP$ of $C\setm e$ into $\Omega'$, there is a unique edge $\te \in (p')\inv(e)$ such that $\tP \cup \{\te\}$ is balanced. $C$ has the form $C_1\cap S$ where $C_1$ is a circle in $\Delta\setm D$.  Lift $C_1\setm S$ to $\tQ \subseteq \tT$.  Then $\tP \cup \tQ$ is a lift of $C_1\setm e$ into $\Omega$, for which there is a unique $\te\in p\inv(e)$ that makes $\tP\cup \tQ\cup \{\te\}$ balanced.  By the definition of contraction, for $\te\in p\inv(e) = (p')\inv(e)$, $\tP\cup \tQ\cup \{\te\}$ is balanced in $\Omega$ if and only if $\tP \cup \{\te\}$ is balanced in $\Omega'$.  This concludes the proof of (b).

(c) and (d) are obvious. 
\end{proof}

There are an analogous construction and lemma for group expansions.

\begin{quote}
\emph{Construction GXM.}  Given $\fG\Delta$, take $S$, $T$, $D$, and $\tT$ as in Construction XM and modify $[\fG\Delta\setm p\inv(D)]/\tT$ as in that construction to form $\Phi'$.
\end{quote}

\begin{lem}\label{L:gexpminor}
\begin{enumerate}
\item[(a)]  The gain graph $\Phi'$ of Construction GXM is an expansion minor of $\fG\Delta$, and every expansion minor of $\fG\Delta$ is formed in this way.
\item[(b)]  $\Phi' \cong \fG\Delta'$, where $\Delta'$ is as in Lemma \ref{L:bexpminor}(b).
\item[(c)]  For every minor $\Delta'$ of $\Delta$, $\fG\Delta'$ is an expansion minor of $\fG\Delta$.
\item[(d)] $\bgr{\Phi'}$ is an expansion minor of $\bgr{\fG\Delta}$, and every expansion minor of $\bgr{\fG\Delta}$ equals $\bgr{\fG\Delta'}$ for a minor $\Delta'$ of $\Delta$.  
\item[(e)]  Construction XM applied to $\bgr{\fG\Delta}$ yields $\bgr{\Phi'}$.  
\end{enumerate}
\end{lem}

\begin{proof}  (a) is proved as in Lemma \ref{L:bexpminor}.  (c) follows from (a) by taking $\tT = p\inv(T)\cap E(\{1\}\Delta)$ in the construction.  (e) is obvious from the constructions.  (d) follows from (e) and (b).

(b)  As a minor of $\fG\Delta$, $\Phi'$ has gains in $\fG$ \cite[Theorem I.5.4]{BG}.  We may assume by prior switching of $\Phi = \fG\Delta$ that $\phi\big|_{\tT} \equiv 1$.  Thus $\phi' = \phi\big|_{E(\Phi')}$, so $\phi'\big|_{(p')\inv(e)}$ is a bijection onto $\fG$.  It follows easily that $\Phi'\cong \fG\Delta'$.  (We do not say $\Phi' =\fG\Delta'$ because the prior switching means that the edge $ge$, in $E(\Phi')$ as a subset of $E(\fG\Delta)$, may not have gain $g$ in $\Phi'$.) 
\end{proof}


\section{Extension of biased expansions}\label{extension}

An \emph{extension} of a biased expansion $\Omega \downarrow_p \Delta$ is a biased expansion $\Omega'\downarrow_{p'}\Delta'$ such that 
\begin{enumerate}
\item[(a)]  $\Delta$ is a spanning subgraph of $\Delta'$, and 
\item[(b)]  $\Omega'\big|_\Delta = \Omega$ (so that $p'\big|_\Delta = p$).  
\end{enumerate} 
We may say $\Omega'$ is an extension of $\Omega$ \emph{to} $\Delta'$, or \emph{to} $E(\Delta')\setm E(\Delta)$.  The extension is \emph{simple} if $\Delta'$ is a simple graph.  It is a \emph{maximal extension} if it has no simple proper extension.

We are interested in two types of extensions.  The first is extension to a link that is parallel to an existing edge of $\Delta$.  

\begin{exam}[Parallel Extension]\label{XX:parallel} 
Suppose $\Delta$ is any graph, $f$ is a link in $\Delta$, and $e$ is an edge parallel to $f$ but not in $\Delta$.  $\Omega$ always extends to $e$.  
Take $(p')\inv(e)$ to be a set in one-to-one correspondence with $p\inv(f)$;  form balanced digons $\{\te,\tf\}$ when $\te$ and $\tf$ correspond; and for a circle $P\cup f$ in $\Delta$, a lift $\tP \cup \te$ is balanced in $\Omega'$ if and only if $\tP \cup \tf$ is balanced in $\Omega$.  
Any selection of edges of $\Delta$ can be reduplicated in this way, as many times as desired.  
\end{exam}

This kind of extension can be technically useful, but the other kind is the more important one:  that is extension by an edge $e_{vw}$ joining nonadjacent nodes of $\Delta$.  The possibility or impossibility of such extension is crucial data about the structure of a biased expansion.

There are four principal extension theorems.  First is uniqueness (Theorem \ref{TX:unique}).  If a biased expansion of a 2-connected graph $\Delta$ extends to one of $\Delta'$, that extension is \emph{unique}, by which we mean unique up to an isomorphism that is the identity on $\Omega\downarrow \Delta$.  We can express this by the existence of a commutative diagram of extensions:
\begin{equation*}
\xymatrix{
 &\Omega \ar@{_{(}->}[dl] \ar@{->}[dd]\ar@{^{(}->}[dr] &\\
\Omega' \ar@{-->}[rr]_>>>>>>\rho \ar@{->}[dd] &&\Omega''\ar@{->}[dd]\\
 &\Delta \ar@{^{(}->}[rd] \ar@{_{(}->}[ld] & \\
\Delta' \ar@{-->}[rr]_-\rho &&\Delta''
}
\end{equation*}
where the maps from $\Omega$ are embeddings and $\rho$ is an isomorphism.  In fact, $\rho$ itself is unique.  
(Recall that a mapping of biased expansions includes a mapping of their base graphs that commutes with projection.)

The second result is the existence of a unique maximal (simple) extension (Theorem \ref{TX:max}).  
The third result says that, if $e$ joins the trivalent nodes of a theta graph in $\Delta$, then $\Omega$ extends to $e$ (Lemma \ref{LX:theta}).  
Last is the theorem that, if $e$ is a chord of a circle $C\subseteq 
\Delta$, and $\Omega\big|_C$ extends to $e$, then $\Omega$ extends to $e$ 
(Proposition \ref{TX:chord}).

\begin{thm}[Uniqueness of Extension]\label{TX:unique}  
Let $\Omega \downarrow \Delta$ be a biased expansion of a 2-connected graph $\Delta$.  If $\Omega'\downarrow \Delta'$ and $\Omega''\downarrow \Delta'$ are two extensions of $\Omega$ to $\Delta'$, then there is a unique biased-expansion isomorphism $\rho : \Omega' \to \Omega''$ such that $\rho\big|_\Omega$ is the identity, provided that $\Delta'$ is simple or, more generally, that $\rho$, the projections, and $\id_{\Delta'}$ commute.  
\end{thm}

\begin{proof}  
Let $e\in E(\Delta') \setm E(\Delta)$ with endpoints $v$ and $w$.  These nodes lie in a common circle in $\Delta$; let $P_0$ and $P$ be the two $vw$-paths constituting the circle.  To define $\rho(\te')$ for $\te' \in (p')\inv(e)$ we choose $\tP_0$ so that $\tP_0 \cup \te'$ is balanced in $\Omega'$, then $\te'' \in (p'')\inv(e)$ so that $\tP_0\cup \te''$ is balanced, and set $\rho_e(\te') = \te''$.  It is clear that $\rho_e$ is a bijection $(p')\inv(e)\to (p'')\inv(e)$ because the roles of $\Omega'$ and $\Omega''$ are reversible.  

We have to prove that $\te''$ is independent of the choice of $\tP_0$.  Take $\tP$ so that $\tP\cup \te'$ is balanced; then $\tP_0 \cup \tP$ is balanced. Since $\tP_0\cup \te''$ and $\tP_0 \cup \tP$ are balanced, so is $\tP \cup \te''$.  Now suppose we change $\tP_0$ to $\tP^1_0$ so that $\tP^1_0 \cup \te'$ is balanced.  Then $\tP^1_0\cup \tP$ is balanced (because $\tP\cup \te'$ is), and since $\tP\cup \te''$ is balanced, so is $\tP^1_0\cup \te''$.  Therefore, $\rho_e(\te')$ is independent of the choice of lift of $P_0$.

Still, we ought to prove $\rho_e(\te')$ is independent of the choice of $vw$-path $P_0$.  Obviously, $P_0$ could be any $vw$-path.  Then suppose $P$ is a $vw$-path such that $P_0\cup P\cup e$ forms a theta graph.  Let $R_0$, $R$, and $R_e$ be the constituent paths of this theta graph that, respectively, lie in $P_0$, lie in $P$, and contain $e$.  Fixing a lift of $R$, one can imitate the previous proof to show that any lifts of $P_0$ and $P$ imply the same bijection $\rho_e$.  

Now take the original $P_0$ and any other $vw$-path $P$.  By the Path Theorem (see Corollary \ref{C:epath}) and the preceding argument, all of $P_0,P_1,\hdots,P_k = P$ induce the same bijection $\rho_e$.  

We now define $\rho(\tf)$, for $\tf \in E(\Omega')$, to be $\tf$ if $f\in E(\Delta)$ and $\rho_f(\tf)$ if $f\not\in E(\Delta)$. It remains to prove that $\rho : E(\Omega') \to E(\Omega'')$ is an isomorphism of biased graphs.  

For that it suffices to show that, if $\tC$ is a balanced circle in $\Omega'$, then $f(\tC)$ is balanced in $\Omega''$, and conversely.  Choose a spanning tree $T$ of $\Delta$ and a lift $\tT$ such that $\tT\cup \tC$ is balanced (possible by Lemma \ref{L:balsub}).  Then $\tC\subseteq \bcl_{\Omega'}\tT$.  By the definition of $\rho$, $\rho(\tC) \subseteq \bcl_{\Omega''}\tT$.  Since the latter is balanced, $\rho(\tC)$ is balanced. This reasoning works in both directions:  if $\tC'' \in \cB(\Omega'')$, then $\rho\inv(\tC'')$ is balanced.  That concludes the proof.
\end{proof}

\begin{thm}[Maximal Extension]\label{TX:max} 
Given any biased expansion $\Omega$ of a 2-connected simple graph $\Delta$, there is a unique maximal extension of $\Omega$; its base graph is $\Delta\cup X$ where 
\begin{equation*}
X = \{ e\notin E(\Delta) : \Omega \text{ extends to } e \}.
\end{equation*}
\end{thm}

Remember that ``uniqueness'' is up to isomorphisms that are the identity on $\Omega$.

\begin{proof}  
No extension $\Omega''\downarrow \Delta''$ can possibly have
$\Delta''\not\subseteq \Delta\cup X$, so we need only produce an extension of $\Omega$ to $\Delta\cup X$ and call upon the Uniqueness Theorem.

For each $e\in X$, let $\Omega_e$ be an extension to $e$.  The major part of the proof is to show that $\Omega_{e_1}$ and $\Omega_{e_2}$ are compatible.  

\begin{lem}[Common Extension]\label{LX:common}  
If $e_1,e_2\in X$, then $\Omega$ extends to $\Delta\cup \{e_1,e_2\}$.  
\end{lem}

\begin{proof}  The core of the lemma is the definition of balance in the common extension $\Omega_{12}$ of $\Omega_1$ and $\Omega_2$ (meaning $\Omega_{e_1}$ and $\Omega_{e_2}$).  The graph $\|\Omega_{12}\|$ is simply $\|\Omega_1\| \cup \|\Omega_2\|$.  
Balance of a circle $\tC$ that covers a circle $C$ in $\Delta_{12} = \Delta\cup \{e_1,e_2\}$ is as in $\Omega_1$ or $\Omega_2$, except when $C$ contains both $e_1$ and $e_2$.  
Then we define $\tC$ to be balanced if and only if there is a path $P$ that forms with $C$ a theta graph whose three constituent paths are $P$ and two paths, $P_1$ and $P_2$, of which $P_1$ contains $e_1$ and $P_2$ contains $e_2$ (we call $P$ a \emph{connecting chordal path} of $C$ because it connects the two components of $C \setm \{e_1,e_2\}$), and $P$ has a lift $\tP$ such that $\tC_1$ and $\tC_2$ are both balanced.  (The notation is that $C_i = P_i\cup P$, $\tP_i$ is the lift of $P_i$ that is contained in $\tC$, and $\tC_i = \tP_i\cup \tP$.) 
Then $\Omega_{12}$ is $\|\Omega_{12}\|$ with balanced circles as just defined.

It is important to know that the definition of balance is independent of the various choices implicit in it.  We need a bit more notation.  For a connecting chordal path path $P$ and edge $e\in P$, let $R = P\setm e$.  For a different connecting chordal path path $P'$ and $e' \in P'$, we define $R'$, $C'_1$, $C'_2$ analogously to $R$, $C_1$, $C_2$.  We begin with a little lemma. 

\begin{lem}\label{LX:chordalpath}  
Let $P$ and $P'$ be two connecting chordal paths of $C$ such that $(C\setm \{e_1,e_2\}) \cup P\cup P'$ contains a unique circle, $D$.  Let $e\in P\setm P'$ and $e' \in P'\setm P$, or let $e = e' \in P\cap P'$.  Let $\tC$ be a lift of $C$ and choose arbitrary lifts $\tR$ and $\tR'$ that agree on $R\cap R'$ and such that $\tD$ is balanced if $e = e' \in P \cap P'$.  Then, for each lift $\te$ such that $\tC_1$ and $\tC_2$ are balanced (in $\Omega_1$ and $\Omega_2$, respectively), there is a unique lift $\te'$ such that $\tC'_1$ and $\tC'_2$ are balanced (in $\Omega_1$ and $\Omega_2$, respectively).
\end{lem}

\begin{proof} 
Let $A = C\setm \{e_1,e_2\}$.  Note that $D \subseteq A \cup R\cup R'$ if $e = e'$, but $e,e' \in D$ if $e\neq e'$.  

Suppose $e = e'$.  In $\Omega_1$, $\tC_1 \cup \tD$ is a theta graph, $\tC'_1 = \tC_1\oplus \tD$, and $\tD$ is balanced, so $\tC_1$ is balanced if and only if $\tC'_1$ is.  Similarly, $\tC_2$ is balanced if and only if $\tC'_2$ is.  Also, the $\te$ that makes $\tC_i$ balanced is unique, by the circle lifting property in $\Omega_i$.  

If $\te\neq \te'$, then for each lift $\te$ there is a unique $\te' = \theta(\te)$ for which $\tD$ is balanced (in $\Omega)$, and $\theta$ is a bijection from $p\inv(e)$ to $p\inv(e')$. Suppose we lift $e$ to $\te$ such that $\tC_1$ is balanced.  Then lifting $e'$ to $\te'$, $\tC'_1$ is balanced $\iff$ $\tD$ is balanced $\iff \te' = \theta(\te)$. A similar argument applies to $\tC_2$ and $\tC'_2$.  
\end{proof}

The next lemma shows, in particular, that the definition of balance of a lift of $C$ is independent of the choice of connecting chordal path.

\begin{lem}\label{LX:balcircle}
Given $C$ containing $e_1$ and $e_2$, any connecting chordal path $P$, any edge $e\in P$, and any lift $\tR$ of $R= P\setm e$, then a lift $\tC$ is balanced if and only if there exists $\te$ such that $\tC_1$ and $\tC_2$ are balanced, and this $\te$ is unique.
\end{lem}

\begin{proof} 
By definition, $\tC$ is balanced if there is $\te$ such that $\tC_1$ and $\tC_2$ are balanced.  

Suppose, conversely, that $\tC$ is balanced; thus, there exist a connecting chordal path $P'$ and a lift $\tP'$ such that $\tC'_1$ and $\tC'_2$ are balanced.  Choose $e' \in P'$.

Since $\Delta$ is inseparable, there exist two connecting chordal paths of $C$, $Q$ and $Q'$, that are internally disjoint.  
By the Path Lemma \ref{C:epath} there is a chain $P' = Q_0,Q_1,\hdots,Q_k = P$ of connecting chordal paths that includes $Q=Q_{m-1}$ and $Q'=Q_m$, such that $(C\setm\{e_1,e_2\}) \cup Q_{i-1} \cup Q_i$ contains a unique circle for all $i$.  
Choose edges $f_i\in Q_i$ so that $f_0 = e'$, $f_k = e$, and $f_{i-1} \in Q_i\Rightarrow f_i = f_{i-1}$.  This is possible because $Q_{m-1}$ and $Q_m$ are edge disjoint.  Thus, at worst we may be forced to take $f_{m-1}=f_0$ and $f_m=f_k$, but there is no necessary relation between $f_0$ and $f_k$.

Let $R_i = Q_i\setm f_i$ and let $C_{i1}$ and $C_{i2}$ be the circles in $C\cup Q_i$ that contain, respectively, $e_1$ and $e_2$ but not both. 
We may apply Lemma \ref{LX:chordalpath} $k$ times to conclude that, for any lifts $\tR_0$ and $\tR_k$, in particular, $\tR_0 \subseteq \tP'$ and $\tR_k = \tR$, and for any $\te'$ such that $\tC_{01}$ and $\tC_{02}$ are balanced, there is a unique lift $\te$ such that$\tC_{k1}$ and$\tC_{k2}$ are balanced.  (If $\theta_i : p\inv(f_i)\to p\inv(f_{i-1})$ is as in the proof of Lemma \ref{LX:chordalpath}, then $\te = (\theta_1\theta_2 \cdots \theta_k)\inv(\te)$.)
\end{proof}

To prove Lemma \ref{LX:common} we need just two more steps: to prove, first, the circle lifting property in $\Omega_{12}$, and second, that $\cB(\Omega_{12})$ is a linear class.

{\textit{Step 1.  Circle lifting.}}  
We need to consider a circle $C\owns e_1,e_2$ and an edge $f\in C$.  Letting $S = C\setm f$, we assume $\tS$ given and must prove there is a unique $\tf$ such that $\tS\cup \tf$ is balanced.  We take $P$, $e\$, and $R\$ as in Lemma \ref{LX:balcircle}, and fix $\tR$.  We may assume $f\in P_2$.  

Choose $\te$ so that $\tC_1\subseteq \tS \cup \tP$ is balanced in $\Omega_1$, then $\tf$ so that $\tC_2\subseteq \tS\cup\tP\cup \tf$ is balanced in $\Omega_2$.  By definition, $\tC = \tS \cup \tf$ is then balanced.  Suppose both $\tf^1$ and $\tf^2$ make $\tS \cup \tf^i = \tC^i$ balanced.  By Lemma \ref{LX:balcircle}, for $i=1$ and $2$, 
\begin{equation*}
(\exists \te^i) \ \tP_1\cup \tR\cup \te^1\ \text{ and }\ 
\tP^i_2\cup \tR\cup \te^i\ \text{ are balanced,}
\end{equation*}
where $\tP^i_2$ is the lift of $P_2$ contained in $\tS\cup \tf^i$.  Comparing $\tC^1_1$ with $\tC^2_1$ in $\Omega_1$, $\te^1 = \te^2$.  Then, comparing $\tC^1_2$ with $\tC^2_2$ in $\Omega_2$, $\tf^1 = \tf^2$.  Thus, $\tf$ is unique.

{\textit{Step 2.  Linearity.}}  
We must examine lifts of a theta graph $\Theta$ that contains both $e_1$ and $e_2$.  There are two cases, according as $e_1$ and $e_2$ are in the same or different paths of $\Theta$.  

If $e_1$ and $e_2$ are in different paths, we can use the notation of Lemma \ref{LX:balcircle}.  Suppose a lift such that $\tC_1$ and $\tC_2$ are balanced:  then $\tC$ is balanced by definition.  On the other hand, suppose $\tC_1$ and $\tC$ are balanced while $\tC_2$ is unbalanced.  By changing $\te_2$ to $\te^2_2$ we get a balanced lift of $C_2$, namely, $\tC^2_2 = (\tC_2\setm \te_2) \cup \te^2_2$.  Then $\tC^2 = (\tC\setm \te_2) \cup \te^2_2$ is balanced.  However, in Step 1 we showed that $\tC$ and $\tC^2$ cannot both be balanced.  Therefore, $\tC_2$ must have been balanced after all.

Suppose now that $e_1$ and $e_2$ lie in the same path of $\Theta$.  Then $\Theta\setm \{e_1,e_2\}$ has two components and contains a unique circle, call it $D$.  Let $B$ and $C$ be the other circles in $\Theta$, and let $P$ be a minimal path in $\Delta$ connecting the components of $\Theta\setm \{e_1,e_2\}$.  We may assume that $P$ has both endpoints in $N(C)$, so that $C\cup P$ is a theta graph with circles $C_1\owns e_1$ and $C_2 \owns e_2$.  If we write the path $B\cap C$ as a concatenation of paths, $R_1e_1Re_2R_2$, we may also assume that $P$, which has one endpoint in $R$, has the other end not in $R_2$.  Therefore, $B\cup C_1$ and $D\cup C_2$ are theta graphs.  In addition, $B\oplus C_1 = D\oplus C_2$.  

Suppose $\tB$ and $\tC$ are balanced.  By Lemma \ref{L:balsub}, there is a lift $\tP$ such that $\tC_1$ and $\tC_2$ are balanced.  Therefore $\tB \oplus \tC_1$ is balanced, and as this equals $\tD\oplus \tC_2$ and $\tC_2$ is balanced, $\tD$ is balanced.

If, however, it is $\tC$ and $\tD$ that are balanced, then $\tD \oplus \tC_2 = \tB \oplus \tC_1$ is balanced, whence $\tB$ is balanced. 

Supposing finally that $\tB$ and $\tD$ are balanced, we choose $\tP$ so that $\tC_1$ is balanced.  Consequently, $\tB \oplus \tC_1$ is balanced.  This being $\tD \oplus \tC_2$, we conclude that $\tC_2$ is balanced, whence $\tC$ is balanced.

Thus in every case linearity is satisfied, and therefore, $\Omega_{12}$ is a biased graph.
\end{proof}

\begin{lem}\label{LX:edge}  
Suppose $\Omega'$ extends $\Omega$ to $\Delta'$ and $e\in X\setm E(\Delta')$; then there is an extension of $\Omega$ to $\Delta'\cup e$.  
\end{lem}

\begin{proof}  
Let 
\begin{equation*}
\cF  = \{\Delta'' \subseteq \Delta' : \Omega'\big|_{\Delta''} \text{ extends to } e\}.
\end{equation*}  
If $\Delta'' \in \cF$ and $f\in E(\Delta')\setm E(\Delta'')$, then $\Omega'\big|_{\Delta''}$ extends both to $e$ and to $f$; by Lemma \ref{LX:common} it extends to $\{e,f\}$, so $\Delta'' \cup f \in \cF$.  This suffices to prove the lemma when $X$ is finite.

Otherwise, we apply Zorn's Lemma in the usual way. Take a minimal chain $\{\Delta_i\}$ in $\cF$; let $\Delta''$ be its union.  Write $\Omega_i$ for the extension of $\Omega'\big|_{\Delta_i}$ to $e$.  By Unique Extension we can regard each $\Omega_i$ for $i < j$ as the restriction $\Omega_j\big|_{\Delta_i\cup e}$.  Therefore $\Omega'' = \bigcup_i \{\Omega_i\}$ is a well defined graph.  It is a biased expansion of $\Delta''\cup e$ because any circle in $\Delta''\cup e$ or theta graph in $\Omega''$ is contained in some $\Delta_i \cup e$ or $\Omega_i$.  It extends $\Omega'\big|_{\Delta''}$ because $\Omega'\big|_{\Delta''} = \bigcup_i\{\Omega'\big|_{\Delta_i}\}$.  Therefore, $\Delta'' \in \cF$.  If $\Delta'' \subset \Delta'$, there is an $f\in E(\Delta')\setm E(\Delta'')$ and, by the first part of the proof, $\Delta'' \cup f \in \cF$.  As that contradicts the maximality of the original chain, $\Delta''$ must be $\Delta'$, so $\Omega'$ extends to $e$.  
\end{proof}

\begin{lem}\label{LX:all}  
Suppose $\Omega'$ extends $\Omega$ to $\Delta' \supset \Delta$; then $\Omega'$ extends to $\Delta\cup X$.  
\end{lem}

\begin{proof}  
Here let 
\begin{equation*}
\cF = \{\Delta'' \subseteq \Delta \cup X\mid \Delta'' \supseteq \Delta'\ \text{ and }\ \Omega'\
\text{ extends to }\ \Delta''\}.
\end{equation*}
There can be only one maximal member of $\cF$, namely, $\Delta \cup X$, since for any other $\Delta''$, taking $e\in X\setm E(\Delta'')$ we know that $\Omega''$, an extension of $\Omega'$ to $\Delta''$, extends to $e$.  This proves the lemma when $X\setm E(\Delta')$ is finite.  

In the infinite case, again we apply Zorn's Lemma. The union of a maximal chain of graphs in $\cF$ is itself in $\cF$, and this union must be $\Delta\cup X$ or the chain could not have been maximal.  
\end{proof}

To complete the proof of the Maximal Extension Theorem we need only appeal to the Unique Extension Theorem.  
\end{proof}

\begin{lem}[Theta Extension]\label{LX:theta}  
Any biased expansion of a theta graph with trivalent nodes $v$ and $w$ extends to the edge $e_{vw}$. 
\end{lem}

\begin{proof}  
Let the theta graph $\Delta$ have constituent paths $P_1$, $P_2$, and $P_3$ and write $e$ for $e_{vw}$.  By Example \ref{XX:parallel} we may assume $v$ and $w$ are nonadjacent in $\Delta$.  Define a set $E_e$ in one-to-one correspondence with some fiber $p\inv(f)$ for $f\in E(\Delta)$.  Letting each $\te \in E_e$ have endpoints $v$ and $w$ defines a graph $\|\Omega'\|$ that covers $\Delta \cup e$.  
The task is to define balance and show it results in a biased expansion $\Omega'\downarrow \Delta\cup e$ extending the original biased expansion $\Omega \downarrow \Delta$.  

Choose a fixed edge $f_1\in P_1$, let $Q_1 = P_1\setm f_1$, and fix a lift $\tQ^0_1$. Choose a bijection $\psi : p\inv(f_1) \to E_e$ and, for $\te \in E_e$, define 
$$
\tQ^0_1 \cup \tf_1 \cup \te \text{ balanced } \iff \te = \psi(\tf_1).
$$

For any lift $\tP_2$ and any $\te \in E_e$, we define 
$$
\tP_2 \cup \te \text{ balanced } \iff \tP_2 \cup \tQ_1^0 \cup \psi\inv(\te) \text{ is balanced.}
$$
For $\tP_3\cup \te$ the definition is similar.  (This leaves balance of $\tP_1\cup \te$ undefined as yet, in general.)  We need to show consistency between the states of balance of $\tP_2\cup \te$ and of $\tP_3\cup \te$.  If both are balanced, $\tP_2 \cup \tQ^0_1 \cup \psi\inv(\te)$ and $\tP_3\cup \tQ^0_1 \cup \psi\inv(\te)$ are balanced, so $\tP_2\cup \tP_3$ is balanced.  Similarly, if only one of $\tP_2\cup \te$ and $\tP_3 \cup \te$ is balanced, $\tP_2\cup \tP_3$ cannot be balanced.  Thus, linearity is satisfied for lifts of $P_2\cup P_3\cup e$.  We call this \emph{23-consistency}.  

Now, for a lift $\tP_1$ we define $\tP_1\cup \te$ to be balanced if $\tP_2 \cup \te$ is balanced for some $\tP_2$ such that $\tP_1\cup \tP_2$ is balanced.  Suppose we took two lifts $\tP^1_2$ and $\tP^2_2$ such that both $\tP_1\cup \tP^j_2$ are balanced, and say $\tP^j_2\cup \te^j$ is balanced.  Pick $\tP_3$ so that $\tP_1\cup \tP_3$ is balanced.  Then each $\tP^j_2\cup \tP_3$ is balanced.  By 23-consistency, $\tP_3 \cup \te^j$ is balanced for $j = 1,2$; thus $\tQ^0_1\cup \tf^j_1$ is balanced for $\tf^j_1 = \psi\inv(\te^1)$,  but since $\tf^1_1 = \tf^2_1$, we see $\te^1 = \te^2$.  Therefore, balance of $\tP_1\cup \te$ is independent of the choice of $\tP_2$.  We call this \emph{12-consistency}.  

We show that, if $\tP_1\cup \tP_3$ is balanced, then $\tP_1\cup \te$ is balanced if and only if $\tP_3 \cup \te$ is balanced.  Take $\tP_2$ so that $\tP_1\cup \tP_2\cup \tP_3$ is balanced.  Then $\tP_1\cup \te$ is balanced $\iff$  (by 12-consistency) $\tP_2\cup \te$ is balanced $\iff$ (by 23-consistency) $\tP_3 \cup \te$ is balanced.  

We still have to prove uniqueness in the circle lifting property.  First, we treat lifts of $e$.  Any $\tP_i$ has a balanced completion $\tP_i \cup \te$, as we have seen.  Suppose $\tP_i\cup \te^1$ and $\tP_i\cup \te^2$ are balanced.  If $i = 2,3$, just take $\tQ^0_1\cup \tf_1$ such that $\tP_i\cup \tQ^0_1\cup \tf_1$ is balanced.  Then $\tf_1 = \psi\inv(\te^j)$ for $j = 1,2$, whence $\te^1=\te^2$.  If $i = 1$, take $\tP_2$ so that $\tP_1\cup \tP_2$ is balanced:  then $\tP_2\cup \te^1$ and $\tP_2\cup \te^2$ are balanced, so $\te^1 = \te^2$.  

Now we treat lifts of $f\in P_i$.  Let $R = P_i\setm f$ and take any $\tR$ and $\te$.  If $i = 2,3$, we know that $\tQ^0_1\cup \psi\inv(\te) \cup \te$ is balanced, and there exists $\tf$ for which $\tQ^0_1\cup \psi\inv(\te) \cup \tR\cup \tf$ is balanced (and it is unique); by definition, for this $\tf$ and no other, $(\tR\cup \tf) \cup \te$ is balanced.  If $i = 1$, we choose any $\tP_2$ such that $\tP_2\cup \te$ is balanced.  Then there is a unique $\tf$ making $\tP_2\cup \tR\cup \tf$ balanced, and by definition that is the only $\tf$ for which $(\tR\cup \tf)\cup \te$ can be balanced.  Thus we have a biased expansion of $\Delta\cup e$.  
\end{proof}

\begin{prop}[Chordal Extension]\label{TX:chord}  
Suppose $\Omega$ is a biased expansion of a $2$-connected graph $\Delta$ and $e\not\in E(\Delta)$.  For any circle $C\subseteq \Delta$ of which $e$ is a chord, $\Omega$ extends to $e$ if and only if $\Omega\big|_C$ extends to $e$.
\end{prop}

\begin{proof}  
We need only prove sufficiency.  Take $C$, of which $e$ is a chord, such that $\Omega\big|_C$ extends to $e$.  Let $P_1$ and $P_2$ be the paths into which $e$ divides $C$.  Let $\Omega_e$ be the extension to $e$ of $\Omega\big|_C$. To define $\Omega'$, the extension of $\Omega$, we set $E(\Omega') = E(\Omega) \cup p\inv_e(e)$ and define a circle $\tP\cup \te$ in $\Omega'$, lifting a circle $P\cup e$ in $\Delta\cup e$, to be balanced if and only if there is a lift $\tP_1$ such that both $\tP\cup \tP_1$ and $\tP_1\cup \te$ are balanced.  It remains to prove that $\Omega'$ is a biased graph and a biased expansion of $\Delta'$. First we show that $P_2$ works as well as $P_1$ in defining balance of $\tP\cup \te$.  

\begin{lem}\label{LX:p1p2}  
$\tP\cup \te$ is balanced if and only if there is a choice of $\tP_2$ so that $\tP\cup \tP_2$ and $\tP_2\cup \te$ are balanced.  
\end{lem}

\begin{proof}  
First, suppose $\tP\cup \te$ is balanced:  then there is a $\tP_1$ such that $\tP\cup \tP_1$ and $\tP_1\cup \te$ are balanced.  Choose $\tP_2$ so that $\tP\cup \tP_1\cup \tP_2$ is balanced.  (That is possible by Lemma \ref{L:balsub}.)  Then $\tP_1\cup \tP_2\cup \te$ is a theta graph in $\Omega_e$, so $\tP_2\cup \te$ is balanced.  Thus, $\tP_2$ exists as desired.

The converse is similar.
\end{proof}

We should prove that different choices of $\tP_1$ give consistent definitions of balance of $\tP\cup \te$.  

\begin{lem}\label{LX:p1indep}  
Suppose $\tP^1_1\cup \tP$ and $\tP^2_1\cup \tP$ are balanced, and $\tP^1_1\cup \te^1$ and $\tP^2_1 \cup \te^2$ are balanced.  Then $\te^1 = \te^2$.  
\end{lem}

\begin{proof}  
Choose $\tP_2$ so $\tP\cup \tP^1_1\cup P_1$ is balanced (by Lemma
\ref{L:balsub}).  Then $\tP_2\cup \te^1$ is balanced in $\Omega_e$, because of the theta graph $\tP^1_1\cup \te^1\cup \tP_2$.

If we remove from $P_1$ the edges of $P$, we are left with $k\geq 1$ segments $S_1,\hdots,S_k$ of positive length.  Choose $e_i\in S_i$ and let $R = P_1\setm \{e_1,\hdots,e_k\}$.  Then 
\begin{enumerate}
\item[(1)]  $R\cup P$ is connected, so $R\cup P\cup P_2$ is connected, and 
\item[(2)]  no edge of any $S_i$ is contained in any circle of $R\cup P$, nor of $R\cup P\cup P_2$.
\end{enumerate}
Consequently, writing $\tR^2$ for the lift of $R$ contained in $\tP^2_1$, 
\begin{enumerate}
\item[(3)]  $\tR^2\cup \tP$ and $\tR^2\cup \tP\cup \tP_2$ are connected, and 
\item[(4)]  $\tR^2\cup \tP\cup \tP_2$ is balanced, because any circle in it lies in $\tP\cup \tP_2$, which is balanced.
\end{enumerate}

Now, $\te^2_i$ lies in a circle in $\tR^2\cup \tP\cup \te^2_i$ by (3), which is balanced because it is in $\tP^2_1\cup \tP$.  Therefore $\te^2_i\in \bcl(\tR^2\cup \tP)$.  So 
\begin{equation*}
\tP^2_1\cup \tP\cup \tP_2\subseteq \bcl(\tR^2\cup \tP\cup \tP_2),
\end{equation*}
which is balanced (Lemma \ref{L:bcl}). Thus $\tP^2_1\cup \tP_2$ is balanced, so $\tP_2\cup \te^2$ is balanced by the theta graph $\tP^2_1\cup \te^2\cup \tP_2$ in $\Omega_e$.  As $\Omega_e$ is a biased expansion and both $\tP_2\cup \te^i$ are balanced, $\te^1=\te^2$.  
\end{proof}

Thus, we have a well defined notion of balance in $\|\Omega'\| = \|\Omega\|\cup \|\Omega_e\|$.

The lemma applies as well to $P_2$ as to $P_1$, of course, due to Lemma \ref{LX:p1p2}.  

\begin{lem}\label{CX:p1choice}  
Suppose $\tP_1$ chosen so that $\tP\cup \tP_1$ is balanced. Then $\tP\cup \te$ is balanced if and only if $\tP_1\cup \te$ is balanced.
\end{lem}

\begin{proof}  
There is a unique $\te^0$ for which $\tP_1\cup \te^0$ is balanced, because $\Omega_e$ is a biased expansion.  Then $\tP\cup \te^0$ is balanced, but by Lemma \ref{LX:p1indep} no other $\tP\cup \te$ can be balanced.
\end{proof}

\begin{lem}\label{CX:cchoice}  
Suppose $\tC$ chosen so that $\tP\cup \tC$ is balanced.  Then $\tP \cup \te$ is balanced (in $\Omega'$) if and only if $\tC \cup \te$ is balanced (in $\Omega_e$).
\end{lem}

\begin{proof}  
Apply Lemma \ref{CX:p1choice} to $P_1$ and $P_2$, the latter requiring Lemma \ref{LX:p1p2}.  
\end{proof}

The rest of the proof shows that $\Omega'$ is a biased expansion. First, the uniqueness of circle lifting.  

\begin{lem}\label{LX:uniqueexpan}  
If $C'$ is a circle in $\Delta'$ and $f\in C'$, and if $\tP'$ is any lift of $C'\setm f$, then there is exactly one lift $\tf$ that makes $\tP' \cup \tf$ balanced.
\end{lem}

\begin{proof}  We may assume $e\in C'$.  When $f=e$, this is a consequence of Lemma \ref{CX:p1choice}.  Otherwise, let $P = C'\setm e$, so $f\in P$.  (We may assume $P\neq P_1,P_2$.)

We shall have need of the graph of $P$, which is $(N(P),P)$, and that of $C$.  Removing $f$, $(N(P),P)$ falls into two connected halves, one containing $v$ and the other $w$; we write $R = (N(P),P)\setm f$ and $R_v$, $R_w$ for the first two halves.  We shall be careless with notation, using $P$, $R$, etc., to denote both the graph and the edge set, trusting that all will be clear.  

A bridge of $C$ in $C\cup P$ is a maximal subpath of $P$ whose internal nodes are in $P\setm C$, and a bridge of $P$ in $C\cup P$ is a maximal subpath of $C$ whose internal nodes lie in $C\setm P$ (excluding edgeless subpaths in both cases).  Call the bridges of $C$ (which are subpaths of $P$) $S_1,S_2,\hdots,S_m$ and choose an edge $s_i\in S_i$ for each bridge.  Let $S = \{s_1,s_2,\hdots,s_m\}$.  Amongst the bridges of $P$ (which are subpaths of $C$), we are interested only in those that connect $R_v$ to $R_w$.  For each such bridge choose an edge $d_i$ in it, and let $D = \{d_1,\hdots,d_k\}$, there being $k$ such bridges. $D$ depends on $f$.  Let $D' = C\setm D$ (as an edge set).

So far, we have two biased expansions:  $\Omega\downarrow \Delta$ and $\Omega_e\downarrow C\cup e$.  The Theta Extension Lemma generates others, which we employ as auxiliary graphs.  Since $C\cup S_i$ is a theta subgraph of $\Delta$, $\Omega\big|_{C\cup S_i}$ extends to a chord $e_i$ of $C$ that joins the endpoints of $S_i$.  Call $\Omega_i$ the resulting biased expansion of $C\cup S_i\cup e_i$, and let $H = \{e_1,e_2,\hdots,e_m\}$.  Taking every $\Omega_i$ separately, we get extensions $\Omega_i\big|_{C\cup e_i}$ of $\Omega\big|_C$.  By the Maximal Extension Theorem, all the extensions of $\Omega\big|_C$, including $\Omega_e$, are compatible; that is, there is an extension $\Omega_P$ of $\Omega\big|_C$ to $H\cup e$, determined by $\Omega_e$ and the subpaths $S_i$.  We now have three groups of biased expansions: $\Omega$, $\Omega_P$ extending $\Omega_e$, and $\Omega_i$ extending $\Omega\big|_{C\cup S_i}$ to $e_i$.  All this is independent of $f$.  

We wish to prove that, given $\tR$ and $\te$, there is a unique $\tf$ such that $\tR \cup \{\te,\tf\}$ is balanced.  First we establish a tool.

\begin{lem}\label{LX:shortcircuit}  
Let $\tP$ and $\tC$ be arbitrary lifts of $P$ and $C$, and let $\te_i$ be the lift that makes $\tS_i \cup \tC$ balanced in $\Omega_i$, and let $\tH = \{\te_1,\hdots,\te_m\}$.  Then $\tP\cup \tC$ is balanced in $\Omega \iff  \tC\cup \tH$ is balanced in $\Omega_P$. 

Furthermore, let $\te$ be any lift of $e$ such that $\tC \cup \te$ is balanced.  Then $\tP \cup \te$ is balanced (in our definition given previously) $\iff \tC\cup \tH\cup \te$ is balanced in $\Omega_P$.   
\end{lem}

\begin{proof}  
For the first part, since $\tP \cup \tC$ is balanced, from $\Omega_i$ we know every $\tC\cup \te_i$ is balanced.  Thus, $\te_i \in \bcl_{\Omega_P} \tC$.  By Lemma \ref{L:bcl}, $\tC \cup \tH$ is balanced.

Conversely, if $\tC \cup \tH$ is balanced, then each $\tC \cup \tS_i$ is balanced.  Therefore, $\ts_i\in \bcl_\Omega(\tC\cup \tP\setm \tS)$.  $\tC\cup \tP\setm \tS$ is balanced because its only circle is $\tC$.  It follows that $\tC \cup \tP$ is balanced.

For the second part, because we assume balance of $\tC\cup \te$, $\tP\cup \te$ is balanced $\iff$ $\tC \cup \te$ is balanced in $\Omega_e$.  We can reformulate the statement as: $\tP \cup \tC$ is balanced (in $\Omega)$ $\iff \tC\cup \tH \cup \te$  is balanced (in $\Omega_P$).  The proof is like that of the first part.  
\end{proof}

Let $A = (P\cup D')\setm f$.  In case $f\in P\setm C$, $f$ lies in a subpath $S_t$ corresponding to a chord $e_t$.  If $f\in C$, we leave $S_t$ and $e_t$ undefined.  Define $B = (D'\cup H)\setm \{f,e_t\}$.  Then each of $A$ and $B$ contains $R$ but, due to the absence of $D$, $f$, and (when appropriate) $e_t$, remains disconnected into a $v$-component and a $w$-component.  Adding in any one of $e$, $d_i$ for $1\leq i\leq m$, or $f$ or $e_t$ makes $A$ and $B$ connected.

We were given $\tR$ and we can extend it (in $\Omega)$ to a balanced lift $\tA$.  Each $\tS_i \subseteq \tR$, except when $i = t$, implies a unique $\te_i$ for which $\tS_i\cup \te_i$ is balanced (in $\Omega_i$).  Thus we have a balanced lift $\tB$ (in $\Omega_P)$ as well, uniquely defined.  Now we add the given $\te$ and take $\bcl_{\Omega_P}(\tB\cup \te)$.  It contains exactly one lift $\tilde d_i$ for each $i$ and one $\tf$ (if $f\in C$) or $\te_t$ (if $f\not\in C$), and it is balanced.  Thus we have a balanced lift $\tB\cup \tD \cup \{\te,\tf\}$ (if $f\in C$) or $\tB \cup \tD\cup \{\te,\te_t\}$ (if not), which in both cases is $\tC\cup \tH \cup \te$.  Moreover, the lifts $\tD$ and (if $f\in C$) $\tf$ are the only ones that give balance, by the circle lifting property in $\Omega_P$.  When $f\in C$, Lemma \ref{LX:shortcircuit} shows that, not only is $\tP\cup \te$ balanced, but $\tf$ is the only lift of $f$ for which this is true.  When $f\not\in C$, we find $\tf$ as the unique edge in $p\inv(f)\cap \bcl_{\Omega_t}(\widetilde{S_t\backslash f} \cup \tC\cup \te_t)$.  Balance of $\tP\cup \te$ follows from the second part of Lemma \ref{LX:shortcircuit}.  In both cases, $\tf$ exists and is unique.
\end{proof}

\begin{lem}\label{LX:chordalbg} $\Omega'$ is a biased graph.  
\end{lem}

\begin{proof}  
We look at a theta graph that contains $e$.  Let $P\cup e$ and $P'\cup e$ be its circles that contain $e$ and $D = P\oplus P'$ the third circle.

Suppose, in a lift of $P\cup P'$, $\tD$ is balanced.  Since $\tP\cup \tP'$ is balanced, we can choose $\tP_1$ so that $\tP\cup \tP'\cup \tP_1$ is balanced.  Then for any $\te$, $\tP\cup \te$ is balanced $\iff \tP_1\cup \te$ is balanced $\iff \tP' \cup \te$ is balanced.  That is, one or three circles in $\tP\cup \tP'\cup \te$ are balanced.  Thus $\te$ is unique due to Lemma \ref{LX:p1indep}.  

Suppose, however, that $\tD$ is not balanced.  Take $f\in P\setm P'$ and replace $\tf \in \tP$ by $\tf^0$ such that $\tP^0\cup \tP^1$ is balanced.  (No other lift edges are altered.)  By the first part, $\tP^0\cup \te$ is balanced $\iff \tP'\cup \te$ is balanced.  If both are balanced, then $\tP\cup \te$ is unbalanced by Lemma \ref{LX:uniqueexpan}, so only one circle is balanced in $\tP \cup \tP'\cup \te$.  On the other hand, if neither is balanced, then $\tD$ and $\tP'\cup \te$ are unbalanced, so at most one circle is balanced in $\tP \cup \tP'\cup \te$.
\end{proof}

The combination of Lemmas \ref{LX:uniqueexpan} and \ref{LX:chordalbg} proves the Chordal Extension Theorem.
\end{proof}

Call a graph \emph{theta-complete} if the trivalent nodes of any theta subgraph are adjacent.  
The \emph{theta completion} $\theta(\Delta)$ of a simple graph $\Delta$ is the smallest theta-complete simple graph that contains $\Delta$.  
The results of this section imply:

\begin{prop}\label{CX:thetacompletion}
 A biased expansion of a simple graph $\Delta$ extends uniquely to 
$\theta(\Delta)$.  If $Delta$ is the base graph of a maximal biased 
expansion, then $\Delta$ is theta-complete.
 \end{prop}


\section{Inescapable groups (3-connection)}\label{group}

Our extension results imply a strong characterization of biased expansions of well-connected graphs.

\begin{thm}\label{T:3connbg}  
Every biased expansion of a $3$-connected graph of order at least four is a group expansion.  The group is unique.
\end{thm}

The graph being expanded may have finite or infinite order.

\begin{lem} \label{L:completebx}  
 A biased expansion of a complete graph of order four or more is a group 
expansion by a unique group.
\end{lem}

\begin{proof}[Proof of Lemma]   
Let $K$ be the complete graph and $\Omega$ its biased expansion.

 In the finite case the lemma is a consequence of the theorem of 
``generalized associativity'' stated by Belousov \cite{ASQ} and proved by 
Hossz\'u \cite{HosszuA}, Acz\'el, Belousov, and Hossz\'u \cite[Theorem 
1]{ABH}, and Belousov \cite{GASQ} (see \cite[pp.\ 76--78]{Latin}), and 
independently proved by Kahn and Kung \cite[Section 7, pp.\ 
490--492]{VCG}.
 ``Generalized associativity'' states that, if a set has four quasigroup 
operations that satisfy $g_1(h_1(x_1,x_2),x_3) = g_2(x_1,h_2(x_2,x_3)),$ 
all four operations are isotopic to the same associative quasigroup. 
 (It follows that, if a finitary quasigroup factors into binary 
quasigroups in all possible ways, then the quasigroup is an iterated group 
isotope.) 
 Acz\'el, Belousov, and Hossz\'u prove this by producing explicit 
isotopisms.  Kahn and Kung construct four quasigroups that satisfy the same equation, 
from combinatorial data equivalent to a biased expansion of $K_4$, in such 
a way that they have identity elements; thus they are obviously equal, 
hence a group. 
 Either way, it follows from generalized associativity that every 
restriction $\Omega\big|_{K'}$ to a $K_4$ subgraph $K'\subseteq K$ is a 
group expansion, say $\bgr{\fG_{K'}K'}$.  Since $\Omega\big|_{K''} \cong 
\bgr{\fG_{K'}K_3}$ for $K''\subseteq K'$ of order three, and $\fG_{K'}$ is 
unique up to isomorphism (by a theorem of Bruck \cite{Bruck}, or proved 
directly by Dowling \cite[Theorem 8]{CGL}), Kahn and Kung deduce that all 
$\fG_{K'}$ are isomorphic and $\Omega = \bgr{\fG K}$ for a group $\fG$.

In essence, these approaches depend on interpreting $\Omega\big|_C$ for a spanning circle $C$ of $K'$ as encoding a quasigroup \emph{multiplication}.  
Partly for completeness' sake and partly because it is such a natural way of deducing the group directly from the biased expansion, we give a new proof that depends on setting up the \emph{division} operation of the group by means of a spanning star subgraph of $K$.

Let $v_0 \in N(K)$, and distinguish a balanced lift $\tK^0$ of $K$.  Take a set $\fQ$ in one-to-one correspondence with each fiber $p\inv(e)$.  Holding $\tK^0\setm v_0$ fixed, and letting $v_1$ be another node of $K$, each choice of edge $\te_{01}$ implies by balance-closure one $\te_{0j}$ for all $j\neq 0,1$ such that the lift $\tK^0 \supseteq \tK^0\setm v_0$ is balanced.   
This determines bijections $\zeta_j : p\inv(e_{01})\to p\inv(e_{0j})$.  
We define $\psi_1$ to be any one bijection $\psi_1 : p\inv(e_{01}) \to \fQ$ and $\psi_j : p\inv(e_{0j}) \to \fQ$ to be the bijection $\zeta\inv_j \circ \psi$.  We also define 
\begin{equation*}
\epsilon = \psi_j(\te^0_{0j})\ \text{ for all }\ j\neq 0.
\end{equation*}
This will serve as the group identity.  

We have now labelled (from $\fQ$) all edges $\te_{0i}$.  The next task is to label all $\te_{ij}$.  We define 
\begin{equation}
\psi_{ij}(\te_{ij}) = \psi_i(\te_{0i})\ \text{ if }\ \te_{0i}\te^0_{0j}\te_{ij}\ \text{is balanced.}\label{E:labelbg}
\end{equation}
Finally, we define division.  Actually, we define an operation $(\alpha/\beta)_{ij}$ for each ordered pair of distinct $i,j\neq 0$ by 
\begin{equation*}
(\alpha/\beta)_{ij} = \psi_{ij}(\te_{ij})\ \text{ if }\
\psi\inv_i(\alpha)\psi\inv_j(\beta)\te_{ij} \ \text{ is balanced}.  
\end{equation*}
These definitions are illustrated in Figure \ref{F:completebg}.  
\begin{figure}[htbp]
\begin{center}
\psfrag{v0}[c]{$v_0$}
\psfrag{vi}[c]{$v_i$}
\psfrag{vj}[c]{$v_j$}
\psfrag{alp}[c]{$\alpha$}
\psfrag{bet}[c]{$\beta$}
\psfrag{gam}[c]{$\gamma$}
\psfrag{eps}[c]{$\epsilon$}
\psfrag{psia}[l]{$\psi_{ij}(\te_{ij})=\alpha$}
\psfrag{psib}[l]{$\psi_{ij}(\te_{ij})=(\alpha/\beta)_{ij}$}
\includegraphics{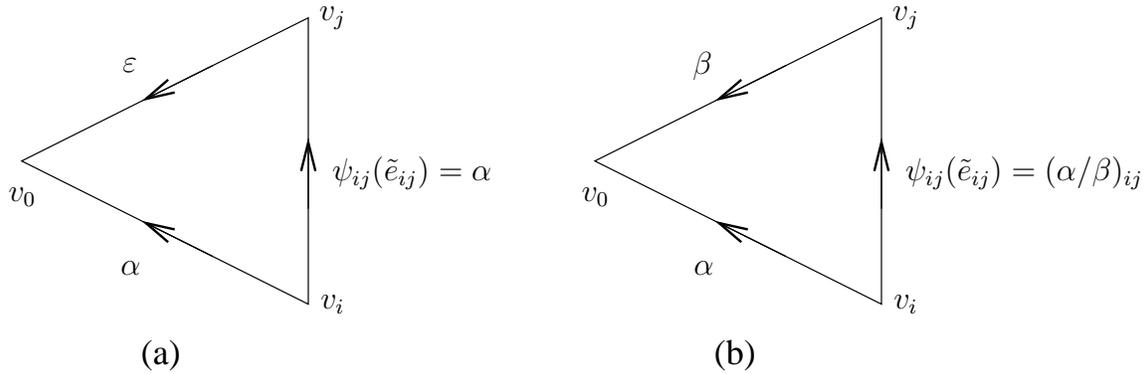}
\end{center}
\caption{(a)  The way an edge $\te_{ij}$ is labelled.  (b)  The definition of $(\alpha/\beta)_{ij}$.  All triangles are balanced.}\label{F:completebg}
\end{figure}

The next step is to prove that division is independent of the first subscript:
\begin{equation}
(\alpha/\beta)_{ik} = (\alpha/\beta)_{jk}.\label{E:indep} 
\end{equation}
Look at Figure \ref{F:indep}(a):  if the $K_4$ is balanced with edges at $v_0$ labelled $\alpha$, $\alpha$, $\beta$, the labels on $\triangle v_iv_jv_k$ are as shown.  Keeping this triangle, change the $v_0$ edges to those labelled as in Figure \ref{F:indep}(b).  The label $\epsilon$ on $\te_{ij}$ implies that $\gamma' = \gamma$.  The definition \eqref{E:labelbg} implies that $\gamma = (\alpha/\beta)_{ik}$ and $\gamma' = (\alpha/\beta)_{jk}$.  Thus \eqref{E:indep} is proved.
\begin{figure}[htbp]
\psfrag{v0}{$v_0$}
\psfrag{vi}{$v_i$}
\psfrag{vj}{$v_j$}
\psfrag{vk}{$v_k$}
\psfrag{alp}{$\alpha$}
\psfrag{bet}{$\beta$}
\psfrag{gam}{$\gamma$}
\psfrag{gamp}{$\gamma'$}
\psfrag{eps}[c]{$\epsilon$}
\psfrag{abik}[l]{$(\alpha/\beta)_{ik}$}
\psfrag{abjk}[c]{$(\alpha/\beta)_{jk}$}
\includegraphics{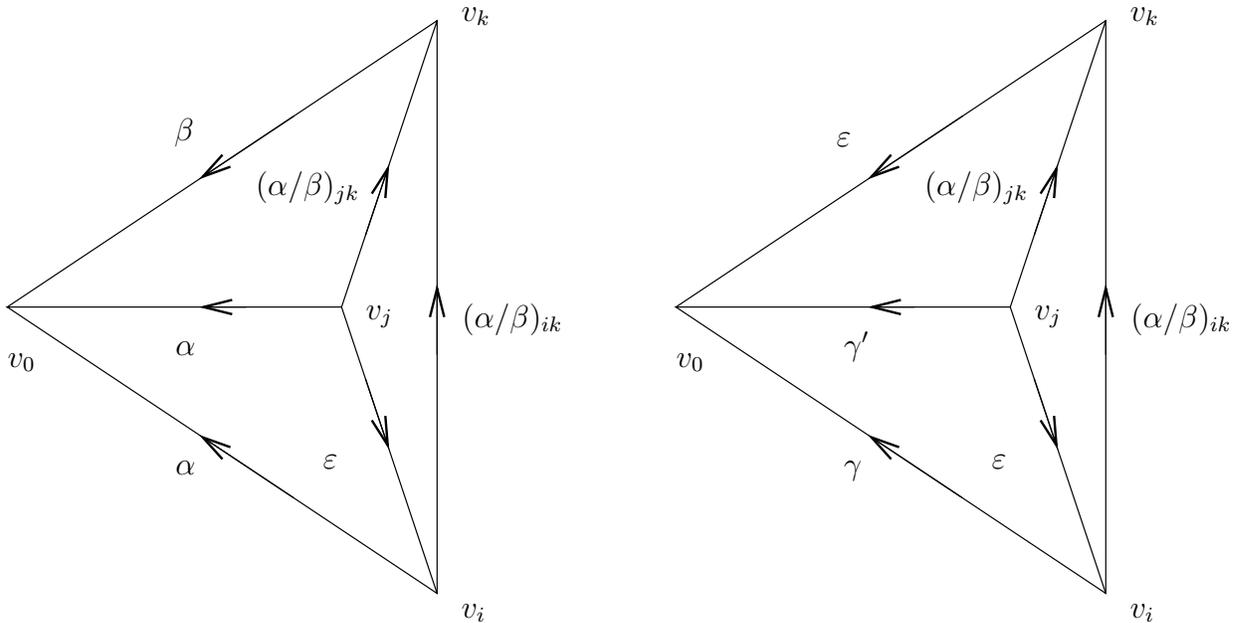}
\caption{Illustrating the proof of \eqref{E:indep}.  The graphs are balanced.}\label{F:indep}
\end{figure}

Another consequence of the definition of division is the reversal property $(\alpha/\beta)_{ij} = (\beta/\alpha)_{ji}$.  Assuming there are at least four nodes and applying \eqref{E:indep} thrice,
\begin{equation*}
(\alpha/\beta)_{ij} = (\alpha/\beta)_{kj} = (\beta/\alpha)_{jk} = (\beta/\alpha)_{ik} =
(\alpha/\beta)_{ki} = (\alpha/\beta)_{ji}.
\end{equation*}
Thus, all $(\alpha/\beta)_{ij}$ are equal:  we have a single well-defined operation $\alpha/\beta$.  

From the definitions, then:
\begin{enumerate}
\item[(L1)]  $\alpha/\alpha = \epsilon$
\item[(L2)]  $\alpha/\epsilon = \alpha$.
\end{enumerate}
By the reversal property, $(\epsilon/(\beta/\gamma))_{ij} = 
((\beta/\gamma)/\epsilon)_{ji} = (\beta/\gamma)_{ji} = (\gamma/\beta)_{ij}$, so 
\begin{enumerate}
\item[(L3)]  $\epsilon/(\beta/\gamma) = \gamma/\beta$.
\end{enumerate}
These are three of the four axioms for a group defined by division, given in \cite[p.\ 6]{Hall}.  It remains to prove that 
\begin{enumerate}
\item[(L4)]  $(\alpha/\gamma)/(\beta/\gamma) = \alpha/\beta$.
\end{enumerate}
\begin{figure}[htbp]
\psfrag{v0}{$v_0$}
\psfrag{vi}{$v_i$}
\psfrag{vj}{$v_j$}
\psfrag{vk}{$v_k$}
\psfrag{alp}{$\alpha$}
\psfrag{bet}{$\beta$}
\psfrag{gam}{$\gamma$}
\psfrag{eps}[c]{$\epsilon$}
\psfrag{ab}[l]{$\alpha/\beta$}
\psfrag{ag}[l]{$\alpha/\gamma$}
\psfrag{bg}[l]{$\beta/\gamma$}
\includegraphics{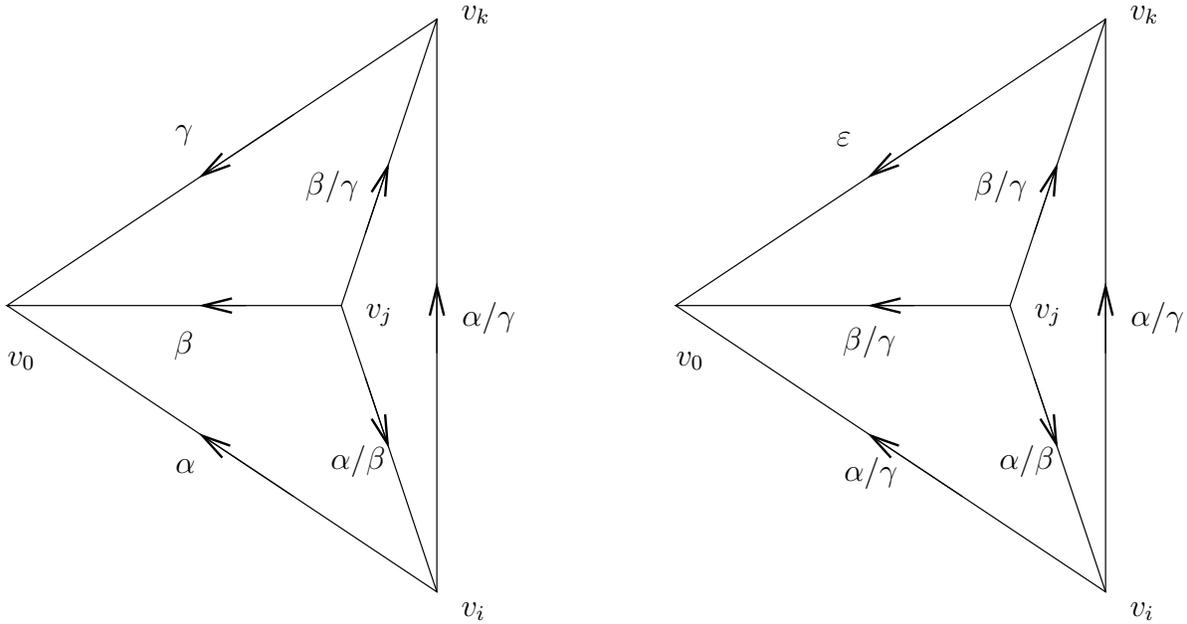}
\caption{Diagrams for the proof of property (L4).}\label{F:L4}
\end{figure}
Again, we use two diagrams:  see Figure \ref{F:L4}.  Diagram (a) is just definitions.  Holding $\triangle v_iv_jv_k$ fixed, we change the edge labels at $v_0$ so that $\te_{0k}$ has label $\epsilon$.  The labels on $\te_{0i}$ and $\te_{0j}$ are from the definition of division. Then $\te_{ij}$ has label $(\alpha/\gamma)/(\beta/\gamma)$, but we already know its label is $\alpha/\beta$.  That proves (L4).

Therefore, $\fQ$ is a group, and it is easy to verify that $\Omega = \bgr{\fQ K}$.
\end{proof}

\begin{proof}[Proof of Theorem \ref{T:3connbg}] Let $\Omega\downarrow 
\Delta$ where $\Delta$ is 3-connected.  By Example \ref{XX:parallel} we 
may assume $\Delta$ is simple.  If $v$ and $w$ are nonadjacent nodes in 
$\Delta$, they are the trivalent nodes of a theta subgraph of $\Delta$.  
By Propositions \ref{LX:theta} and \ref{TX:chord} and Theorem 
\ref{TX:max}, $\Omega$ extends to $\Omega'$, an expansion of the 
complete graph on $N(\Delta)$.  By Lemma \ref{L:completebx}, then, 
$\Omega'$ is a group expansion; hence, so is $\Omega$. \end{proof}


\section{Amalgamation (2-separation)}\label{amalg}

Biased expansions of the same multiplicity can be assembled by an analog of the ordinary graphical operation of edge amalgamation, which means combining two graphs by identifying an edge from each of them. 
This operation is essential to the structure theory of biased expansions.  
Besides that, it enables us to produce nongroup expansions out of group expansions, in two different ways.  The easy way is to combine group or quasigroup expansions by different (quasi)groups of the same order.  For instance, in multiplicity $4$ we can assemble a $\bbZ_4$-expansion and a $\fV_4$-expansion, $\fV_4$ being the Klein four-group.  A more sophisticated kind of application combines expansions by the same group but with a nasty twist.  

The first task is to define and justify the method of combination.

If a graph $\Delta$ is the union of two subgraphs, $\Delta_1$ and
$\Delta_2$, that have in common only a link $e$ and its endpoints, i.e.\
$\Delta_1 \cap \Delta_2 = (N(e),\{e\})$, we say $\Delta$ is the \emph{edge
amalgamation} (or \emph{parallel connection}) of $\Delta_1$ and $\Delta_2$ \emph{along} $e$, written $\Delta_1\cup_e \Delta_2$, 
and we call $\Delta\setm e$ the \emph{edge sum} (or \emph{$2$-sum}) of
$\Delta_1$ and $\Delta_2$ along $e$, written $\Delta_1 \oplus_e \Delta_2$.  Another way
to look at edge amalgamation or edge sum is as identification or cancellation of
distinct links $e_1 \in E(\Delta_1)$ and $e_2 \in E(\Delta_2)$.  We shall sometimes take this point of view.  

These constructions can be modelled in biased expansions.
Suppose $\Omega_1$ and $\Omega_2$ are biased expansions of $\Delta_1$ and
$\Delta_2$.  We construct an \emph{expanded edge amalgamation} of $\Omega_1$ and
$\Omega_2$ \emph{along} $e$, written $\Omega_1 \cup_e \Omega_2$ or in full
$\Omega_1\cup_{e,\beta} \Omega_2$, by choosing a bijection $\beta : p\inv_1(e) \to
p\inv_2(e)$ and using it to identify $p\inv_1(e)$ with $p\inv_2(e)$.  
The edge set of $\Omega_1 \cup_e \Omega_2$ is thus the disjoint union of $E(\Omega_1)$ and $E(\Omega_2)$ with $p\inv_1(e)$ and $p\inv_2(e)$ identified by $\beta$.  
A circle $\tilde C$ in $\Omega_1 \cup_e\Omega_2$ is balanced if it belongs to $\cB(\Omega_1) \cup \cB(\Omega_2)$ or it has the form $\tC_1 \cup \tC_2\setm [p\inv_1(e) \cup p\inv_2(e)]$ where $\tC_i \in \cB(\Omega_i)$ and $e \in p(\tC_i)$ for $i=1,2$ and $\beta(\tC_1 \cap p\inv_1(e)) = \tC_2 \cap p\inv_2(e)$.  (We may write $\tC$ more simply as $\tC_1+\tC_2$ if we bear in mind the identification of $p\inv_1(e)$ with $p\inv_2(e)$.)  
The \emph{expanded edge sum} along $e$ is $\Omega_1 \oplus_e\Omega_2 = \Omega_1 \oplus_{e,\beta} \Omega_2 = (\Omega_1\cup_e \Omega_2)\setm p\inv(e)$, $p$ being the projection mapping of $\Omega_1 \cup_e\Omega_2$.  
Both constructions apply to group expansions $\fG_1\Delta_1$ and $\fG_2\Delta_2$ by taking $\Omega_i = \bgr{\fG_i\Delta_i}$.

An example is any biased expansion $\Omega$ of $\Delta_1 \cup_e \Delta_2$.  
If $\Omega_i = p\inv(\Delta_i)$ and $\beta$ is the identity map, then $\Omega_1\cup_e \Omega_2 = \Omega$.  On the other hand, a biased expansion of $\Delta_1 \oplus_e \Delta_2$ need not be an expanded edge sum $\Omega_1 \oplus_e \Omega_2$:  see Example \ref{XD:nonsum}.

\begin{thm}\label{T:amalg}
Let $\Delta = \Delta_1 \cup_e\Delta_2$ or $\Delta_1 \oplus_e \Delta_2$, the amalgamation or sum along $e$ of graphs $\Delta_1$ and $\Delta_2$, and let $\Omega_1$ and $\Omega_2$ be biased expansions of $\Delta_1$ and $\Delta_2$ such that $\# p\inv_1(e) = \# p\inv_2(e)$.  Any expanded edge amalgamation $\Omega_1 \cup_e \Omega_2$ or expanded edge sum $\Omega_1 \oplus_e \Omega_2$ is a biased expansion of $\Delta$.  
If $\Delta$ is the edge amalgamation, then $\Omega_1$ and $\Omega_2 \subseteq \Omega_1\cup_e \Omega_2$.  
If $\Delta$ is the edge sum and $\Delta_1$ and $\Delta_2$ are 2-connected of order at least $3$, then $\Omega_1$ and $\Omega_2$ are expansion minors of $\Omega_1 \oplus_e \Omega_2$.
\end{thm}

\begin{proof} 
We show first  that $\Omega = \Omega_1\cup_e \Omega_2$ is a biased graph and
a biased expansion of $\Delta_1\cup_e\Delta_2$.  For convenience of notation we
assume that the identification prescribed by $\beta$ has been carried out.

Suppose $\tC_1\cup \tC_2$ is a theta graph in $\Omega$ and $\tC_1$ and $\tC_2$ are
balanced in $\Omega$.  We want $\tC_1 + \tC_2$ to be balanced.  If $\tC_1\cup \tC_2 \subseteq \Omega_i$, this will be so.  There are two ways $\tC_1 \cup \tC_2$ may not be in $\Omega_1$ or $\Omega_2$:  one of its three constituent paths may be an edge $\te \in p\inv(e)$, or $\tC_1 \cup \tC_2$ may be disjoint from $p\inv(e)$.
In the first case $\tC_i \setm \te$ is a path in $\Omega_i$ and $\tC_1 +
\tC_2$ is balanced by the definition of $\cB(\Omega)$.  In the second case, one
circle is contained in an $\Omega_i$, say $\tC_1\subseteq \Omega_1$; then $\tC_2$
lies partly in $\Omega_1$ and partly in $\Omega_2$.  Because $\tC_2$ is balanced,
it must be the sum $\tC'_1 + \tC'_2$ of balanced circles  $\tC'_i\subseteq \Omega_i$ that contain an edge $\te \in p\inv(e)$.  Then $\tC_1\cup \tC'_1$ is a theta graph in $\Omega_1$ and is the union of balanced circles; thus $\tC_1 +\tC'_1$ is balanced.  Hence $(\tC_1+\tC'_1) + \tC'_2$ is balanced, and this equals $\tC_1+\tC_2$.  We have proved that $\Omega$ is a biased graph.

Given a circle $C$ in $\Delta$, $f\in C$, and a lift $\tP$ of $P = C\setm f$
into $\Omega$, we want to prove there is one and only one $\tf \in p\inv(f)$ that
makes $\tP\cup \{\tf\}$ balanced.  If $C \subseteq \Delta_i$ there is nothing to prove, so we assume $C= P_1\cup P_2$ where $P_i$ is a path in $\Delta_i$ with endpoints
$N(e)$ and that $f\in P_2$.  Let $C_i = P_i\cup \{e\}$.  Then $P_1$ lifts to $\tP_1 \subseteq \tP$ and $P_2\setm f$ lifts to $\tQ_2 \subseteq \tP$.  There is a unique $\te$ for which $\tP_1 \cup \{\te\}$ is balanced. Then there is just one $\tf$ for which $\tQ_2 \cup \{\te,\tf\}$ is balanced.  Now we have two balanced circles, $\tP_1 \cup \{\te\}$ and $\tQ_2 \cup \{\te,\tf\}$, whose union is a theta graph with $\te$ as one constituent path; the other paths form a circle $\tP\cup \{\tf\}$, balanced by the definition of $\cB(\Omega)$, that projects to $C$.  Hence $\tf$ exists as desired.  Its uniqueness is obvious.

The remaining part that is not obvious is that $\Omega_1$ is a minor of $\Omega_1 \oplus_e \Omega_2$.  By 2-connectedness of $\Delta_2$, there is a circle $C$ in $\Delta_2$ that contains $e$ and an arbitrary other link $f_2$.  Let $C = eQ'_2f_2Q''_2$, lift $Q'_2\cup Q''_2$ arbitrarily to $\tQ_2$, and form the subgraph $\Omega_{f_2} = \Omega_1 \cup p\inv(f_2)\cup \tQ_2$.  We prove that $\Omega_1 \cong (\Omega_{f_2} \setm p\inv(e))/\tQ_2$ by the isomorphism $\epsilon_1$ that is the identity on $\Omega_1 \setm p\inv(e)$ and is defined on $p\inv(e)$ by
$\epsilon_1(\te) =$ that edge $\tf_2$ for which $\{\te,\tf_2\}\cup \tQ_2$ is balanced.  What has to be proved is that, for a circle $\tC \subseteq E(\Omega_1)$, $\tC$ is balanced if and only if $\epsilon_1(\tC)$ is balanced.  Let $\tC \cap p\inv(e) = \{\te\}$.  
Then $\tC\cup \tQ_2 \cup \{\epsilon_1(\te)\}$ is a theta graph in which $\tQ_2 \cup \{\te\} \cup \{\epsilon_1(\te)\}$ is balanced.  The conclusion follows.  
\end{proof}

Theorem \ref{T:amalg} allows us to produce arbitrarily large biased
expansions that are not group expansions, of any multiplicity $\gamma \geq 4$.  

\begin{exam}\label{XD:nongroup}
Let $\Delta = \Delta_1 \oplus_e \Delta_2$, where $\Delta_1$ and $\Delta_2$ are
2-connected, and let $\fG_1$ and $\fG_2$ be different groups of the same order
$\gamma$.  Form $\Omega_i = \bgr{\fG_i\Delta_i}$.  Any bijection $\fG_1\to
\fG_2$ induces a bijection $\beta : p\inv_1(e) \to p\inv_2(e)$ by which we can
form an expanded edge sum $\Omega_1 \oplus_e \Omega_2$.  The sum has as minors both
$\bgr{\fG_1\Delta_1}$ and $\bgr{\fG_2\Delta_2}$, and these in
turn have minors $\bgr{\fG_1K_3}$ and $\bgr{\fG_2K_3}$.  If
$\Omega = \bgr{\fG\Delta}$, then all triangular minors are isomorphic to
$\bgr{\fG K_3}$, but this is impossible.  Therefore $\Omega$ is a nongroup
regular biased expansion of $\Delta$.  Note that this construction cannot be
carried out for prime multiplicities $\gamma$.
\end{exam}

\begin{exam}\label{XD:qnongroup}
In the preceding construction take $\Delta_2 = K_3$ and let $\Omega_2$ be any
quasigroup expansion of $K_3$ having multiplicity $\gamma$ but not isomorphic
to $\bgr{\fG_1K_3}$.  Then $\Delta = \Delta_1 \oplus_e K_3$ has a regular
biased expansion with nonisomorphic triangular minors $\bgr{\fG\Delta_1}$
and $\Omega_2$, so it is a nongroup regular biased expansion of $\Delta$.  This
construction can be carried out for all multiplicities $\gamma\geq 4$.
\end{exam}

The technique of summing with quasigroup expansions of a triangle yields highly nongainable biased expansions of series-parallel graphs, just to mention a sizeable class to which it applies.  The reason is that every series-parallel graph $\Delta$ is constructed by doubling edges in parallel, an operation that is trivial to reproduce in a biased expansion of $\Delta$ (see Example \ref{XX:parallel}), and by subdividing edges, which is equivalent to taking an edge sum with a triangle.  
On the other hand, the methods of Example \ref{XD:nongroup} and \ref{XD:qnongroup} together still do not give ungainable biased expansions with all multiplicities $\gamma\geq4$ of all 2-separable inseparable graphs.  For that see Corollary \ref{C:basenongroup}.

\begin{exam}\label{XD:nonsum}
Here is an example of a biased expansion of $\Delta_1 \oplus_e \Delta_2$ that is not an expanded edge sum of expansions of $\Delta_1$ and $\Delta_2$. In the example, $\Delta_1 \cong \Delta_2 \cong K_3$.

Take $C_4 = (N,E)$ where $N = \{v_1,v_2,v_3,v_4\}$ and 
$E = \{ e_{12},e_{23},e_{34},e_{41} \}$.  $C_4$ is an edge sum in two different ways:  it is $\Delta_{123} \oplus_{e_{13}} \Delta_{134}$ and 
$\Delta_{124} \oplus_{e_{24}} \Delta_{234}$.  Here $\Delta_{ijk}$ denotes the triangle with node set $\{v_i,v_j,v_k\}$.  Let $\gamma \geq 4$ and let $\gamma\cdot \Delta_{123}$ and $\gamma\cdot \Delta_{134}$ be biased expansions that are not both group expansions by the same group.  (That is, one or both is not a group expansion, or $\gamma\cdot \Delta_{123} = \bgr{\fG\Delta_{123}}$ and $\gamma\cdot \Delta_{134} = \bgr{\fH\Delta_{134}}$ where $\fG \not\cong \fH$.)  Then $\Omega = (\gamma\cdot \Delta_{123}) \oplus_{e_{13}} (\gamma\cdot \Delta_{134})$ is a biased expansion of $C_4$; also, 
$\Omega_{13} = (\gamma \cdot \Delta_{123}) \cup_{e_{13}} (\gamma\cdot \Delta_{134})$ is a biased expansion of $K_4 \setm e_{24}$.  

Although $C_4 = \Delta_{124} \oplus_{e_{24}} \Delta_{234}$, $\Omega$ cannot be an expanded edge sum of the form $(\gamma\cdot\Delta_{124}) \oplus_{e_{24}} (\gamma\cdot\Delta_{234})$.  We prove this by contradiction.  
Suppose it were; then (as we shall demonstrate) the union of $\Omega_{13}$ with $\Omega_{24} = (\gamma\cdot\Delta_{124})\cup_{e_{24}} (\gamma\cdot \Delta_{234})$ would be a $\gamma\cdot K_4$ having as minors both $\gamma\cdot \Delta_{123}$ and $\gamma\cdot \Delta_{134}$.
These are not isomorphic, but by Lemma \ref{L:completebx} $\gamma\cdot K_4$ is a group expansion and therefore all its triangular minors are isomorphic.  We have a contradiction.  

To prove that $\Omega_{13}$ and $\Omega_{24}$ can be contained in a biased expansion $\Omega_4$ of $K_4$ we need to define balance of circles not contained in either part of the union.  Such a circle has to be a quadrilateral that contains edges $\te_{13}$ and $\te_{24}$; let us say it is $\tQ = \te_{13}\te_{34}\te_{24}\te_{12}$.  
The decision about balance of $\tQ$ is made by representing it as the sum of circles in $\cC(\Omega_{13}) \cup \cC(\Omega_{24})$ of which one or two are balanced: in the latter case $\tQ$ is balanced, in the former case it is not.  
We find such circles by taking a chord of $\tQ$, either
$\te_{23}$ or $\te_{14}$, that makes at least one balanced triangle with $\tQ$.
There will be a biased graph $\Omega$ if and only if all the ways of deducing
balance of $\tQ$ yield the same conclusion.  There are essentially two ways to
choose the chord.  Suppose $\te_{23}$ chosen to make $\te_{12}\te_{13}\te_{23}$
balanced and $\te'_{23}$ chosen to make $\te_{24}\te_{34}\te'_{23}$ balanced. If
$\te_{23} = \te'_{23}$, $\tQ$ must be balanced.  If $\te_{23} \neq \te'_{23}$,
$\tQ$ must not be balanced.  Similarly, choose $\te_{14}$ and $\te'_{14}$ so that
$\te_{12} \te_{24}\te_{14}$ and $\te_{13}\te_{34}\te'_{14}$ are balanced.  $\tQ$
should be balanced if and only if $\te_{14} =\te'_{14}$.  

To show that $\Omega_4$ is well defined it will suffice to prove that $\te_{23}
=\te'_{23}$ implies that $\te_{13}\te_{34}\te_{14}$ is balanced.  The two balanced
triangles we are assuming in $\Omega_{24}$ imply that
$\te_{12}\te_{23}\te_{34}\te_{14}$ is balanced.  This balanced quadrilateral and
the balanced triangle $\te_{12}\te_{13}\te_{23}$ in $\Omega_{13}$ imply that
$\te_{13}\te_{34}\te_{14}$ is balanced.  Thus, our two criteria are consistent.

Note that we have no choice in deciding whether $\tQ$ is balanced; we can only have a consistent result or an inconsistency.  As there is no inconsistency,
$\Omega_4$ is well defined and uniquely defined.  Consequently, by the previous
argument, we have contradicted the hypothesis that 
$\Omega = (\gamma\cdot\Delta_{124}) \oplus_{e_{24}} (\gamma \cdot \Delta_{234})$. 
\end{exam} 

One wants to know that a multiple edge amalgamation is independent of the order of amalgamation.  It suffices to treat two amalgamations.

\begin{thm}\label{T:amalgorder}  
Let $\Omega_i \downarrow \Delta_i$ for $i = 1,2,3$, where
\begin{equation*}
E(\Delta_1 \cap \Delta_2) = \{e\}, 
\quad E(\Delta_2\cap \Delta_3) = \{f\}, 
\quad E(\Delta_1\cap \Delta_3) = \{e\} \cap \{f\},
\end{equation*}
and $\Delta_1$, $\Delta_2$, $\Delta_3$ are pairwise node-disjoint except as required by shared edges.  Suppose given bijections $\alpha : p_1\inv(e)\to p_2\inv(e)$ and $\beta : p_3\inv(f) \to p_2\inv(f)$.  Then 
\begin{equation}\label{E:amalgorder}
\Omega_1 \cup_{e,\alpha} (\Omega_2 \cup_{f,\beta} \Omega_3) = (\Omega_1 \cup_{e,\alpha} \Omega_2) \cup_{f,\alpha} \Omega_3.
\end{equation}
\end{thm}

\begin{proof}  
The only question is the balance of circles in the amalgamation.  
Let $\Omega_L$ and $\Omega_R$ be the biased expansions of 
$\Delta_1\cup_e \Delta_2 \cup_f \Delta_3$ on the left and right sides of \eqref{E:amalgorder}.  Consider a circle $\tC$ that meets both 
$\Omega_1\setm p_1\inv(e)$ and $\Omega_3\setm p_3\inv(f)$; 
thus $\tC$ consists of $\tP_1 = \tC \cap E(\Omega_1)$, $\tP_3 = \tC\cap E(\Omega_3)$ (both of which are paths) and $\tQ = \tC \cap E(\Omega_2)$. The latter may consist of two, one, or (if $N(e) = N(f))$ no paths.

We may use $\alpha$ and $\beta$ to identify $p_1\inv(e)$ with $p_2\inv(e)$ and $p_2\inv(f)$ with $p_3\inv(f)$.  

The case $e = f$ is trivial, since we are really looking at $\Omega_1 \cup_{e,\alpha\circ\beta\inv} \Omega_3$ in both $\Omega_L$ and $\Omega_R$.  

When $e\neq f$, choose $\te$ and $\tf$ so $\tP_1\cup \te$ and $\tP_3 \cup \tf$ 
are balanced in $\Omega_1$ and $\Omega_3$, respectively.  Then $\tC$ is balanced in $\Omega_L\iff \tP_3 \cup \tQ \cup \te$ is balanced in 
$\Omega_2\cup_{f,\beta} \Omega_3$ (because $\tP_1\cup \te$ is balanced) $\iff \tQ \cup \{\te, \tf\}$ is balanced in $\Omega_2$ (because $P_3 \cup \tf$ is balanced).  Similarly, $\tC$ is balanced in $\Omega_R \iff \tQ\cup \{\te,\tf\}$ is balanced in $\Omega_2$.  It follows that balance of $\tC$ is the same in $\Omega_L$ and $\Omega_R$.  
\end{proof}

The theorem implies that one can define a multiple expanded edge amalgamation directly, even one with an infinite number of amalgamations, because defining balance of any particular circle $\tC$ in the result only involves a finite number of amalgamations, so is order independent by Theorem \ref{T:amalgorder} and induction.  For more on the definition of multiple amalgamation see after Theorem \ref{T:structuremax}.

\bigskip

When we amalgamate two group expansions, whether we get a group expansion or not
depends on the nature of the identification function $\beta$.  As a mapping
$p\inv_1(e) \to p\inv_2(e)$, $\beta$ induces a mapping of groups by
composition with the gain functions, namely
\begin{equation*}
\bbeta = \big[\phi_1\big|_{p\inv_1(e)}\big]\inv \circ \beta \circ
\big[\phi_2\big|_{p\inv_2(e)}\big] : \fG_1\to \fG_2,
\end{equation*}
or in a more compact expression, 
\begin{equation}\label{ED:betadef}
\bbeta(\phi_1(\te)) = \phi_2(\beta(\te)) \qquad\text{for}\quad \te \in p\inv_1(e).
\end{equation}
Note that $\bbeta$ depends on the choice of gains; if we used different gain functions $\phi'_i$ we would get a different bijection $\bbeta'$.

We shall need to know the effect on $\bbeta$ of switchings $\eta_i$ and group
automorphisms $\alpha_i$ applied to $\Phi_1$ and $\Phi_2$.  We write $\phi'_i =
\phi_i^{\eta_i\alpha_i}$. The definition of $\bbeta'$, in full, is
\begin{align*}
\phi'_2(\beta(\te)) &= [\eta_2(v)\inv\phi_2(\beta(\te))\eta_2(w)]^{\alpha_2}\\
&=\eta_2(v)^{-\alpha_2}\bbeta(\phi_1(\te))^{\alpha_2}\eta_2(w)^{\alpha_2}.
\end{align*}
Since $\phi_1^{\eta_1\alpha_1}(\te) = [\eta_1(v)\inv\phi_1(\te)\eta_1(w)]^{\alpha_1}$, we can substitute
$$
\phi_1(\te) = \eta_1(v)[\phi'_1(\te)]^{\alpha_1\inv}\eta_1(w)\inv
$$
in the previous equation, getting 
$$
\bbeta'(\phi'_1(\te)) = \phi'_2(\beta(\te)) =
\eta_2(v)^{-\alpha_2}\bbeta[\eta_1(v)\phi'_1(\te)^{\alpha_1\inv}\eta_1(w)\inv]^{\alpha_2}\eta_2(w)^{\alpha_2}.
$$
Here $\phi'_1(\te)$ can be any group element; therefore we can rewrite the equation as
\begin{equation}\label{ED:betasw}
\bbeta'(g) = \eta_2(v)^{-\alpha_2} \bbeta[\eta_1(v)g^{\alpha_1\inv}\eta_1(w)\inv]^{\alpha_2} \eta_2(w)^{\alpha_2}.
\end{equation}

A \emph{pseudoisomorphism} of groups (or quasigroups) is any mapping 
$\fG_1\to\fG_2$ that has the form $g \mapsto g^\alpha c$ where $\alpha: 
\fG_1\to\fG_2$ is an isomorphism and $c\in\fG_2$.  The (quasi)groups must 
be isomorphic for such a mapping to exist.  The pseudoautomorphisms of a 
group form a group, which we denotate $\Aff\fG$.

\begin{thm}\label{T:twist}
Let $\Delta = \Delta_1\cup_e\Delta_2$, where $\Delta_1$ and $\Delta_2$ are
2-connected simple graphs, and let $\Omega_1$ and $\Omega_2$ be biased expansions
of $\Delta_1$ and $\Delta_2$ with the same multiplicity.  The expanded edge
amalgamation $\Omega = \Omega_1\cup_{e,\beta} \Omega_2$ and the expanded edge sum
$\Omega_0 = \Omega_1 \oplus_{e,\beta} \Omega_2$ are group expansions (of
$\Delta$ and $\Delta \setm e$, respectively) if and only if $\Omega_1 =
\bgr{\fG_1\Delta_1}$ and $\Omega_2 = \bgr{\fG_2\Delta_2}$ and $\beta$, after suitable switching of $\fG_1\Delta_1$ and
$\fG_2\Delta_2$, induces an isomorphism $\fG_1\to \fG_2$.

If $\fG_1 \cong \fG_2$, the condition on $\beta$ is equivalent to 
$\bbeta$'s being a pseudoisomorphism $\fG_1\to \fG_2$.
 \end{thm}

What we mean by suitable switching is that there exist switching functions 
$\theta_1$ and $\theta_2$ such that when $\beta$ is applied to 
$(\fG_1\Delta_1)^{\theta_1}$ and $(\fG_2\Delta_2)^{\theta_2}$, then 
$\beta$ induces an isomorphism $\fG_1\to \fG_2$.  In terms of the 
original, unswitched gains, the induced mapping is 
$\big[\phi^{\theta_1}_1\big|_{p\inv_1(e)}\big]\inv \circ \beta \circ 
\big[\phi_2^{\theta_2}\big|_{p\inv_2(e)}\big]$. We call $\beta$ 
\emph{twisted} if no such switchings exist, or equivalently if $\bbeta$ is 
not a pseudoisomorphism, and in particular if the groups are not 
isomorphic in the first place.

\begin{proof} 
 This is one of those theorems that seem obvious but have a complicated 
proof.  The beginning is easy:  according to Theorem \ref{T:amalg} the 
expanded edge amalgamation or sum can be a group expansion only if 
$\Omega_1 = \bgr{\fG\Delta_1}$ and $\Omega_2 =\bgr{\fG\Delta_2}$ for some 
group.  Let us therefore assume this is so and write $\Phi_i = 
\fG\Delta_i$.  What we need to prove is the equivalence of the following 
properties:
 \begin{enumerate}
 \item[(i)]\label{twist0}
 $\Omega_0$ is a $\fG$-expansion of $\Delta_0 = \Delta\setm e$: that is, 
$\Omega_0\cong \bgr{\fG\Delta_0}$.  
 \item[(ii)]\label{twistb}
 $\bbeta(g) = g^\alpha c$ for some $c \in \fG$ and $\alpha \in \Aut\fG$.
 \item[(iii)]\label{twists}
 There are switchings $(\fG\Delta_1)^{\theta_1}$ and
$(\fG\Delta_2)^{\theta_2}$ such that $$\alpha =
\big[\phi_1^{\theta_1}\big|_{p\inv_1(e)}\big]\inv \circ \beta\circ
\big[\phi_2^{\theta_2}\big|_{p\inv_2(e)}\big]$$ is an automorphism of $\fG$.
 \item[(iv)]\label{twist}
 $\Omega\cong \bgr{\fG\Delta}$.
 \end{enumerate}
 We show that (i)$\implies$(ii)$\implies$(iii)$\implies$(iv).  

Assume, then, that $\Omega_0 \cong \bgr{\Phi_0}$ where $\Phi_0 = \fG\Delta_0$, and also, by prior switchings $\eta'_i$ of $\Phi_i$ to $\Phi'_i$ (for $i=1,2$), that
\begin{equation}\label{ED:betanorm}
\phi'_1(1e) = 1 \quad\text{ and }\quad \phi'_2(\beta(1e)) = 1.
\end{equation}
We may choose $\eta'_1\equiv 1$ and $\eta'_2(v)=1$.  
Now, if we take paths $\tP_i$ in $\Phi_i$ such that $\tP_1\cup \{1e\}$ and 
$\tP_2 \cup \{\beta(1e)\}$ are balanced circles, then 
$\phi'_1(\tP_1) = \phi'_1(1e) =1$ and $\phi'_2(\tP_2) =\phi'_2(\beta(1e)) = 1$.  
(Here we orient $1e$, $\beta(1e)$, $\tP_1$, and $\tP_2$ similarly, from one endpoint $v$ of $e$ to the other endpoint $w$.)  
By construction, $\tP_1\cup \tP_2$ is balanced in $\Omega_0$; thus 
$\phi'_0(\tP_1) = \phi'_0(\tP_2)$; consequently, we may assume by prior switching of $\Phi_0$ that $\phi_0\big|_{\tP_1} \equiv 1$ and $\phi_0\big|_{\tP_2} \equiv 1$.  Note, though, that $\phi_0$ need not agree with $\phi'_1$ even though 
$\Omega_1\setm p\inv_1(e)\subseteq \Omega_0$, and the same for $\phi'_2$.

Nevertheless, $\Omega_1$ is isomorphic to a minor $\Omega_{01}$ of $\Omega_0$ that
can be found by following the proof of Theorem \ref{T:amalg}.  In that proof
choose $\tf_2 \in \tP_2$ and $\tQ_2 =\tP_2\setm \tf_2$.  The proof constructs
$\Omega_{01}$ with underlying graph $\|\Omega_0\| =\|\Omega_1\|\setm
p\inv_1(e) \cup p\inv_2(f_2)$.  Because $\tP_2$ has all identity gains, the
correspondence $\epsilon_1$ preserves gains.
Therefore, $\Omega_{01} = \bgr{\Phi_{01}}$ where $\Phi_{01}$ is a minor of $\Phi_0$ with gains $\phi_{01} = \phi_0\big|_{E_{01}}$.

Since $\bgr{\Phi_1} \cong \bgr{\Phi_{01}}$, by uniqueness of gains \cite[\uniquegains(c)]{BG} 
$\phi_{01} =(\phi'_1)^{\eta_1\alpha_1}\circ\epsilon_1$,
where $\eta_1$ is a switching function and $\alpha_1 \in \Aut\fG$.  Without
loss of generality we may assume that $\eta_1(v) = 1$.  Then 
$$
1 = \phi'_0(\tP_1) = (\phi'_1)^{\eta_1\alpha_1} (\tP_1) =
[\eta_1(v)\inv\phi'_1(\tP_1)\eta_1(w)]^{\alpha_1} 
= [1\cdot 1\cdot \eta_1(w)]^{\alpha_1}
$$
implies $\eta_1(w) = 1$.

Similarly we construct $\Omega_{02} =\bgr{\Phi_{02}}$, a minor of $\Omega_0$ that
is isomorphic to $\Omega_2$ with $\phi_{02} = (\phi'_2)^{\eta_2\alpha_2}\circ
\epsilon_2$ where $\eta_2(v) = \eta_2(w) = 1$. 

Applying Equation \eqref{ED:betasw} to the special circumstances of $\Phi^{\eta_i\alpha_i}_i$ where $\eta_i(v) = \eta_i(w) = 1$, we see that 
$$
\bbeta'' = \alpha\inv_1\circ\bbeta'\circ \alpha_2.
$$  
Consequently, $\bbeta'\in \Aut\fG \iff \bbeta'' \in \Aut\fG.$

The next step is to prove the (surprising) fact that $\bbeta''$ is the identity.
If $\tP_1\cap p\inv_1(f_1) = \{\tf_1\}$, and if we write $\tf^*_1 =
\epsilon_2(\te)$ for $\te\in p\inv_1(e)$ and $\tP^*_1 = \tP_1\setm
\{\tf_1\}\cup \{\tf^*_1\}$, then $\tP^*_1\cup \{\te\}$ is balanced.  Similarly,
$\tP^*_2\cup \{\beta(\te)\}$ is balanced, so $\tP^*_1\cup \tP^*_2$ is also
balanced.  It follows, from balance of each of these circles in turn, that 
\begin{gather*}
\phi''_1(\te) = \phi_{01}(\tf^*_1) = \phi_0(\tf^*_1),\\
\phi''_2(\beta(\te)) =\phi_{02}(\tf^*_2) = \phi_0(\tf^*_2),\\
\phi_0(\tf^*_1) = \phi_0(\tP^*_1) = \phi_0(\tP^*_2) =\phi_0(\tf^*_2),
\end{gather*}
where $\phi''_i = (\phi'_i)^{\eta_i\alpha_i}$.
Hence, $\phi''_1(\te) = \phi''_2(\beta(\te))$.  This means that $\bbeta''$ is the identity mapping.  

Therefore $\bbeta'$ is an automorphism of $\fG$; in fact, $\bbeta' =
\alpha_1\circ \alpha\inv_2$. 

The course of the proof so far may be summarized in a diagram.  In it, $\teta_i$ is the fibered permutation of $E_i$ induced by $\eta_i$; that is, $(\phi'_i)^{\eta_i} = \teta_i \circ \phi'_i$.
(In the description and diagram $\phi'_i$, $\teta_i$, etc.\ stand for $\phi'_i\big|_{p\inv_i(e)}$, etc.; so that all maps are bijections.) 
The first square is commutative because $\teta_i$ is the identity on $p\inv_i(e)$, a consequence of having $\eta_i(v) = \eta_i(w) = 1$.  
The triangles commute by the definition of $\teta_i$.  The square $\phi'_1\circ\bbeta'$ vs.\ $\beta \circ \phi'_2$ commutes
by the definition of $\bbeta'$ and the rectangle
$(\phi'_1)^{\eta_1}\circ\alpha_1\circ\bbeta''$ vs.\
$\beta\circ(\phi'_2)^{\eta_2}\circ\alpha_2$ commutes by the definition of
$\bbeta''$.  From this it follows that the entire diagram commutes; since
$\bbeta''\in \Aut\fG$, then $\bbeta'\in \Aut\fG$. 
$$
\xymatrix{
p\inv_1(e) \ar@/^1.5pc/[rr]^{(\phi'_1)^{\eta_1}_1} \ar@{->}[r]_-{\teta_{1}} \ar@{->}[d]_-\beta
&p\inv_1(e)\ar@{->}[r]_-{\phi'_1} \ar@{-->}[d]_-\beta &\fG
\ar@{->}[r]_-{\alpha_1} \ar@{-->}[d]_-{\bbeta'} &\fG\ar@{-->}[d]_-{\bbeta''}\\
p\inv_2(e) \ar@/_1.5pc/[rr]_{(\phi'_2)^{\eta_2}}
\ar@{->}[r]^-{\teta_{2}} &p\inv_2(e)\ar@{->}[r]^-{\phi'_2}
&\fG\ar@{->}[r]^-{\alpha_2} &\fG
}
$$

The reason $\bbeta'$ is an automorphism is that we did the right kind of switching.
First we switched $\Phi_i$ by $\eta'_i$ so that $\phi'_1(1e) = 1$ and $\phi'_2(\beta(1e)) = 1$, then we switched $\Phi'_i$ by $\eta_i$.  
The overall effect is that of switching $\Phi_1$ by $\theta_1=\eta_1$ and
$\Phi_2$ by $\theta_2=\eta'_2\eta_2$.  
We also switched $\Phi_0$, but that is unimportant because the gains on $\Phi_0$ were not given in advance like those on $\Phi_1 = \fG\Delta_1$ and $\Phi_2 =\fG\Delta_2$.  

Expressed in terms of the original gains $\phi_i$, the definition of $\bbeta'$ is
$\bbeta'(\phi_1(\te)) = \phi_2^{\eta'_2}(\beta(\te)).$
Substituting the values of $\eta'_2(v)$ and $\eta'_2(w)$, this becomes 
\begin{equation*} \bbeta'(\phi_1(\te)) = \phi_2(\beta(\te))\eta'_2(w) = \bbeta(\phi_1(\te))\eta'_2(w).
\end{equation*}
Setting $\te=ge$, we see that $\bbeta(g) = g^\alpha c$ for $\alpha = \bbeta' \in \Aut\fG$ and $c = \eta'_2(w)\inv \in \fG$, thereby proving (ii) from (i).

We know (ii)$\implies$(iii) because we can produce the necessary switching functions: $\theta_1\equiv1$ for $\Phi_1$ and $\theta_2$ with $\theta_2(v)=1$ and $\theta_2(w)=c\inv$ for $\Phi_2$.

Proving (iii)$\implies$(iv) is easy.  We may assume $\fG\Delta_1$ and
$\fG\Delta_2$ switched and $\alpha$ previously applied to $\Phi_1$ so that, in effect, $\alpha$ becomes the identity.  Then 
$\beta : p_1\inv(e)=\fG\times\{e\} \to p_2\inv(e)=\fG\times\{e\}$ 
is the identity, so the amalgamation is
$\bgr{\fG_1\Delta_1} \cup_{e,\text{id}} \bgr{\fG_2\Delta_2}$, which is simply $\bgr{\fG\Delta}$.
\end{proof}

Theorem \ref{T:twist} helps answer some questions about the existence of biased expansions that do not have gains.  One question is whether an expanded edge amalgamation or sum of two $\fG$-expansions is itself a group expansion.  That depends in part on whether or not $\Aut\fG$ is the full symmetric group of $\fG\setm \{1\}$.

\begin{lem}[{\cite[Corollary V.\{CD:fiberaut\}]{BG}}]\label{C:fiberaut} 
Assuming $\Delta$ is a block of order at least 3, $\Aut_p\bgr{\fG\Delta}$ acts as the symmetric group on a fiber $p\inv(e)$ if and only if $\fG = \bbZ_\gamma$ for $\gamma \leq 3$ or $\fG = \fV_4$.  
\end{lem}

\begin{cor}\label{C:gamalg}
Suppose $\Delta_1$ and $\Delta_2$ are 2-connected simple graphs.  An expanded edge amalgamation or expanded edge sum of group expansions $\fG\Delta_1$ and $\fG\Delta_2$ is necessarily a group expansion if and only if $\fG = \bbZ_\gamma$ for $\gamma \leq 3$ or $\fG = \fV_4$; and then it is a $\fG$-expansion. 
\end{cor}

\begin{proof}  
This is immediate from Lemma \ref{C:fiberaut}, which tells us that it is possible to find a bijection $\beta$ for which $\bbeta$, after suitable switching, is still not an automorphism if and only if $\fG$ is any group other than $\bbZ_\gamma$, $\gamma\leq 3$, and $\fV_4$.
\end{proof}

The application to multary quasigroups is Corollary \ref{C:gqgroup}.

Another question resolved by Theorem \ref{T:twist} is whether it might
be possible to ensure that an edge amalgamation or sum is a group expansion by
putting a restriction on triangular expansion minors.  For any group expansion, all expansion minors are expansions by the same group (Proposition
\ref{P:expminor}).  We might conjecture a kind of converse:  that 
$\Omega = \bgr{\fG\Delta_1} \cup_e\bgr{\fG\Delta_2}$ or 
$\bgr{\fG\Delta_1} \oplus_e \bgr{\fG\Delta_2}$ is a $\fG$-expansion if every triangular expansion minor is isomorphic to $\bgr{\fG K_3}$.  However, in general this is false.

\begin{cor}\label{C:gminor}
It is possible to have a biased expansion $\gamma\cdot\Delta$, where $\Delta$ is a 2-connected but 2-separable simple graph, such that every triangular expansion minor is isomorphic to $\bgr{\fG K_3}$ for a fixed group $\fG$ but $\gamma\cdot\Delta$ is not a group expansion, except when $\fG = \bbZ_\gamma$ for $\gamma \leq 3$ or $\fG = \fV_4$.  
Furthermore, $\gamma\cdot\Delta$ can be taken to be an edge amalgamation of group expansions.
\end{cor}

\begin{lem}\label{L:gminor}
Let $\fG$ be a fixed group.  Suppose $\Omega$, a biased expanison of a 2-connected graph $\Delta$, is obtained by expanded edge summations and amalgamations from $\fG$-expansions of inseparable graphs.  Then every triangular expansion minor of $\Omega$ is isomorphic to $\bgr{\fG K_3}$.
\end{lem}

\begin{proof}
We use induction on the order of $\Omega$.  Suppose in the construction of $\Omega$ that the last step is to assemble $\Omega_1\downarrow\Delta_1$ and $\Omega_2\downarrow\Delta_2$ into $\Omega = \Omega_1 \cup_{e,\beta} \Omega_2$ (or $\Omega = \Omega_1 \oplus_{e,\beta} \Omega_2$, but it suffices to consider the case of amalgamation).  
Consider a triangular expansion minor $\Omega_3$ of $\Omega$ whose edge set is $E_3 = p\inv(\{e_1,e_2,e_3\})$.  Let $\Omega_3'$ be the corresponding subgraph of $\Omega$; that is, the subgraph induced by the edge set $E_3$.  
If $e_1,e_2,e_3 \in E(\Delta_i)$, then $\Omega_3$ is an expansion minor of $\bgr{\fG\Delta_i}$ and the desired conclusion follows from Proposition \ref{P:expminor}.  
Otherwise, we may assume $e_1,e_2\in E(\Delta_1\setm e)$ and $e_3\in E(\Delta_2\setm e)$.  By the definition of a minor, there is a circle $C$ in $\Delta$ that contains all three edges such that $R = C \setm \{e_1,e_2,e_3\}$ has a lift $\tR$ for which $\Omega_3 = (\Omega_3' \cup \tR)/\tR$.  
Let $E_{30} = E_3 \cup p\inv(e)$, let $\Omega_{30}'$ be the corresponding subgraph of $\Omega$, and let $\Omega_{30} = (\Omega_{30}' \cup \tR)/\tR$.  
Then $\Omega_{30}$ is a biased expansion of the graph $\Delta_{30}$ consisting of the triangle $\{e_1,e_2,e_3\}$ and an edge $e$ parallel to $e_3$.  By Example \ref{XX:parallel}, $\Omega_3 = \Omega_{30} \setm p_{30}\inv(e)$ is isomorphic to $\Omega_{30} \setm p_{30}\inv(e_3)$.  The latter is an expansion minor of $\Omega_1$, hence isomorphic to $\bgr{\fG K_3}$.
\end{proof}

\begin{proof}[Proof of Corollary \ref{C:gminor}.]
The exceptional cases are covered by Corollary \ref{C:gamalg}.
For other groups, by Corollary \ref{C:gamalg} $\Omega$ need not be a group
expansion.  However, by the lemma, every triangular expansion minor is a $\fG$-expansion.
\end{proof}

Corollary \ref{C:gminor} might suggest that it is difficult to say from a criterion based on small minors whether $\Omega$ is or is not a group expansion.  But that is not correct; minors of order four suffice; see Theorem \ref{T:gminor4}.

\bigskip

A question that is not answered so far is that of reducibility of arbitrary biased expansions of 2-connected, 2-separable graphs.  
The methods of Theorems \ref{T:amalg} and \ref{T:twist} produce only nongroup expansions that are 2-separable and have a 2-separation whose nodes are adjacent or can be made adjacent in an extended biased expansion.  They will not give an example in which no 2-separating node pairs can be made adjacent: for instance, a biased expansion $4\cdot C_4$ in which it is not possible to add a chord of the $C_4$.  
Any irreducible $n$-ary quasigroup $\fQ$ with $n \geq 3$ provides such an example in the form of the expansion $\fQ C_{n+1}$.  By the results of Section \ref{bx}, that is the only way.

\bigskip

We want criteria to decide when a biased expansion of a 2-separable graph $\Delta$ is an expanded edge amalgamation or sum along an edge (not necessarily in $\Delta$) whose endpoints separate $\Delta$.  

\begin{cor}[Test for Decomposability across a 2-Separation]\label{C:test}
Suppose $\Omega\downarrow\Delta$, where $\Delta$ is 2-connected, and $\{v,w\}$ is a 2-separation of $\Delta$ into subgraphs $\Delta_1$ and $\Delta_2$.  Let $\Omega_i=\Omega\big|_{\Delta_i}$.  
If $v$ and $w$ are adjacent by an edge $e_{vw}$, then $\Omega$ has the form of an expanded edge amalgamation $\Omega_1 \cup_{e_{vw}} \Omega_2$.  
If they are not adjacent, choose an arbitrary circle $C \subseteq E(\Delta)$ through $v$ and $w$.  
Then $\Omega$ is an expanded edge sum $\Omega_1' \oplus_{e_{vw}} \Omega_2'$ if and only if $\Omega\big|_C$ extends to $e_{vw}$.
\end{cor}

\begin{proof}
The first part is obvious.  In the second part, if $\Omega$ is an edge 
sum, then it extends to $\Omega' \downarrow \Delta \cup e_{vw}$, formed by 
amalgamating instead of summing.  Conversely, if $\Omega\big|_C$ extends 
to $e_{vw}$, then $\Omega$ extends, by Proposition \ref{TX:chord}, and 
therefore is an expanded edge sum. \end{proof}

Belousov and Sandik have a criterion for extendibility of $\Omega\big|_C$ to a chord $e_{vw}$, expressed in terms of factorizability of a multary quasigroup (which is equivalent by Proposition \ref{P:qgroup}).  Let $P$ and $Q$ be the paths into which $v$ and $w$ divide $C$.  Translated to biased expansions, the criterion says:

\begin{prop}[{\cite[Lemma 6]{BelSand}}]\label{P:sandik}
If there exist lifts $\tP$, $\tP^*$, $\tQ$, and $\tQ^*$ such that $\tP\cup\tQ$, $\tP^*\cup\tQ$, $\tP\cup\tQ^*$ but $\tP^*\cup\tQ^*$ is not, then $\Omega\big|_C$ does not extend.  Otherwise, it extends.
\end{prop}


\section{The structure of biased expansions}\label{bx}

We have two main structure theorems.  One is about maximal biased expansions, and translates directly into a structural description of multary quasigroups (Theorem \ref{T:quasidecomp}).  The other describes all biased expansion graphs.  We want to make it very clear that these theorems are proved only for expansions of base graphs that are \emph{2-connected} and have \emph{finite order}.  
The former is an insignificant restriction in general:  when expanding an arbitrary graph, the expansion of each block is unrelated to that of any other block, so it is inevitable that a theorem can only refer to 2-connected graphs (but for regular expansions see Proposition \ref{T:sepmax}).  
The restriction to finite order is due to the absence of a 3-decomposition theory of infinite graphs.  (I see no reason why such a theory should not exist.)  Another necessity for our structural theorems is a Menger theorem for 2-separation of nodes in infinite graphs; but as this was proved for denumerably infinite graphs by Erd\H{o}s (see \cite[Ch.\ XIV, \S4]{Konig}) it is not an obstacle in that case.  Our results should follow for arbitrary infinite cardinalities once Menger's and 3-decomposition theorems are proved.

Now, here are the main results, beginning with a simple regularity property.  

\begin{prop}\label{T:sepmax}
A regular biased expansion that is maximal is necessarily inseparable.
\end{prop}

\begin{proof}
Suppose a regular biased expansion $\Omega\downarrow\Delta$ has a cutpoint $v$, so that $\Omega = \Omega_1\cup\Omega_2$ and $\Omega_1\cap\Omega_2 = \{v\}$; let $\Delta_i = p(\Omega_i)$.  Choose $e_i \in E(\Delta_i)$ incident with $v$ and take any biased expansion $\Omega_3\downarrow K_3$ whose multiplicity equals that of $\Omega$.  Identify $e_1$ and $e_2$ with different edges of the $K_3$ and amalgamate edges to form, first, $\Omega_1 \cup_{e_1} \Omega_3$ and then $(\Omega_1 \cup_{e_1} \Omega_3) \cup_{e_2} \Omega_2$.  This is a proper extension of $\Omega$.  The disconnected case is similar.
\end{proof}

\begin{thm}[Structure of Maximal Biased Expansions]\label{T:structuremax}  
Any $2$-connected maximal biased expansion graph $\Omega \downarrow \Delta$ of finite order $n\geq 3$ is obtained by expanded edge amalgamation of group expansions of complete graphs of order at least $3$ and irreducible, nongroup circle expansions of order at least $3$, all of which are restriction subgraphs $\Omega\big|_{\Delta'}$ of $\Omega$.  
The group expansions and circle expansions are uniquely determined as the maximal complete subgraphs and the maximal chordless circle expansions contained in $\Omega$.

Any such edge amalgamation is a biased expansion.  It is maximal if and only if, for any two group expansions that are amalgamated along an expanded edge, the attachment map is twisted.
\end{thm}

The last part calls for explanation.  Twist was defined at Theorem \ref{T:twist}.  
Let $\Delta$ be the base graph of $\Omega$.  The theorem is saying, in part, that $\Delta$ is obtained by amalgamating circles and complete graphs.  In the second half, several complete graphs may be amalgamated along the same edge, either one at a time or all at once (to be explained momentarily).  Call these $\Delta_1,\ldots,\Delta_r$ and the common edge $e$, and let $\Omega_i = \Omega\big|_{\Delta_i}$.  
There are many ways to amalgamate one step at a time, each described by a rooted binary tree with leaves $\Omega_1,\ldots,\Omega_r$.  
We might amalgamate first, for instance, $\Omega_1$ and $\Omega_2$ by way of a bijection $\beta_{12} : p_1\inv(e)\to p_2\inv(e)$, then $\Omega_3$ and $\Omega_6$ by $\beta_{36}$, then $\Omega_5$ to $\Omega_1\cup_e \Omega_2$ via $\beta_{15}$, etc.  
All these ways have the same outcome, by Theorem \ref{T:amalgorder}.  
Instead, we could amalgamate all at once by means of commuting bijections $\beta_{ij} : p_i\inv(e) \to p_j\inv(e)$, that is, $\beta_{ik} = \beta_{ij} \circ \beta_{jk}$ and $\beta_{ij}\inv = \beta_{ji}$.  The theorem means that, if $\Omega_i = \fG_i\Delta_i$ and $\Omega_j = \fG_j\Delta_j$ for groups $\fG_i \cong \fG_j$, then $\beta_{ij}$ should not have the form that, according to Theorem \ref{T:twist}, makes $\Omega_i \cup_{e,\beta_{ij}}\Omega_j$ into a group expansion. 
(We discuss this further at Corollary \ref{C:basemax}.)

\begin{thm}[Structure of Biased Expansions]\label{T:structureexp}  
Any $2$-connected biased expansion $\Omega$ of a simple graph $\Delta$ of finite order at least $3$ is obtained by operations of expanded edge sum and amalgamation from $3$-connected group expansions and nongroup irreducible quasigroup expansions of circles, each of which is uniquely determined and is an expansion minor of $\Omega$.
\end{thm}

Call the group and circle expansions the \emph{3-constituents} of $\Omega$.  (In Theorem \ref{T:structureexp} they may not be uniquely determined.)  Note that $K_3$ is considered to be 3-connected.  
In the construction of $\Omega$ it may be that an edge $e$ in $\Delta$ belongs to several 3-constituents.  Then $p\inv(e)$ is the subject of several expanded amalgamations, and we could carry them all out at once as described previously.  Similarly, if an edge $e$ not in $\Delta$ belongs to several 3-constituents, then it is the subject of several edge sums; which means that all of the copies of $p\inv(e)$, except one, are amalgamated, and the last one is summed with the amalgamation of the others.  Then we could carry out, instead, a multiple (expanded) edge sum, similar to the multiple edge amalgamation we described.

For the proofs we need Tutte's theory of decomposition of an inseparable graph into 3-blocks.  We outline this theory (from \cite[Chapter IV, Sections 3 and 4]{Tbook}, originally in \cite{Tconn}).  Let $\Delta$ be a 2-connected graph.   If $\Delta$ is 3-connected, it is its own unique 3-block.  If it is not 3-connected, we define a \emph{cleavage} to be a 2-separation $\{x,y\}$ together with a bridge $B$ of $\{x,y\}$, such that $B$ is inseparable and not a single edge.  
(Then the complement of $B$ has at least two edges, since $\Delta$ is 2-connected and 2-separable.)  
Choose a cleavage, and split $\Delta$ into two graphs:  $B\cup e_{xy}$ and $B^c \cup e_{xy}$, where $B^c$ is the union of the other bridges of $\{x,y\}$ and $e_{xy}$ is a new edge, called a \emph{virtual edge}.  One continues this process on the resulting graphs until one obtains graphs $\Delta_1,\ldots,\Delta_k$ without cleavages.  These are the \emph{3-blocks} of $\Delta$.  Each virtual edge appears exactly twice and represents an edge sum; if all the indicated sums are carried out, the 3-blocks are reassembled into $\Delta$.  Each 3-block is either 3-connected, or a circle graph of order three or more, or a \emph{multilink} of size three or more (that is, a graph consisting of at least three parallel links and their two nodes).  
There is a graph of 3-blocks, in which the nodes are the 3-blocks and two 3-blocks are adjacent when they share a virtual edge.  
Tutte's theorem is, first, that the 3-blocks are uniquely determined by $\Delta$, and second, that the graph of 3-blocks is a tree, called the \emph{3-block tree} of $\Delta$.  

Suppose $\Delta$ is simple.  Then a multilink $\Delta_0$ contains at most one real edge (i.e., an edge of $\Delta$).  Suppose $\Delta_0$ does contain a real edge, $e$, and virtual edges $e_1,\ldots,e_k$.  
If $\Delta_1,\ldots,\Delta_k$ are the 3-blocks that contain the other copies of $e_1,\ldots,e_k$, then $\Delta_0\oplus_{e_1} \Delta_1 \oplus_{e_2} \cdots \oplus_{e_k} \Delta_k$ is the same as the amalgamation $\Delta_1\cup_e \Delta_2 \cup_e \cdots \cup_e\Delta_k$ if we treat all the $e_i$ as copies of $e$.  Thus, by amalgamating rather than summing we can dispense with $\Delta_0$.  
If $\Delta_0$ contains only virtual edges, then $\Delta_0 \oplus_{e_1} \Delta_1 \oplus_{e_2}\cdots \oplus_{e_k} \Delta_k$ is the same as $(\Delta_1\cup_{e_k}\cdots \cup_{e_k}\Delta_{k-1})\oplus_{e_k} \Delta_k$ if we treat all the $e_i$ as copies of $e_k$; so again we can dispense with $\Delta_0$.  (Or, again, we can treat this as a simultaneous edge sum.)  
The conclusion is that, for simple graphs $\Delta$, the multilinks are not needed if we modify the 3-blocks and permit amalgamation.  This is what we shall do.

For the proof of Theorem \ref{T:structuremax} we need the definition of 
a theta-complete graph from Section \ref{extension}.

\begin{lem}\label{L:thetacomplete}  
 Any theta-complete simple, 2-connected graph $\Delta$ is obtained by edge 
amalgamation of complete and circle subgraphs of $\Delta$, and conversely 
such an amalgamation is theta-complete.
 \end{lem}

\begin{proof}  
Consider the 3-blocks of $\Delta$ in Tutte's unmodified system. We show that every multilink 3-block $\Delta_0$ contains a real edge.  That is the same as saying that the two nodes of a cleavage are adjacent.  This comes from theta-completeness and a sublemma.

\begin{lem}\label{L:hingetheta}  
In Tutte's $3$-decomposition of any $2$-connected graph $\Delta$, two nodes $x$, $y$ of a cleavage are the trivalent nodes of a theta subgraph.  
\end{lem}

\begin{proof}  
If $\{x,y\}$ has more than two bridges, this is trivial.  If it has only two bridges, $B$ and $B^c$, then we know (by definition of a cleavage) that $B$ is 2-connected. Hence, $\Delta$ contains two internally disjoint $xy$-paths in $B$ and one more in $B^c$.
\end{proof}

Since every multilink 3-block does contain a real edge, it can be eliminated in favor of amalgamation.  And, because every 2-separating pair of nodes is adjacent, every virtual edge lies in a 3-block that is a multilink.  Consequently, when we modify Tutte's 3-decomposition all edge sums are replaced by amalgamations.

Conversely, we have to prove the amalgamation is theta-complete.  This is obvious.
\end{proof}

\begin{proof}[Proof of Theorem \ref{T:structuremax}]  
Assume $\Omega\downarrow \Delta$ is maximal.  From the results of the preceding section we know $\Delta$ is theta-complete.  The rest is obvious.

Conversely, suppose $\Omega\downarrow \Delta$ is the result of expanded edge amalgamations applied to group expansions $\fG_1 K_{n_1},\ldots,\fG_rK_{n_r}$ and nongroup irreducible circle expansions $\Omega_1\downarrow C_{l_1}, \ldots, \Omega_s\downarrow C_{l_s}$.  These are the 3-constituents of $\Omega$ and the $K_{n_i}$, $C_{l_j}$ are the 3-constituents of $\Delta$.  By Tutte's 3-decomposition theorem they are unique.  We have to prove $\Omega$ cannot be extended to any edge $e$ not in $\Delta$, the base graph constructed by the amalgamations.  

Suppose it did, for some $e\not\in E(\Delta)$, and let $\Omega'\downarrow \Delta\cup e$ be the extension.  The endpoints of $e$ cannot be contained within one 3-constituent, because each $K_{n_i}$ is complete, and if $\Omega$ extended to a chord of $C_{l_j}$, then $\Omega_j$ would be reducible (by Theorem \ref{L:q}).  It follows that, if we take a path in the 3-block tree of $\Delta$ joining a 3-block containing $x$ to a 3-block containing $y$, the path has positive length.  
Let $\Delta_1,\ldots,\Delta_r$ be the shortest such path, with $x$ in $\Delta_1$ and $y$ in $\Delta_r$, and set 
\begin{equation*}
\Delta'' = \Delta_1 \cup \cdots \cup \Delta_r \cup e.
\end{equation*}
Then $x$ and $y$ are connected by two internally disjoint paths in $\Delta'' \setm e$ and therefore by three in $\Delta''$.  

If $\Delta_1,\ldots,\Delta_r$ are all complete graphs, then $\Delta''$ is 3-connected, because the only 2-separations of $\Delta''$ are those at cleavages of $\Delta$ where a $\Delta_{h-1}$ and $\Delta_h$ share an edge.  But if $\Delta''$ is 3-connected, then $\Omega'\big|_{\Delta''}$ is a group expansion, and therefore $\Omega\big|_{\Delta_1}$ and $\Omega\big|_{\Delta_2}$ are group expansions, amalgamated by an attaching bijection that makes $\Omega\big|_{\Delta_1\cup \Delta_2}$ a group expansion, contrary to hypothesis.  So, some $\Delta_h$ is a circle of length $l\geq 4$.  

We may assume by choice of indices that $h > 1$, so that $\Delta_h$ amalgamates with $\Delta_{h-1}$ along an edge $uv$.  Also, either $y\in N(\Delta_h)$, or $h < r$ and $\Delta_h$ shares with $\Delta_{h+1}$ an edge $u'v'$.  It is easy to verify that one can name the nodes so that $u$ and $y$, in the former case, or $u$ and $u'$, in the latter, are not adjacent.  In the former case let $u' = y$.  Consider the two internally disjoint $xy$-paths in $\Delta''\setm e$.  One must pass through $u$ but not $v$; call $P_1$ its portion from $x$ to $u$.  One must pass through $u'$ but not $v'$; call $P_2$ its portion from $u'$ to $y$. (This is a trivial path if $u' = y$.)  $P_1$ and $P_2$ are internally disjoint from $\Delta_h$.  Consequently, the $uu'$-path $P_1\cup e\cup P_2$ is internally disjoint from $\Delta_h$, and in combination with the two $uu'$ paths in the circle $\Delta_h$, it forms a theta graph with trivalent nodes $u, u' \in N(\Delta \cup e)$.  By the previous section, then, $\Omega'$ extends to $\Delta\cup e \cup e_{uu'}$.  Because $u$ and $u'$ are not adjacent in $\Delta_h$, hence not in $\Delta$ either, we have contradicted the irreducibility of $\Delta_h$.  

Since in either case we deduce a contradiction, $\Omega$ is indeed maximal.
\end{proof}

As an example, the expanded edge amalgamation of two maximal biased expansions is maximal if (but not only if) for any group $\fG$ the two expansions contain at most one 3-constituent that is a $\fG$-expansion.

\begin{proof}[Proof of Theorem \ref{T:structureexp}]  
The trick is to extend $\Omega$ to edges $e_{xy}$ for all cleavages $\{x,y\}$.  We know from Lemma \ref{L:hingetheta} that this is possible, but we also need to know that the cleavages are the same in the extended base graph $\Delta'$.  Clearly, $\Delta'$ has all the cleavages of $\Delta$.  On the other hand, in a cleavage $(\{x,y\},B')$ of $\Delta'$, $\{x,y\}$ is a 2-separation of $\Delta$ and $B = B'\cap \Delta$ is connected, is a bridge of $\{x,y\}$, and has at least two edges; and the same holds for any other bridge $B'_1$ and $B_1= B'_1 \cap \Delta$ unless $B'_1$ is an edge.  These facts are a consequence of Lemma \ref{L:hingetheta}.  The conclusion is that $(\{x,y\},B)$ is a cleavage of $\Delta$.  That is, $\Delta$ and $\Delta'$ have the same cleavages.

Consequently, they have the same 3-blocks (in Tutte's sense) except that 
the 3-blocks in $\Delta'$ may have additional edges.  $\Omega'$ is obviously obtained from its 3-constituents by expanded edge amalgamation, and $\Omega$ is the same except for deletion of the amalgamated fibers $(p')\inv(e)$ for each additional edge $e$.  This deletion simply converts an amalgamation to a sum; thus $\Omega$ is obtained by edge sum and amalgamation from 3-connected group expansions and circle expansions.  Each circle expansion, if reducible, is an edge sum of smaller circle expansions; thus $\Omega$ does have the form stated in the theorem.

That all the 3-constituents are expansion minors follows from Theorem \ref{T:amalg}.
\end{proof}

Two questions remain.  First, what are the graphs that support maximal expansions?  Second, which graphs have nongroup expansions (a question raised in \cite[Example III.3.8]{BG}).  Theorems \ref{T:structuremax} and \ref{T:structureexp} suggest the answers, but there are details to attend to.  Let us call a complete graph {\it large} if it has at least four nodes.

\begin{cor}\label{C:basemax}  
 A finite simple graph $\Delta$ has a biased expansion that is maximal if 
and only if it is inseparable and is obtained by edge amalgamation of complete graphs and circles.  

Let $N_1$, resp.\ $N_0$, be the maximum number of large, resp.\ all, 
complete 3-constituents of $\Delta$ that contain any one edge.  The 
possible multiplicities of a maximal finite biased expansion 
$\gamma\cdot\Delta$ include every composite number $\gamma \geq 5$ such 
that $(\gamma-1)! \geq 2N_1$, as well as $\gamma=4$ if $N_0 \leq 3$.
 \end{cor}

\begin{proof}  
The form of $\Delta$ is entailed by Theorem \ref{T:structuremax}, but it is necessary to produce examples.  The general idea is to expand each 3-constituent $\Delta_i$ and amalgamate.  We assume $N \geq 2$.  Belousov and Sandik \cite{BelSand}, Frenkin \cite{Frenkin}, and Borisenko \cite{Borisenko} demonstrated the existence of an irreducible $k$-ary quasigroup with $\gamma$ elements for every $k \geq 3$ and composite $\gamma \geq 4$.  (See \cite{AG}.)  We also know there is a binary quasigroup of every order $\gamma \geq 5$ that is not isotopic to a group (by \cite[Theorem 1.5.1]{Latin} for $\gamma\neq 6$, \cite[Figure 1.3.1]{Latin} for $\gamma=6$).  As for a complete graph, it has group expansions of every multiplicity.  The difficulty is to assemble the expansions into a maximal expansion.

Consider some complete 3-constituents $\Delta_1,\hdots,\Delta_r$ that 
share an amalgamating edge $e$.  Expand them all by a group $\fG$ of order 
$\gamma$ to construct $\Omega_i =\bgr{\fG\Delta_i}$.  Now we need 
attachment maps $\beta_{ij} : p_j\inv(e)\to p_j\inv(e)$.  (We include 
$\beta_{ii} = \id$.)  Since we want the amalgamated expansion to be 
maximal, none of the $\bbeta_{ij}$ can be a pseudoautomorphism of $\fG$, 
except of course for the $\bbeta_{ii}$.  Factoring $\bbeta_{ij} = 
\bbeta_{1i}\inv \circ \bbeta_{1j}$, we conclude that the mappings 
$\bbeta_{1i}$ for $i = 1,2,\ldots,r$ must belong to different cosets of 
$\Aff\fG$, the group of pseuodautomorphisms, in the symmetric group of 
$\fG$.  This condition is necessary and sufficient for maximality of the 
amalgamation.

In the simplest case we expand every large complete 3-constituent on $e$ by the cyclic group $\bbZ_\gamma$.  The number of cosets of $\Aff \bbZ_\gamma$ is $(\gamma-1)!/2$, so we can accommodate $r\leq (\gamma-1)!/2$ different complete 3-constituents.  If $\gamma\geq5$ we expand the $K_3$ 3-constituents by binary nongroup isotopes so we can take $r=N$.  If $n\leq3$ we can expand every complete 3-constituent by $\bbZ_4$ and take $r=n$.  The corollary follows easily.
\end{proof}

The list of achievable multiplicities can be improved in special cases.  If all 3-constituents are complete, they can all be expanded by a group so $\gamma$ need not be composite; however, then we have to take $r=n$.  If $\Delta$ is a circle, $\gamma$ can be any composite number $\geq 4$.  If $\Delta$ is complete, $\gamma$ can be any positive integer.  In some situations we could handle larger $n$ or $N$ by using more than one gain group.

One would have liked to say that any two maximal biased expansions, 
$\gamma\cdot\Delta_1$ and $\gamma\cdot\Delta_2$, with a common base edge 
$e$ and the same multiplicity, can be amalgamated into a maximal expansion 
by choosing $\beta$ appropriately, but this is not true.  For one reason, 
there could be a group $\fG$ that is the gain group of several 
3-constituents, for which the combined number of 3-constituents in both 
graphs that are $\fG$-expansions and cover $e$ exceeds the number of 
cosets of $\Aff\fG$.  It is possible to describe the exact conditions 
under which an expanded amalgamation is maximal, in terms of double cosets 
of pseudoautomorphism groups of groups of order $\gamma$, but the 
description is excessively complicated.

\begin{cor}\label{C:basenongroup}
 A finite simple graph $\Delta$ has a regular biased expansion that is not 
a group expansion if and only if it is not a forest and is not 
3-connected.

The possible finite multiplicities of a regular nongroup expansion $\gamma\cdot\Delta$ include every $\gamma \geq 4$, except that when every block is 3-connected with at least four nodes $\gamma$ cannot be prime.
 \end{cor}

\begin{proof}  
If $\Delta$ is separable we can expand two different blocks by two different groups of order $\gamma$ with the exception noted.  In a 2-separable block we can expand every 3-constituent by $\bbZ_\gamma$ and make sure to attach one of them, whether by edge summation or edge amalgamation, so as to produce a nongroup expansion.
\end{proof}


\section{Four-node minors and thin expansions}\label{4minors}

A biased expansion graph may have gains for fairly special reasons.  
As we mentioned in connection with Corollary \ref{C:gminor}, gainability of minors of order four suffices to imply that $\Omega$ is a group expansion.  
Partially for that reason, a biased expansion may be forced to have gains in a group simply because its multiplicity is very small.

\begin{lem}\label{L:gminor4} 
If $\Omega\downarrow C_{n+1}$, where $n\geq 3$, and all expansion minors of order four that contain a specific edge fiber $p\inv(e_i)$ are group expansions (not necessarily of the same group), then $\Omega$ is a group expansion of $C_{n+1}$.
\end{lem}

\begin{proof}
 We assume the reader is acquainted with contraction of gain and biased 
graphs (see \cite[Sections I.2 and I.5]{BG}). 
 We write $C = C_{n+1} = e_0e_1 \cdots e_n$, with $N(e_i) = 
\{v_i,v_{i+1}\}$ where $v_0 = v_{n+1}$.  The special edge in the statement 
of the lemma will be $e_0$.  The case $n=3$ being trivial, we assume 
$n\geq 4$.

Fix a balanced lift $\tC^0$.  Some notation that will be convenient:  $\tC^0(\te_i,\te_j)$ is $\tC^0$ with $\te_i$ and $\te_j$ replacing $\te_i^0$ and $\te_j^0$.  $\Omega_{ijk}$ is the expansion minor of $\Omega$ whose edge set is $p\inv(\{e_i,e_j,e_k\})$ that is obtained by contracting $\tC^0 \setminus p\inv(\{e_i,e_j,e_k\})$; similarly, we write $\Omega_{0ijk}$, $\Omega_{ij}$.

The hypothesis is that each $\Omega_{0ijk} \cong \bgr{\fG_{ijk}C_4}$ for a group $\fG_{ijk}$.  $\Omega_{0ij} $ is an expansion minor of both $\Omega_{0ijk}$ and $\Omega_{0ijl}$.  In the former capacity it is isomorphic to $\bgr{\fG_{ijk}C_3}$ and in the latter to $\bgr{\fG_{ijl}C_3}$ (by contracting $\te^0_k$ and $\te_l^0$, respectively).  Since the gain group of a group expansion is unique, $\fG_{ijk} \cong \fG_{ijl}$. It follows that all groups $\fG_{ijk}$ are isomorphic to a single group $\fG$.  

In the rest of the proof we construct a gain graph $\Phi = \fG C$ and prove that $\bgr{\Phi} = \Omega$.  For this purpose we consider $e_i$ to be oriented from $v_{i}$ to $v_{i+1}$.

\emph{Step 1.}  We define the gain mapping $\phi$.  Its identity-gain edge set will be $\tC^0$.  Any isomorphism $\Omega_{0123} \cong \bgr{\fG C_4}$ defines gains $\phi$ on $p\inv(\{e_0,e_1,e_2,e_3\})$; we choose $\phi$ so it is $1$ on $\{\te^0_0,\te_1^0,\te_2^0,\te_3^0\}$.  We extend $\phi$ to $p\inv(e_i)$ for $i>3$ by $\phi(\te_i) = \phi(\te_0)\inv$, where $\te_0$ is the lift of $e_0$ that makes $\tC^0(\te_0,\te_i)$ balanced.  This rule can be expressed as choosing $\phi \big|_{p\inv(e_i)}$ so that $\bgr{\Phi_{0i}} = \Omega_{0i}$.

Note that, if we want to change $\phi\inv(1)$ to be a different balanced lift, $\tC^1$, we can do it by switching $\phi$.

\emph{Step 2.}  We next show that $\Phi$ is valid on expansion minors of order four that include $p\inv(e_0)$; that is, $\bgr{\Phi_{0ijk}} = \Omega_{0ijk}$.  

For $\Omega_{0123}$ that is a matter of definition.

For $\Omega_{012k}$ (where $k>3$), because $\Omega_{012k} \cong \bgr{\fG C_4}$ we can choose gains in $\fG$ for $\Omega_{012k}$, and we may choose them so that, contracted by $\te_3^0$ to $\Omega_{012}$, they agree with $\Phi_{012}$.  Then the gains on $\Omega_{012k}$ are forced by the $\Omega_{0k}$ minor to be as in $\Phi$.  Thus, $\bgr{\Phi_{012k}} = \Omega_{012k}$.  We infer that $\bgr{\Phi_{01k}} = \Omega_{01k}$.  

Considering $\Omega_{01jk} \cong \bgr{\fG C_4}$, the gains can be chosen to agree on $\Omega_{01j}$ with those of $\Phi_{01j}$.  The minor $\Omega_{0k}$ forces $\Omega_{01jk}$ to have gains as in $\Phi$, so $\bgr{\Phi_{01jk}} = \Omega_{01jk}$.  We further conclude from this and the previous cases that $\Omega_{ij} = \bgr{\Phi_{ij}}$ for every pair $\{i,j\}$, and that $\bgr{\Phi_{0jk}} = \Omega_{0jk}$.

Finally, $\Omega_{0ijk} \cong \bgr{\fG C_4}$ and the gains on $\Omega_{0ij}$ can be chosen to agree with those of $\Phi_{0ij}$.  Again $\Omega_{0k}$ forces the gains of $\Omega_{0ijk}$ to be as in $\Phi_{0ijk}$, so $\bgr{\Phi_{0ijk}} = \Omega_{0ijk}$.

\emph{Step 3.}  We prove by induction on $n$ that $\bgr{\Phi} = \Omega$.  The task is to prove that every lift $\tC^*$ is well behaved: it is balanced in $\Omega$ if and only if $\phi(\tC^*) = \phi(\te_0^*)\phi(\te_1^*)\cdots\phi(\te_n^*) = 1$.

If $\tC^*$ has an edge $\te^*_i = \te_i^0$ with $i\neq 0$, then we contract $\Omega$ and $\Phi$ by $\te^0_i$ and discard loops.  This gives expansion minors $\Omega' \downarrow C_n$ and $\Phi' = \fG C_n$, in which $\tC^*$ becomes $\tC^*/\te_j^0$ and $\phi'$ is the restriction of $\phi$.  The process of constructing gains in $\Omega'$ in Step 1 produces the gain function $\phi'$ if the various choices are made in agreement with those defining $\Phi$.  By induction, therefore, $\tC^*/\te_i^0$ is balanced in $\Omega'$ if and only if $\phi'(\tC^*/\te_i^0) = 1$.  However, $\phi'(\tC^*/ \te^0_i) = \phi(\tC^*)$ because $\phi(\te^0_i) = 1$, and by definition of contraction $\tC^*$ is balanced in $\Omega$ if and only if $\tC^*/\te^0_i$ is balanced in $\Omega'$.  Therefore, $\tC^*$ is well behaved.

If $\tC^*$ fails to contain an edge $\te^0_i$ with $i\neq 0$, we replace $\tC^0$ by a different balanced circle $\tC^1$ that does have an edge $\te_i^1$ in common with $\tC^*$. We choose $\tC^1 = \tC^0(\te_0^1,\te^1_1)$ where $\te_1^1 = \te^*_1$ and $\te^1_0$ is the edge that makes $\tC^1$ balanced; that is, $\phi(\te_0^1) = \phi(\te^1_1)\inv$, since $\tC^1$ is well behaved.  Changing $\tC^0$ to $\tC^1$ alters the gain mapping $\phi$, but under control: we simply switch it by a suitable switching function $\eta$.  A valid choice for $\eta$ is $\eta(v_i) = 1$ except $\eta(v_1) = \phi(\te^*_1)\inv = \phi(\te_0^1)$.  Then $(\phi^\eta)\inv(1) = \tC^1$, and because $\bgr{\Phi}$ is invariant under switching $\phi^\eta$ is a suitable gain function with respect to $\tC^1$ in Step 1.  By the previous case with $\tC^1$ in place of $\tC^0$, $\tC^*$ is well behaved.
\end{proof}

\begin{thm}\label{T:gminor4}
Suppose $\Omega$ is a finite, 2-connected biased expansion graph of order at least four. If every expansion minor of order four is a group expansion, then so is $\Omega$.
\end{thm}

\begin{proof}
The lemma demonstrates that all 3-constituents of $\Omega$ are group expansions.  If $\Omega $ is not a group expansion, then at some point in the process of amalgamation and summation two group expansions, $\fG_1\Delta_1$ and $\fG_2\Delta_2$, are summed (or amalgamated, which is treated similarly) along an edge $e$ by a twisted attachment map $\beta$ to form a nongroup biased expansion.  
There are expansion minors $\bgr{\fG_1C_3}$ of $\bgr{\fG_1 \Delta_1}$ and $\bgr{\fG_2 C_3}$ of $\bgr{\fG_2\Delta_2}$ that contain $p\inv(e)$, and $\bgr{\fG_1 C_3} \oplus_{e,\beta} \bgr{\fG_2 C_3}$ is a minor of $\bgr{\fG_1\Delta_1} \oplus_{e,\beta} \bgr{\fG_2 \Delta_2}$.  
It is not a group expansion of $C_3 \oplus_e C_3 = C_4$ because $\beta$ is twisted and twistedness is unaltered by taking minors.  
But $\bgr{\fG_1 C_3} \oplus_{e,\beta} \bgr{\fG_2 C_3}$ is one of the four-node expansion minors of $\Omega$ that, by hypothesis, are group expansions.  
This contradiction demonstrates that $\beta$ cannot be twisted.
\end{proof}

By Theorem \cite[V.2.1(a)]{BG}, in Theorem \ref{T:gminor4} we may assume just that every expansion minor of order four is gainable, or more simply that every minor of order four, all of whose edges are links, is gainable.

\begin{prob}\label{Pr:gminor4}
Can the list of order-four expansion minors in the hypotheses of Theorem \ref{T:gminor4} be reduced?  Can the list in Lemma \ref{L:gminor4} be reduced?
\end{prob}

We now examine the case of small multiplicity.

\begin{thm}\label{T:thin}  
Let $\Delta$ be a finite graph and $\Omega = \gamma\cdot\Delta$ a $\gamma$-fold biased expansion.  Then $\Omega = \bgr{\pm\Delta}$ if $\gamma=2$ and $\Omega = \bgr{\bbZ_3\Delta}$ if $\gamma=3$.  
\end{thm}

\begin{proof}[Proof of the Case $\gamma=2$.]  
Choose a balanced lift $\tDelta$.  Label $+$ any edge in $\tDelta$ and $-$ any edge in $E(\Omega) \setm E(\tDelta)$; this defines a signature $\sigma$ on $\Omega$.  We have to show that a circle $\tC$ in $\Omega$ is balanced if and only if it has an even number $f(\tC)$ of negative edges.  This follows from the labelling if $f(\tC) = 0$.  If $f(\tC) > 0$, let $\tC = \te_1 \cdots \te_l$ where $\te_l$ is negative and let $\te^*_l$ be the other edge that projects to $e_l$.  By the definition of a biased expansion, exactly one of $\tC$ and $\te_1\cdots \te_{l-1}\te^*_l$ is balanced.  Moreover, $f(\tC) = f(\tC^*)+1$.  Thus, by induction on $f(\tC)$, $f(\tC)$ is even $\iff f(\tC^*)$ is odd $\iff \tC^*$ is unbalanced $\iff \tC$ is balanced.
\end{proof}

\begin{proof}[Proof of the Case $\gamma=3$.]
For $3\cdot C_3$ the result follows from Proposition \ref{P:qgroup} and the fact that $\bbZ_3$ is, up to isotopy, the unique quasigroup of order 3.  

For $3\cdot C_4$ we define a gain function $\phi : E(\Omega) \to \bbZ_3$ and prove that $\bgr{\Phi} = \Omega$; for this purpose we consider $e_i$ to be oriented from $v_{i-1}$ to $v_i$.  We employ the notation of Lemma \ref{L:gminor4} except that the group is additive, with identity 0.

In the first step we define gains on $p\inv(\{e_0,e_1,e_2\})$ by means of an isomorphism $\Omega_{012} \cong \bgr{\bbZ_3 C_3}$.  We choose the gains so that $\tC^0 \cap p\inv(\{e_0,e_1,e_2\})$ has identity gains.  We define gains on $p\inv(e_3)$ so that $\bgr{\Phi_{03}} = \Omega_{03}$.

Consider $\Omega_{0i3}$ for $i = 1,2$.  It has gains in $\bbZ_3$.  Let $\phi'$ be these gains, chosen so that they are 0 on $\tC^0$ and agree with $\phi$ on $p\inv(e_0)$.  Then $\bgr{\Phi'_{0i}} = \Omega_{0i} = \bgr{\Phi_{0i}}$ implies that $\phi$ and $\phi'$ agree on $p\inv)e_i)$ and similarly they agree on $p\inv(e_3)$.  Therefore $\bgr{\Phi_{0i3}}=\Omega_{0i3}$.

Consequently, any lift $\tC$ that has in common with $\tC^0$ an edge $\te_i^0$ with $i\neq0$ is balanced in $\Omega$ if and only if it has gain 0 in $\Phi$.  One proves this by contracting $\te_i^0$ and by the fact that all $\bgr{\Phi_{0ij}}=\Omega_{0ij}$.

Let us identify $\tC = \te_0\te_1\te_2\te_3$ with its gain sequence
$(\phi(\te_0),\phi(\te_1),\phi(\te_2),\phi(\te_3))$.  Suppose $\tC$ has gain sequence $(a_0,a_1,a_2,a_3)$ that sums to 0, yet $\tC$ is unbalanced; we derive a contradiction.  We may assume $a_1,a_2,a_3 \neq 0$.  There are two cases, up to permutation of $e_1,e_2,e_3$ and negation of gains.

If the last three gains are not all equal, the gains are $(a_0,1,1,-1)$, hence $(-1,1,1,-1)$.  Since $(-1,0,1,-1)$ is unbalanced (because of $\bgr{\Phi_{023}}$), by the Circle Lifting Property $(-1,-1,1,-1)$ is balanced.  Then $(-1,-1,0,-1)$ is unbalanced; but this is impossible by $\bgr{\Phi_{013}}$.

The other case is that of gains $(a_0,1,1,1)$, i.e., $(0,1,1,1)$.  Here $(0,0,1,1)$ is unbalanced by $\bgr{\Phi_{023}}$, so $(0,-1,1,1)$ is balanced.  But then $(0,-1,0,1)$ is unbalanced, which is impossible due to $\bgr{\Phi_{013}}$.

We have shown that any circle whose gain is 0 is balanced.  It is then clear that a circle with nonzero gain is unbalanced.  Thus, $\bgr{\Phi}=\Omega$.

This solves the case $n=3$.  We conclude by Lemma \ref{L:gminor4} that $\Omega = 3\cdot C_{n+1}$ has gains in $\bbZ_3$ for all $n>3$.

It remains to solve the case in which $\Delta$ is 2-connected but not a circle.  Let $\Omega = 3\cdot \Delta$.  By the preceding case and Theorem \ref{T:structureexp}, $\Omega$ is the expanded edge amalgamation and sum of various $\bbZ_3$-expansions.  
By Corollary \ref{C:gamalg}, $\Omega$ is a $\bbZ_3$-expansion.
\end{proof}

There are two reasons why the theorem is limited to multiplicities below four.  
The simpler is that in each order $\gamma > 3$ there exists a (binary) quasigroup that is not isotopic to a group.  
The other is that, for many graphs, one can combine expansions by the same group of order at least four so as to make a nongroup expansion (Corollary \ref{C:basenongroup}).
Still, all counterexamples are separable since it is impossible to have a nongroup biased expansion, exept of $K_3$, that is 3-connected (Theorem \ref{T:3connbg}).


\section{Factorization and construction of multary quasigroups}\label{qgroup}

Let us discuss the consequences of our results for multary quasigroups.  
From a $k$-ary quasigroup $\fQ$ with operation $f$ (which we shall 
sometimes denote by $\fQ_f$) construct the \emph{factorization graph} $\Delta(\fQ)$:  this is the circle graph $C_{k+1}$, whose edges $e_i = e_{i-1,i}$ we call \emph{sides}, together with an added chord $e_{ij}$ whose endpoints are $v_i$ and $v_j$ whenever $f$ has a factorization 
\begin{equation}\label{E:fact}
f(x_1,\hdots,x_n) = g(x_1,\hdots,h(x_{i+1},\hdots,x_j),\hdots,x_k).
\end{equation}
 Clearly, $\Delta(\fQ) = K_{k+1}$ if $\fQ$ is isotopic to an iterated 
group, and the converse has long been known (see Lemma 
\ref{L:completebx}).  A stronger converse follows from Theorem 
\ref{T:3connbg}; that is Theorem \ref{T:3connq}. 
 From Theorem \ref{T:structuremax} we further deduce a structural 
description of multary quasigroups (Theorem \ref{T:quasidecomp}) due 
to Belousov.

To obtain our results we need the connection between the factorization graph and the maximal extension of $\bgr{\fQ C_{k+1}}$.

\begin{thm}\label{L:q}  
The maximal extension $\Omega(\fQ)$ of the biased graph $\bgr{\fQ C_{k+1}}$ corresponding to an $n$-ary quasigroup $\fQ$ is a biased expansion of the factorization graph $\Delta(\fQ)$.
\end{thm}

\begin{proof}  By the theorems of Section \ref{extension} it suffices to prove that $\bgr{\fQ C_{k+1}}$ extends to every chord in $\Delta(\fQ)$ but to no other chord of $C_{k+1}$; i.e., the last part of Proposition \ref{P:qgroup}.  

Suppose $\bgr{\fQ C_{k+1}}$ extends to a chord $e_{ij}$.  Call the extension $\Omega$.  Let $C'$ and $C''$ be the circles formed by the chord, with $e_0 \in C'$.  Then $\Omega' = \Omega\big|_{C'}$ and $\Omega'' = \Omega\big|_{C''}$ define operations $g$ and $h$ satisfying \eqref{E:fact} by the construction described in Section \ref{expan}.  Thus, $e_{ij}$ belongs to $\Delta(\fQ)$. 

Suppose on the other hand that $f$ factors as in \eqref{E:fact}.  Then $\bgr{\fQ_g C'} \cup_{e_{ij},\beta} \bgr{\fQ_h C''}$, which we call $\Omega$, is a biased expansion of $C_{k+1} \cup e_{ij}$, where we take the amalgamating mapping $\beta : {p'}\inv(e_{ij}) \to {p''}\inv(e_{ij})$ to be the identity function $\beta (xe_{ij}) = xe_{ij}$. A circle $\{x_0e_0,x_1e_1,\hdots,x_ke_k\}$ is balanced in $\Omega$ if there is an edge $xe_{ij}$ that makes $\{xe_{ij},x_{i+1}e_{i+1},\hdots,x_je_j\}$ and \newline $\{xe_{ij},x_0e_0,x_1e_1,\hdots,x_ie_i,x_{j+1}e_{j+1},\hdots,x_ke_k\}$ both balanced.  In terms of $g$ and $h$, this means that 
\begin{equation*}
x = h(x_{i+1},\hdots,x_j)
\end{equation*}
and
\begin{equation*}
x_0 = g(x_1,\hdots,x_i,x,x_{j+1},\hdots,x_k).
\end{equation*}
It follows that $x_0 = f(x_1,\hdots,x_k)$, so $\Omega\big|_{C_{k+1}} = \bgr{\fQ C_{k+1}}$.  Thus, $\bgr{\fQ C_{k+1}}$ extends to every chord $e_{ij}$ in $\Delta(\fQ)$.
\end{proof}

We see in the proof of Theorem \ref{L:q} that expanded edge amalgamation is the analog of functional composition.  

We immediately obtain from Theorem \ref{T:3connbg} the promised strong characterization of iterated group isotopes.

\begin{thm}\label{T:3connq} 
If $\fQ$ is a $k$-ary quasigroup with $k\geq 3$ and $\Delta(\fQ)$ is $3$-connected, then $\fQ$ is isotopic to an iterated group.
\hfill\qedsymbol
\end{thm}

Therefore, if $\Delta(\fQ)$ is 3-connected it is complete.  We mentioned at Lemma \ref{L:completebx} the long-known fact that completeness of $\Delta(\fQ)$ implies that $\fQ$ is an iterated group isotope.  The new result amounts to saying that one need not know $\Delta(\fQ)$ completely to arrive at the same conclusion.  

\begin{exam}\label{X:leftright}
Suppose $2n-2$ binary quasigroups satisfy the identity
$$
f_{n-1}(f_{n-2}(\cdots(f_2(f_1(x_1,x_2),x_3),\ldots,x_{n-1}),x_n) = g_1(x_1,g_2(x_2,\ldots,g_{n-1}(x_{n-1},x_n)\cdots)) .
$$
We see immediately that the $n$-ary operation defined by either side of this equation has 3-connected factorization graph.  Therefore, all $f_i$ and $g_i$ are isotopic to one group.
\end{exam}

\begin{exam}[Multary groups]\label{X:multgroup}
For instance, consider a $k$-ary group (with $k\geq 3$), where Equation \eqref{E:assoc} holds.  Let $\hat \fQ$ be the $(2k-1)$-ary quasigroup with operation $\hat f$ defined by \eqref{E:assoc}.  Multary associativity means that $\Delta(\hat \fQ)$ contains diameters $e_{i,i+k}$ for $i = 0,1,\hdots,k-1$, so it is 3-connected.  By Theorem \ref{T:3connq}, $\hat \fQ$ is an iterated group isotope.  It follows that $\fQ$ is an iterated group isotope, either by an easy algebraic argument or by combinatorial reasoning: $\bgr{\fQ C_{k+1}}$ is a subgraph of $\bgr{\hat \fQ C_{2k}}$ extended to the chords; the latter is a group expansion; therefore the former is a group expansion; therefore $\fQ$ is isotopic to an iterated group.  This is the easy part of the theorem of Hossz\'u and Gluskin mentioned in the introduction; their complete result is much stronger and is not an immediate corollary of our work.

\begin{figure}[htbp]\label{F:multarygroup}
\begin{center}
\psfrag{v0}[c]{$v_{0}$}
\psfrag{v1}[c]{$v_{1}$}
\psfrag{vk-1}[c]{$v_{k-1}$}
\psfrag{vk}[c]{$v_{k}$}
\psfrag{v2k-2}[c]{$v_{2k-2}$}
\psfrag{v2k-1}[c]{$v_{2k-1}$}
\psfrag{e0}[c]{$e_{0}$}
\psfrag{e1}[c]{$e_{1}$}
\psfrag{ek-2}[c]{$e_{k-2}$}
\psfrag{ek-1}[c]{$e_{k-1}$}
\psfrag{e2k-1}[c]{$e_{2k-1}$}
\psfrag{e}[c]{$e_{k-1,2k-1}$}
\includegraphics{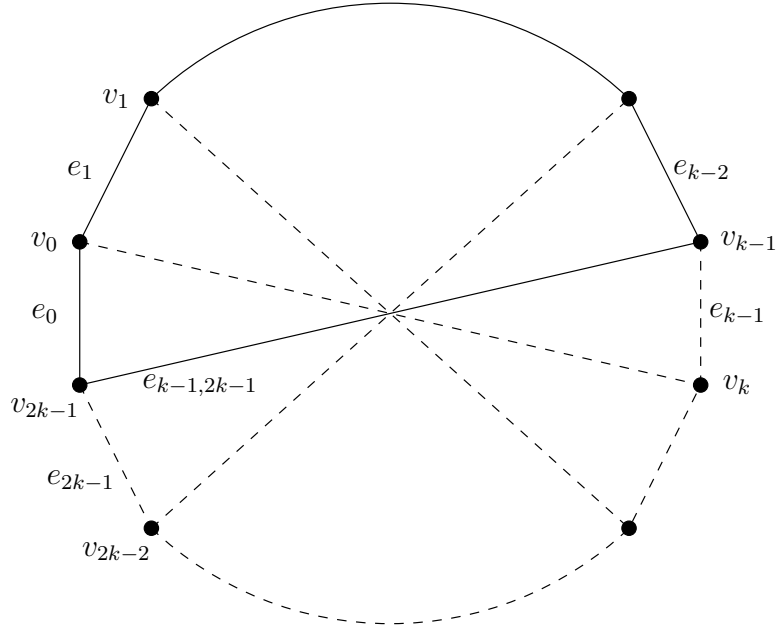}
\end{center}
\caption{The sides and diametric chords of $\Delta(\hat \fQ)$, showing the embedded $C_{k+1}$ corresponding to $\fQ$ (solid lines).}
\end{figure}

The various appearances of $f$ need not represent the same operation; it is only necessary that for $k$-ary quasigroups $f_1, g_1,\ldots,f_k, g_k$ the $k$ compositions
$$
f_i(x_1,\ldots,g_i(x_{i},\ldots,x_{i+k-1}),\ldots,x_{2k-1})
$$
be independent of $i$.  Then all $2k$ operations are isotopic to the $k-1$-fold iteration of a single group operation; this is part of a theorem of U\v{s}an \cite{Usan}.
\end{exam}

The basis for Theorem \ref{T:3connq} is that the factorization graph of 
any multary quasigroup is theta-complete.  (This is the quasigroup version 
of Proposition \ref{CX:thetacompletion}.)  We are able to characterize 
factorization graphs completely.

\begin{thm}\label{L:qtheta}
For a simple graph $\Delta$ to be the factorization graph of a multary quasigroup, a necessary and sufficient condition is that $\Delta$ be theta-complete and have a Hamiltonian circle.  A second necessary and sufficient condition is that $\Delta$ be obtained by edge amalgamations of circles and complete graphs and have a Hamiltonian circle.
\end{thm}

\begin{proof}
Apply Theorem \ref{T:structuremax} in view of Theorem \ref{L:q}.
\end{proof}

The amalgamation in Theorem \ref{T:structuremax} corresponds to a decomposition of $f$ into iterated group isotopes, irreducible multary quasigroups of arity greater than 2, and nongroup binary quasigroups.  Furthermore, the decomposition of $f$ is unique because the 3-constituents of $\Omega(\fQ)$ are unique.  Thus we have:

\begin{cor} [Belousov; see {\cite[Section V.4]{BelnK}}]\label{T:quasidecomp}  
Every multary quasigroup is in a unique way (up to isotopy) the composition of iterated group isotopes and irreducible, nongroup multary quasigroups.
\end{cor}

Belousov deduces this through the algebra of multary quasigroup composition.  His key result about such composition is our next corollary, which we prove by another application of theta completeness.  Bear in mind that a 1-ary quasigroup is merely a permutation of the set $\fQ$.

\begin{cor} [{\cite[Theorem 2.1]{BIAQ}, \cite[Chapter IV]{BelnK}}]\label{C:qfactor}
Suppose an $n$-ary quasigroup $\fQ$ has an $(i,k)$ factorization,
$$
f(x_1,\ldots,x_n) = g(x_1,\ldots,x_{i-1},h(x_i,\ldots,x_k),\ldots,x_n),
$$
and a $(j,l)$ factorization, 
$$
f(x_1,\ldots,x_n) = g'(x_1,\ldots,x_{j-1},h'(x_j,\ldots,x_l),\ldots,x_n),
$$
where $i \leq k$ and $j \leq l$ and $g,h,g',h'$ are multary quasigroups.

(a)  If $i \leq j \leq k \leq l$, then
$$
f(x_1,\ldots,x_n) = g(x_1,\ldots, a(x_i,\ldots,h'(x_j,\ldots,x_l),\ldots,x_k), \ldots,x_n),
$$
where $a$ is a multary quasigroup.

(b)  If $k\leq k \leq j \leq l$, then
$$
f(x_1,\ldots,x_n) = b(x_1,\ldots, h(x_i,\ldots,x_k), x_{k+1},\ldots, h'(x_j,\ldots,x_l), \ldots,x_n),
$$
where $b$ is a multary quasigroup.

(c) If $i<j\leq k<l$, then
$$
f(x_1,\ldots,x_n) = c(x_1,\ldots,d(x_i,\ldots,x_{j-1})\circ d'(x_j,\ldots,x_k)\circ d''(x_{k+1},\ldots,x_l),\ldots,x_n),
$$
where $c,d,d',d''$ are multary quasigroups and $\circ$ is a group multiplication.
\end{cor}

\begin{proof} 
(a) and (b) are immediate from Theorem \ref{L:q}.  (c) is from that theorem, Theorem \ref{L:qtheta}, and the case of order four in Lemma \ref{L:completebx}.
\end{proof}

Other of our results on biased expansions also have quasigroup interpretations.
Corollary \ref{C:gamalg} applied to multary quasigroups is the following statement:

\begin{cor}\label{C:gqgroup}
A composition of multary quasigroups, all isotopic to iterates of a group $\fG$, is necessarily isotopic to an iterated group if $\fG = \bbZ_\gamma$ for $\gamma \leq 3$ or $\fG = \fV_4$, but not otherwise; and then it is isotopic to an iteration of $\fG$. 
\end{cor}

The quasigroup version of Corollary \ref{C:gminor} requires a 
definition.  Take a $k$-ary quasigroup $\fQ$.  Lemma 
\ref{L:bexpminor}(a) implies that expansion minors of $\bgr{\fQ 
C_{k+1}}$ of order $r+1$ correspond to residual $r$-ary quasigroups.  
Apply Construction XM, taking $S \cup T = C_{k+1}$.  The choice of lift 
$\tT$ signifies fixing the values of the variables corresponding to edges 
of $T$.  The variables of the residual quasigroup are the variables that 
correspond to edges of $S$.

\begin{cor}\label{C:qminor}
It is possible to have an $n$-ary quasigroup of any order $\gamma \geq 4$ and any arity $n\geq3$ that is not isotopic to an iterated group but whose residual binary quasigroups are all isotopic to the same arbitrary group of order $\gamma$, except when the group is $\fV_4$.
\end{cor}

But raising the residual arity yields quite a different result.  The 
quasigroup interpretation of Lemma \ref{L:gminor4} is: 

\begin{thm}\label{T:qinduced3}
 If every residual ternary quasigroup of an $n$-ary quasigroup $\fQ$ with 
arity $n\geq3$ is isotopic to an iterated group (not necessarily the same 
group), then $\fQ$ is an iterated group isotope.
 \end{thm}

The interpretation of Theorem \ref{T:thin} for multary quasigroups (from 
the case $\gamma\cdot C_n$) is a result that was stated by Belousov, who 
published a proof only for order $\gamma=2$ because of the length of the 
proof for order three (this information obtained by Dudek \cite{Dpc}; I 
have not been able to find references).

\begin{thm}\label{C:thin}
A multary quasigroup of order three or less is isotopic to an iterated group.
\end{thm}


\section{Postscript}\label{ps}

\subsection{Nontopological homotopy?}\label{homotopy}

There is a perceptible flavor of homotopy about our combinatorial arguments.  We treat balanced circles in a manner reminiscent of contractible circles in a topological space.  
A way of making this similarity exact is to embed the underlying graph in a topological space so that the balanced circles are precisely the graph circles that are contractible.  That is possible if and only if the graph has gains \cite{UTG} so it cannot be used to justify our reasoning.  
Nevertheless the analogy is suggestive.  One has to wonder what lies behind it.

\subsection{Formulas and bijections}\label{bijections}

In characterizing multary groups (Example \ref{X:multgroup}) our method yields a description up to isotopy, and this is typical of our results.  Hossz\'u and Gluskin, however, found an exact formula for any multary group operation.  
Belousov \cite[Chapter 3]{BIAQ} strenthened this by completely charactering $(i,j)$-associative operations: those that satisfy the hypothesis of Corollary \ref{C:qfactor} with $g=h=g'=h'$.  
I believe their formulae and some of the many generalizations can be reproduced and perhaps further extended if the expansion-graph method is supplemented by careful attention to the exact bijections between the set $\fQ$ and the edge fibers.


\section*{Acknowledgements}

I am grateful to Wies{\l}aw Dudek, Vladislav Goldberg, and Victor Shcherbakov for 
assistance in locating and interpreting previous work on multary quasigroups and groups. I 
am indebted to Jeff Kahn for suggesting the design interpretation of Section \ref{td} 
and to Lori Koban for advice on the definition of homomorphism of biased graphs.  I thank 
Marge Pratt for speedy and dependable typing assistance.


\end{document}